\documentclass[dvips,12pt]{article} 
\usepackage{amsmath,amssymb,amsfonts,amscd,rotating} %
\def\simarrow{\mathrel{\raise -0.5mm\hbox{$\sim$}}\hspace{-1.8mm}{\rightarrow} } 

\def\bsimarrow{\leftarrow\hspace{-0.7mm}\mathrel{\raise -0.5mm\hbox{$\backsim$}} }

\def\bt{\begin{tabular}}
\def\te{\end{tabular}}

\def\lettrine#1#2#3{\noindent\hangindent#1\hangafter-#2
\hskip-#1\smash{\hbox to #1{#3\hfill}}\ignorespaces}

\newcommand{\To}[1]{\mathop{\to}\limits_{#1}}

\def\BM{\begin{pmatrix}}
\def\EM{\end{pmatrix}}

\def\txt{\textstyle}
\def\ds{\displaystyle}

\def\d=f{\buildrel\hbox{\scriptsize d\'{e}f}\over \Longleftrightarrow}

\def\cit{\text{\it I\hskip -6ptC\/}}

\def\square{\hfill\hbox{\vrule height .9ex width .8ex depth -.1ex}}
\def\rit{\text{\it I\hskip -2pt  R}}

\def\qit{\text{\it I\hskip -5.5pt  Q}}

\def\nit{\text{\it I\hskip -2pt  N}}

\def\rl {\rit^{\hskip 1pt\ell}}

\def\As{{\cal A}}

\def\Bd{{\text B}}

\def\Ed{{\text E}}
\def\Es{{\cal E}}

\def\Hs{{\cal H}}

\def\Ms{{\cal M}}

\def\Os{{\cal O}}

\def\Ws{{\cal W}}

\def\ung{\hbox{1\hskip -4.2pt \rm 1}}

\def\be{\begin{equation}}
\def\ee{\end{equation}}
\def\beqn{\begin{eqnarray}}
\def\eeqn{\end{eqnarray}}
\def\nobeqn{\begin{eqnarray*}}
\def\noeeqn{\end{eqnarray*}}
\def\ba{\left(\begin{array}}
\def\ea{\end{array} \right) }

\def\bpr{\paragraph{Proof.}}

\def\epr{\square\vskip 6pt}
\def\eop{\hbox{\vrule height .9ex width .8ex depth -.1ex}}

\def\o{\overline}
\def\and{\; \mbox{and} \;}

\newcommand{\half}{\frac{1}{2}}

\def\hfl#1#2{\smash{\mathop{\hbox to 12mm{\rightarrowfill}}
\limits^{\scriptstyle #1}_{\scriptstyle #2}}}

\def\Ker{\mathop{\rm Ker}\nolimits}

\def\mod{\mathop{\rm mod}\nolimits}

\def\Be{\begin{enumerate}}
\def\Ee{\end{enumerate}}

\def\Bena{\begin{enumerate}
\def\labelenumi{\theenumi)}
\def\theenumi{\arabic{enumi}}
\def\labelenumii{\theenumii)}
\def\theenumii{\alph{enumii}}}

\def\Bean{\begin{enumerate}
\def\labelenumii{\theenumii)}
\def\theenumii{\arabic{enumii}}
\def\labelenumi{\theenumi)}
\def\theenumi{\alph{enumi}}}

\def\Bero{\begin{enumerate}
\def\labelenumii{\theenumii)}
\def\theenumii{\arabic{enumii}}
\def\labelenumi{(\theenumi)}
\def\theenumi{\roman{enumi}}}

\def\BeRo{\begin{enumerate}
\def\labelenumii{\theenumii)}
\def\theenumii{\arabic{enumii}}
\def\labelenumi{(\theenumi)}
\def\theenumi{\Roman{enumi}}}

\def\Bi{\vskip 11pt\begin{itemize}\itemsep=18pt}
\def\bi{\begin{itemize}}
\def\Ei{\end{itemize}\vskip 11pt}
\def\ei{\end{itemize}}

\def\Bd{\begin{description}}
\def\Ed{\end{description}}

\def\R{\right}
\def\L{\left}
\def\F{\frac}

\usepackage[latin1]{inputenc} 

\newcommand{\Aa}{\mathbb{A}\,}
\newcommand{\NN}{\mathbb{N}\,}
\newcommand{\ZZ}{\mathbb{Z}\,}

\newcommand{\CC}{\mathbb{C}\,}
\newcommand{\RR}{\mathbb{R}\,}
\renewcommand{\cit}{\mathbb{C}\,}
\renewcommand{\rit}{\mathbb{R}\,}
\newcommand{\FF}{\mathbb{F}\,}
\newcommand{\WW}{\mathrm{W}\,}
\newcommand{\QQ}{\mathbb{Q}\,}
\renewcommand{\qit}{\mathbb{Q}\,}

\def\sum{\mathop{\Sigma}\limits}
\def\prod{\mathop{\Pi}\limits}
\def\bigoplus{\mathop{\oplus}\limits}

\def\bbf{\boldmath\bf}

 \def\To{\longrightarrow}
\def\RL{_{R,L}}
\def\SRL{_{S_{R,L}}}

\def\RTL{_{R\times L}}
\def\SRTL{_{S_{R\times L}}}

\def\LR{_{L,R}}
\def\3{\perp\joinrel\perp\joinrel\perp}
 \def\lr{left (resp. right) }
\def\rl{right (resp. left) }
 \newcommand{\gfrak}{\mathfrak{g}}
\def\tr{\operatorname{trace}}
\def\det{\operatorname{det}}
\def\Gal{\operatorname{Gal}}
\def\GL{\operatorname{GL}}
\def\SL{\operatorname{SL}}
\def\Frob{\operatorname{Frob}}
\def\End{\operatorname{End}}
\def\mod{\operatorname{mod}}
\def\Im{\operatorname{Im}}
\def\Re{\operatorname{Re}}
\def\Hom{\operatorname{Hom}}

\def\groth{\operatorname{Groth}}
\def\Cusp{\operatorname{cusp}}
\def\Irr{\operatorname{Irr}}
\def\Lie{\operatorname{Lie}}
\def\Rep{\operatorname{Rep}}
\def\Repsp{\operatorname{Rep\ sp}}
\def\EIS{\operatorname{EIS}}
\def\eis{\operatorname{eis}}
\def\ELLIP{\operatorname{ELLIP}}

\begin{document} \thispagestyle{empty}   \null\vfill \begin{center} {\LARGE From global class field concepts and modular representations  to the conjectures of Shimura-Taniyama-Weil, Birch-Swinnerton-Dyer and Riemann  } \vfill {\sc C. Pierre\/} \vskip 11pt  Institut de Mathématique pure et appliquée\\ Université de Louvain\\ Chemin du Cyclotron, 2\\ B-1348 Louvain-la-Neuve,  Belgium  \\pierre@math.ucl.ac.be   \vfill Mathematics subject classification (1991): 11G07, 11F11, 11F30. \eject \end{center}  \null\vfill \begin{center} {\LARGE From global class field concepts and modular representations  to the conjectures of Shimura-Taniyama-Weil, Birch-Swinnerton-Dyer and Riemann  } \vfill {\sc C. Pierre\/} \vfill

\begin{abstract} Based upon new global class field concepts leading to two-dimensional global Langlands correspondences, a modular representation of cusp forms is proposed in terms of global elliptic bisemimodules which are (truncated) Fourier series over $\RR$~.  As application, the conjectures of Shimura-Taniyama-Weil, Birch-Swinnerton-Dyer and Riemann are analyzed. \end{abstract} \thispagestyle{empty} \vfill \eject \end{center}

\setcounter{page}{1} \section*{Introduction} 
The main objective of this paper consists in providing a modular representation of cusp 
forms in terms of global elliptic semimodules which are (truncated) Fourier series over $\RR$ whose terms correspond
to the places of the considered (semi)fields. \vskip 11pt 

This challenge, directly connected to global  Langlands correspondences \cite{Lan}, depends on the following central result:  

At every place of a (CM) (semi)field corresponds a conjugacy class of the Borel subgroup of the general (linear)group.  This conjugacy class is associated with a sublattice of Hecke and is in one-to-one correspondence with the corresponding class of representation of the global Weil-Deligne group defined in this paper in such a way that the Hecke eigenvalues and the Frobenius eigenvalues coincide.  

 More concretely, the symmetric splitting semifields $ \widetilde F_L$ and $\widetilde F_R$ are introduced as symmetric algebraic extensions of a global number field $k$ of characteristic zero. In the complex case, $\widetilde F_L$ and $\widetilde F_R$ are respectively composed of a set of complex and conjugate complex simple roots of a polynomial ring over $k$~. In the real case, the left and right symmetric semifields are noted respectively $\widetilde F_L^+$ and $\widetilde F^+_R$~.\vskip 11pt

The infinite places of these semifields are characterized by the algebraic extension degrees of their completions. For example, the $n$-th place of $F_L^+$ is defined with respect to the algebraic extension degree
\[ [\widetilde F^+_{v_n}:k] \simeq f_{v_n}\cdot N=n\cdot N\] 
of the pseudo-ramified extension $\widetilde F^+_{v_n}$~, associated with the completion $F^+_{v_n}$ 
where $f_{v_n}=n$ is its global class residue degree and where $N$ is the degree of an irreducible algebraic extension and also
the order of the global
inertia subgroup $I_{F_{v_n}^+}$~. This global inertia subgroup can be viewed as the subgroup of inner automorphisms of Galois with respect to the Galois subgroup $\Gal(\widetilde F^+_{v_n}/k)$ which will be considered as a subgroup of modular automorphisms of Galois of the algebraic general linear subsemigroup $T_2(\widetilde F^+_{v_n})$~.\vskip 11pt

Let $T_2(F^+_{v})$ be the algebraic group of upper triangular matrices over the set of completions $F^+_v$ of $F^+_L$ at the set of real places $v=\{V_1,\cdots,V_n,\cdots,V_s\}$~. \vskip 11pt 

We are interested in the enveloping (semi)algebra $B_R\otimes B_L\simeq T^t_2(F^+_{\o v}) \times T_2(F^+_{v})$ of the division (semi)algebra $B_L\simeq T_2(F^+_v)$ because $T^t_2(F^+_{\o v})\times T_2(F^+_v)$~, denoted $\GL_2(F^+_{\o v}\times F^+_v)$~, is an algebraic bilinear semigroup ``having a representation" in the tensor product $M_{F^+_{\o v}}\otimes M_{F^+_v}$ of a right $T^t_2(F^+_{\o v})$-semimodule $M_{F^+_{\o v}}$ by a left $T_2(F^+_{v})$-semimodule $M_{F^+_{v}}$~.\vskip 11pt

Indeed, we show that a linear algebraic group $\GL_2(F^+)$ over a symmetric field $F^+=F^+_R\cup F^+_L$~, having as	representation space a $2^2$-dimensional vector space, is covered under the conditions of proposition 1.11 by the bilinear algebraic semigroup $\GL_2(F^+_R\times F^+_L)$~.
\vskip 11pt

The $\GL_2(F^+_{\o v}\times F^+_v)$-bisemimodule $M_{F^+_{\o v}}\otimes M_{F^+_v}$ decomposes into the direct sum of 
$\GL_2(F^+_{\o v_n}\times F^+_{v_n})$-subbisemimodules $M_{F^+_{\o v_{n,m_n}}}\otimes M_{F^+_{v_{n,m_n}}}$~, having $m_n$ representatives, according to the decomposition of a Hecke bilattice $\Lambda ^{(2)}_{\o v}\otimes \Lambda ^{(2)}_{v}$ into subbilattices $\Lambda ^{(2)}_{\o v_{n,m_n}}
 \otimes \Lambda ^{(2)}_{v_{n,m_n}}$ at the biplaces $\o v_n\times v_n$ of $F^+_{\o v}\times F^+_v$~. And, the ring of endomorphisms of the $\GL_2(F^+_{\o v}\times F^+_v)$-bisemimodule $M_{F^+_{\o v}}\otimes M_{F^+_v}$~, decomposing it into the set of subbisemimodules following the subbilattices, is generated by the products $(T_{q_R}\otimes T_{q_L})$ of Hecke operators $T_{q_R}$ and $ T_{q_L}$ for $q \nmid N$ and by the products $( U_{q_R} \otimes U_{q_L})$ of $U_{q_R}$ and $U_{q_L}$ for $q\mid N$~: it is noted $T_H(N)_R\otimes T_H(N)_L$~.\vskip 11pt

The coset representative of $U_{q_L}$ (resp. $U_{q_R}$~), referring to the upper (resp. lower) half plane, is chosen to be upper (resp. lower) triangular and given by the integral matrix $\left( \begin{smallmatrix} 1 & b_N \\ 0&q_N\end{smallmatrix}\right)$ (resp. $\L( \begin{smallmatrix} 1 & 0\\ b_N & q_N\end{smallmatrix}\R)$~) where the $b_N$ and $q_N$ are integers modulo $N$~.\vskip 11pt

So, $U_{q_R}\otimes U_{q_L}$ has the coset representative \[ g_2(q^2_N,b^2_N) = \L[ \BM 1 & b_N \\ 0&1\EM\BM 1&0\\ b_N & 1\EM\R]\ \BM 1&0\\ 0&q_N^2\EM\] where \[ D_{q^2_N,b^2_N}=u_2(b_N)\times u^t_2(b_N)=\BM 1&b_N \\ 0&1\EM\BM 1&0\\ b_N&1\EM\] is interpreted as the decomposition group element associated with the split Cartan subgroup element $\alpha_{q^2_N}=\L(\begin{smallmatrix} 1&0\\ 0&q^2_N\end{smallmatrix}\R)$~.\vskip 11pt

As a consequence, we have the following proposition:\vskip 11pt

{\em There is an explicit irreducible semisimple (pseudo-)ramified representation $\rho _{\lambda _{\pm}}$~, associated with a weight two cusp form, \[ \rho_{\lambda_{\pm}}:\quad \Gal(\widetilde F^+_{\o v}/k) \times \Gal(\widetilde F^+_{v}/k) \To
\GL_2(T_H(N)_R\otimes T_H(N)_L)\]
having eigenvalues:
\[ \lambda_{\pm}(q^2_N,b^2_N) =
\F{(1+b^2_N+q^2_N)\pm [(1+b^2_N+q^2_N)^2-4q^2_N]^{\half}}2\]
verifying:
\Bi
\item \quad $\ds\tr\rho_{\lambda_{\pm}}=1+b^2_N+q^2_N$~;
\item \quad $\det \rho_{\lambda_{\pm}}= \lambda_+(q^2_N,b^2_N) \cdot \lambda_-(q^2_N,b^2_N) 
=q^2_N$~.
\Ei}\vskip 11pt

Let us introduce the (pseudo-)ramified lattice bisemispace
\[ X_{S_{R\times L}}=\GL_2( \widetilde F_{R }\times \widetilde F_L)\big/ \GL_2(\ZZ\big/N\ \ZZ)^2\]
where
\Bi
\item the general bilinear semigroup $\GL_2(\widetilde F_{R}\times \widetilde F_L)$ is taken over the
product
$ \widetilde F_{R}\times \widetilde F_L$ of the complex number semifields $\widetilde F_R$ 
and $\widetilde F_L$~;
\item $\GL_2((\ZZ\big/N\ \ZZ)^2)$ has entries in squares of integers modulo $N$~.
\Ei
The toroidal compactification of $X_{S_{R\times L}}$ of the Borel-Serre type can be considered as a toroidal projective 
isomorphism
\[ \gamma^c_{R\times L}:\quad X_{S_{R\times L}}\To \o X_{S_{R\times L}}\;, \]
with
\[ \o X_{S_{R\times L}} =\GL_2(F^T_{R }\times F^T_L )\big/\GL_2(\ZZ\big/N\ \ZZ)^2
\approx \GL_2(F^T_{\o\omega }\times F^T_\omega )\]
where
\Bi
\item $F^T_L $ is the toroidal compactified equivalent of $\widetilde F_L $~;
\item $F^T_\omega $ denotes the set of toroidal complex completions.
\Ei

Its boundary is given by:
\[ \partial \o X_{S_{R\times L}} =\GL_2(F^{+,T}_{R}\times F^{+,T}_L )\big/\GL_2(\ZZ\big/N\ \ZZ)^2
\approx \GL_2(F^{+,T}_{\o v }\times F^{+,T}_v )\]
where \Bi
\item $F^{+,T}_L$ is the  toroidal compactified equivalent of $\widetilde F^+_L$~;
\item $F^{+,T}_v$ denotes the set of toroidal real completions;
\Ei
in such a way that
\Be
\item the number of complex places is equal to the number of real places;
\item the multiplicity of the complex places is equal to one;
\item the complex conjugacy class representatives are covered by the multiples of their real equivalents.
\Ee
\vskip 11pt

The cosets of $\o X_{S_{R\times L}}$ (resp. $\partial \o X_{S_{R\times L}} $~) correspond to the 
conjugacy classes of\linebreak $\GL_2(F^T_{\o \omega }\times F^T_\omega )$ (resp. $\GL_2(F^{+,T}_{\o v }\times F^{+,T}_v )$~).\vskip 11pt

The analytic representation of the algebraic bilinear semigroup $ \GL_2(F^T_{\o \omega }\times F^T_\omega )$ is given, by means of a Langlands global correspondence, by the product (of the Fourier development) of  a  right cusp form of weight $k=2$ restricted to the lower half plane by its classical left analogue restricted to
the upper half plane.\vskip 11pt

Similarly, in the real case, a Langlands global correspondence relates the algebraic 
bilinear semigroup $\GL_2(F^{+,T}_{\o v }\times F^{+,T}_v )$ to its analytic representation given by 
the product, right by left, of two (truncated) Fourier series over $\RR$ called global elliptic $A_R$ (resp. $A_L$~)-semimodules.\vskip 11pt

These Langlands global correspondences are based upon the representation spaces of the algebraic
bilinear semigroups given in terms of their conjugacy class representatives characterized by increasing ranks. The double sets of conjugacy class representatives of the envisaged algebraic bilinear semigroups generate two symmetric towers in such a way that the analytic representation of each conjugacy class representative is given by a term of the Fourier development of the cusp form entering into the cuspidal representation of the considered algebraic bilinear semigroup. This allows to give an algebraic-geometric interpretation to the Fourier series.
\vskip 11pt

Let
\begin{align*}
\phi _L (s_L) &= \sum_n \sum_b \phi (s_L)_{n,b}\ q^n\big/ \QQ_L\;, \quad n\le \infty\;, 
\quad q=e^{2\pi ix}\;, \quad x\in\rit\;,\\[11pt]
\text{(resp.} \quad
\phi _R (s_R) &= \sum_n \sum_b \phi (s_R)_{n,b}\ q^{-n}\big/ \QQ_R\;)\end{align*}
be respectively a \lr global elliptic $A_L$ (resp. $A_R$~)-semimodule where:
\Bi
\item $A_L$ (resp. $A_R$~) is the ring of sections $s_L$ (resp. $s_R$~) of the semisheaf of rings over the $T_2(F ^{+,T}_v)$-semimodule 
$M_{F^{+,T}_v}$ (resp. $T^t_2(F ^{+,T}_{\o v})$~)-semimodule $M_{F^{+,T}_{\o v}}$~);
\item $\phi (s_L)_{n,b}\equiv\lambda(n^2_N,b^2_N)$ such that $\lambda(n^2_N,b^2_N)$ is an 
eigenvalue of the coset representative $g_2(n^2_N,b^2_N)$ of product of Hecke operators.
\vskip 11pt
\Ei

Then, the inclusion of $\GL_2(F^{+,T}_{\o v }\times F^{+,T}_v )$ into 
$\GL_2(F^{T}_{\o \omega }\times F^{T}_\omega )$ in the sense of the Borel-Serre compactification 
implies the following results: \vskip 11pt

Let $S_L(f)$ (resp. $S_R(f)$~) be the left [semi-]algebra (resp. right [semi-]algebra) of modular forms
\[ f_L(z)=\sum_n a_{n,L}\ q^n_L\qquad \text{(resp.} \quad 
f_R(z)=\sum_n a_{n,R}\ q^n_R\;)\]
\[ q_L=e^{2\pi iz} \qquad \txt{(resp.} \quad q_R=e^{-2\pi iz}\;),\quad z\in\cit\;,\]
\Bi
\item being normalized eigenforms of Hecke operators related to the congruence subgroup 
$\Gamma _L(N)$ (resp. $\Gamma _R(N)$~) introduced in Section 2.12,
\item characterized by a weight two and a level $N$~.
\Ei\vskip 11pt

On the other hand, let $S_L(\phi )$ (resp. $S_R(\phi )$~) denote the left [semi-]algebra (resp. right [semi-]algebra) 
of global elliptic $A_L$-semimodules $\phi _L(s_L)$ (resp. $A_R$-semimodules $\phi _R(s_R)$~) 
in the sense that $f_L(z)\simeq \phi _L(s_L)$(resp. $f_R(z)\simeq \phi _R(s_R)$~).

Then, we have the following inclusions of left [semi-]algebras (resp. right [semi-]algebras):
\begin{align*}
S_L (\phi ) \quad &\hookrightarrow \quad S_L(f)\\
\text{(resp.} \quad S_R (\phi ) \quad &\hookrightarrow \quad S_R(f)\;).\end{align*}
And, if the bisemialgebras of modular forms and of global elliptic semimodules are introduced respectively
by $S_R(f)\otimes S_L(f)$ and by $S_R(\phi )\otimes S_L(\phi )$~, then the inclusion of these 
bisemialgebras can be stated by: 
\[ S_R (\phi )\otimes S_R (\phi ) \quad \hookrightarrow \quad S_R(f)\otimes S_L(f)\;.\]
Remark that a global elliptic $A_L$-semimodule $\phi _L(s_L)$ (resp. $A_R$-semimodule $\phi _R(s_R)$~) constitutes an automorphic representation of a modular form $f_L(z)$ (resp. $f_R(z)$~) in  the sense that the modular representation of $f_L(z)$ (resp. $f_R(z)$~) can be given by a set of $n$~, $1\le n\le \infty $~, two-dimensional semitori $T^2_L[n]$ (resp. $T^2_R[n]$~), restricted to the upper (resp. lower) half plane and covered each one by $m_n$ semicircles of the ``~$n$-th class'' of the global elliptic semimodule $\phi _L(s_L)$ (resp. $\phi _R(s_R)$~).
\vskip 11pt

This kind of modular representation is used to analyze the conjectures of Shimura-Taniyama-Weil, Birch-Swinnerton-Dyer and Riemann.
\vskip 11pt

The hyperbolic uniformization of arithmetic type of an elliptic curve $E(\QQ)$ is studied
 by envisaging its modular representation by means of the surjective mapping 
$\Hs^{\rm cusp\ (res)}_{\GL^{\rm res}_{2_\RR}\to E(\QQ)}$ of the restricted global elliptic $A_{R-L}$-bisemimodule $\phi _R(s_R)_{\rm res}\otimes_D \phi _L(s_L)_{\rm res}$ into $E(\QQ)$~, where $\otimes_D$ denotes a ``diagonal'' tensor product (see section 2.20).
\vskip 11pt

The modular representation of $E(\QQ)$ given by
\[ \Hs^{\rm cusp\ (res)}_{\GL^{\rm res}_{2_\RR}\to E(\QQ)}: \quad
\phi _R(s_R)_{\rm res}\otimes_D \phi _L(s_L)_{\rm res}  \to E(\QQ)\]
can be worked out from the $p$ sets of surjective mappings:
\[\{E_f(p_N,m_p)_R\otimes E_f(p_N,m_p)_L\}_{m_p}\to E(\FF_p)\]
where $E_f(p_N,m_p)$ is a semicircle of rank $p_N=p\cdot N$ entering into the covering of the semitorus $T^2[p]$~, for all prime $p$ taken into account in the restricted eulerian product of $L_R(s_-,E(\QQ)\otimes_DL_L(s_+,E(\QQ))$ in such a way that the orbit space of 
$\{E_f(p_N,m_p)_R\otimes E_f(p_N,m_p)_L\}_{m_p}$ is associated with the elliptic curve $E(\FF_p)$~.
\vskip 11pt

This hyperbolic uniformization of arithmetic type of the elliptic curve $E(\QQ)$ then corresponds to the  Shimura-Taniyama-Weil conjecture and is related to the problem of Diophantine equations by means of the Mordell-Weil group of $E(\QQ)$~.
\vskip 11pt

The Birch-Swinnerton-Dyer conjecture is analyzed in the same context.  Indeed, it consists in the fact that, if $L_R(s_-,E(\QQ))$ and $L_L(s_+,E(\QQ))$ are the $L$-subseries attached to an elliptic curve $E(\QQ)$ and introduced in chapter 3, then the ``pseudo-unramified'' rank of $E(\QQ)$ is the order of vanishing of these $L$-subseries at $s=1$~.  The trivial zeros and non-trivial zeros of these $L$-subseries, defined with respect to the set of $p$ primes envisaged above, are evaluated with respect to the Riemann conjecture which can be treated as follows:

Taking into account that the trivial zeros of the classical zeta function are equal at a sign to 
the global class residue degrees multiplied by a factor 2 and that a one-to-one correspondence 
must exist between trivial zeros and pairs of non-trivial zeros, we are led to formulate the 
proposition: 
 Let $D_{4n^2,i^2}\cdot \varepsilon _{4n^2}$ be a coset representative of the Lie algebra of the 
decomposition group $D_{i^2}(\ZZ)$ and let $\alpha _{4n^2}$ be the corresponding split 
Cartan subgroup element.  Then, the products of the pairs of the trivial zeros of the Riemann zeta functions 
$\zeta _R(s_-)$ and $\zeta _L(s_+)$ are mapped into the products of the corresponding pairs of the 
non-trivial zeros following: \begin{align*} 
D_{4n^2,i^2}\cdot \varepsilon _{4n^2} : \quad 
\det (\alpha _{4n^2}) \quad
&\To \quad \det (D_{4n^2,i^2}\cdot \varepsilon_{4n^2}\cdot \alpha _{4n^2})_{ss}\;,\\
\{(-2n)\cdot (-2n)\} \quad
&\To \quad \{ \lambda_+(4n^2,i^2,E_{4n^2})\cdot \lambda_-(4n^2,i^2,E_{4n^2})\}\;,
\quad \forall\ n\in\nit\;,\end{align*}
where $\varepsilon _{4n^2}=\L(\begin{smallmatrix} E_{4n^2} & 0 \\ 0&1\end{smallmatrix}\R)$ is the 
infinitesimal generator of the Lie algebra 
$(\gfrak\ell_2(F^{+,T,nr}_{\o v}\times F^{+,T,nr}_v)$~.
\vskip 11pt

I would like to thank especially P. Cartier, B. Mazur and D. Zagier for helpful comments.
\vskip 11pt

\section{From global class field concepts to modular representations of general bilinear semigroups}  

\subsection{Preliminaries on semiobjects}

The developments of this paper will essentially concern symmetric objects in such a way that the considered mathematical objects can be cut into two symmetric semiobjects $\Os _R$ and $\Os _L$~.  The left semiobject $\Os _L$ will be localized in (or will refer to) the upper half space while the corresponding right semiobject $\Os _R$ will be localized in (or will refer to) the lower half space.

The right semiobject $\Os _R$ is then the dual of the left semiobject $\Os_L$ and the interest of considering a symmetric object ``~$\Os $~'', decomposed into two dual semiobjects $\Os _R$ and $\Os _L$~, is that the informations concerning the internal mathematical structure of  ``~$\Os $~'' can be obtained from the product $\Os _R\times \Os _L$ of the semiobjects $\Os _R$ and $\Os _L$~.  Indeed, every endomorphism $E$ of the object ``~$\Os $~'' can be decomposed into the product $E_R\times E_L$ of a right endomorphism $E_R$ acting on the right semiobject $\Os _R$ by the opposite left endomorphism $E_L$ acting on the corresponding left semiobject $\Os _L$ such that $E_R=E_L^{-1}$~.
\vskip 11pt 

The existence of symmetric objects can be established by the following considerations on function fields.

Let $k$ be a global number field of characteristic zero.  Let $k[x_1,\cdots,x_r]$ be a polynomial ring over $k$ 
and let $R=R_R\cup R_L$ be a symmetric finite  extension of $k$ as introduced in the following.

\Bi
\item Let {\bbf $I_L=\{P_\mu (x_1,\cdots,x_r)\mid P_\mu (V_L)=0\}$ be the ideal of $k[x_1,\cdots,x_r]$\/} in such a way that:
\Bean
\item $P_\mu (V_L)$ be the polynomial function in $k[V_L]$ represented by $P_\mu (x_1,\cdots,x_r)$~;
\item $V_L\subset R_L$ be an affine semispace restricted to the upper half space.
\Ee
\vskip 11pt

\item Let {\bbf $I_R=\{P_\mu (-x_1,\cdots,-x_r)\mid P_\mu (V_R)=0\}$ be the symmetric ideal of $I_L$\/} obtained under the involution
\[\tau : \quad P_\mu (x_1,\cdots,x_r)\To P_\mu (-x_1,\cdots,-x_r)\]
in such a way that:
\Bean
\item $P_\mu (V_R)$ be the polynomial function in $k[V_R]$
 represented by $P_\mu (-x_1,\cdots,-x_r)$~;
 \item $V_R\subset R_R$ be an affine semispace restricted to the lower half space, symmetric of $V_L$ and disjoint of $V_L$ or possibly connected to $V_L$ on a symmetry axis plane.
 \Ee
 \vskip 11pt
 
 \item {\bf The quotient ring\/} obtained modulo the ideal $I_L$ (resp. $I_R$~) is $Q_L=k[x_1,\cdots,x_r]\big/ I_L$ (resp. $Q_R=k[x_1,\cdots,x_r]\big/ I_R$~) \cite{Wat}.
 
 $Q_L$ and $Q_R$ are quotient algebras characterized by the corresponding homomorphisms:
 \[ \phi _L: \quad Q_L\To R_L\;, \qquad \phi _R: \quad Q_R\To R_R\]
 where $R_L$ and $R_R$ are commutative (division) semirings localized respectively (or referring to) the upper and lower half spaces (commutative semirings are recalled (or introduced) in section 1.3).

 So {\bbf the pair of homomorphisms $\phi _L$ and $\phi _R$\/} sends the general solution to a pair of symmetric solutions respectively in $R_L$ and in $R_R$.
 \vskip 11pt
 
 \item On the other hand, let $\widetilde F= \widetilde F_R\cup \widetilde F_L$ be a symmetric ``algebraic'' finite extension of $k$~.
 
 Let $T_r(\widetilde F_L)$ (resp. $T^t _r(\widetilde F_R)$~) $\subset \GL_n(\centerdot)$ be the group of matrices of dimension $r$ over $\widetilde F_L$ (resp. $\widetilde F_R$~) viewed as an operator sending $\widetilde F_L$ (resp. $\widetilde F_R$~) into the affine semispace $T^{(r)}(\widetilde F_L)$ (resp. $T^{(r)}(\widetilde F_R)$~) of dimension $r$~:
 \[ T_r(\centerdot) : \quad \widetilde F_L\To T^{(r)}(\widetilde F_L) \qquad
 \text{(resp.} \quad 
 T^t_r(\centerdot) : \quad \widetilde F_R\To T^{(r)}(\widetilde F_R) \ )\]
 in such a way that to the indeterminates $(x_1,\cdots,x_\ell,\cdots,x_r,\cdots,x_{1\ell},\cdots,x_{1r})$ of $Q_L$ (resp. 
 $(-x_1,\cdots,-x_\ell,\cdots,-x_r,\cdots,-x_{1\ell},\cdots,-x_{1r})$ of $Q_R$~), $\forall\ x_{\ell r}=x_{\ell}\times x_r$~, corresponds the homomorphism:
\begin{align*}
\phi' _L: \quad Q_L &\To \widetilde F_L(x_1\to e_{11},\cdots,x_\ell=e_{\ell\ell},\cdots, x_{r\ell}\to e_{r\ell})\\
\text{(resp.} \quad 
\phi' _R: \quad Q_R &\To \widetilde F_R(-x_1\to -e_{11},\cdots,-x_\ell=-e_{\ell\ell},\cdots, -x_{r\ell}\to e_{r\ell}\ )\end{align*}
where
\begin{align*}
T_r(\widetilde F_L) &= \{e_L=(e_{\ell r})\in T_r(\widetilde F_L)\mid P_{T_\mu }(e_{\ell r})=0\}\\
\text{(resp.} \quad 
T^t_r(\widetilde F_R) &= \{e_R=(e_{\ell r})\in T^t_r(\widetilde F_R)\mid P_{T_\mu }(e_{\ell r})=0\}\ )\end{align*}
with the polynomials $P_{T_\mu }(e_{\ell r})\in k[x]$ (resp. $P_{T_\mu }(e_{\ell r})\in k[x]$~).
\vskip 11pt

\item Finally, let $X_L$ (resp. $X_R$~) be the functor from the quotient ring $Q_L$ (resp. $Q_R$~) to the affine semispace $V_L$ (resp. $V_R$~) in such a way that the diagram:
\[\begin{CD}
Q_L @>{\phi' _L}>> \widetilde F_L(e_{11},\cdots,e_{rn})\\
@VV{X_L}V @VV{T_r(\centerdot)}V\\
V_L @>{\psi _L}>\sim> T^{(r)}\widetilde F_L)\end{CD}\qquad \left(\text{resp.} \quad
\begin{CD}
Q_R @>{\phi' _R}>> \widetilde F_R(-e_{11},\cdots,e_{rn})\\
@VV{X_r}V @VV{T^t_r(\centerdot)}V\\
V_R @>{\psi _R}>\sim> T^{(r)}\widetilde F_R)\end{CD}\ \right)\]
commutes.

Such functor $X_L$ (resp. $X_R$~) from the $k$-algebra $Q_L$ (resp. $Q_R$~) to the affine semispace $V_L$ (resp. $V_R$~), homeomorphic to $T^{(r)}(\widetilde F_L)$ (resp. $T^{(r)}(\widetilde F_R)$~), is then representable and is a \lr affine semigroup scheme over $k$~.
\Ei
\vskip 11pt

The left and right affine semigroup schemes $X_L$ and $X_R$ are said to be symmetric if every element $a_L\in T^{(r)}(\widetilde F_L)$~, localized in the upper half space, is symmetric to every element $a_R\in T^{(r)}(\widetilde F_R)$~, localized in the lower half space, with respect to an axis or a plane\ldots
\vskip 11pt 

So, the consideration of a sufficiently large polynomial ring as $k[x_1,\cdots,x_n]$ allows to envisage generally objects as being symmetric and being able to be cut into two symmetric semi-objects introduced in a general way in the following section.
\vskip 11pt

\subsection{Semistructures}

Having shown as a general rule the existence of symmetric objects, we shall now recall or introduce the semistructures which will be used in the following developments.

The condensed notation $R,L$ means ``right, resp. left'' (and $L,R$ means ``left, resp. right'').

\Bi
\item A \rl {\bf semigroup\/} $G\RL$ is a nonempty set of \rl elements, localized in (or referring to) the lower (resp. upper) half space, together with a binary operation on $G\RL$~, i.e. a function $G\RL\times G\RL\to G\RL$ or $G_{L,L}\times G\RL\to G\RL$~.

\item A \rl {\bf monoid\/} is a \rl semigroup $G\RL$ which contains an identity element $a\RL\in G\RL$ such that:
\[a_R\cdot e_R=a_R \quad \text{(resp.}\quad e_L\cdot a_L=a_L\ ),\qquad \forall\ a\RL\in G\RL\;.\]

\item A \rl {\bf semiring\/} is a nonempty set $R\RL$ together with two binary operations (addition and multiplication) such that:
\Bean
\item $(R\RL,+)$ is an abelian \rl semigroup.
\item $(a\RL\cdot b\RL)\ c\RL=a\RL\ (b\RL\cdot c\RL)$~, \quad $\forall\ a\RL,b\RL,c\RL\in R\RL$ (associative multiplication).
\item $a\RL\ (b\RL+c\RL)=a\RL\ b\RL+a\RL\ c\RL$ and $(a\RL+b\RL)\ c\RL=a\RL\ c\RL+b\RL\ c\RL$ (left and right distribution).
\Ee

\item If $R\RL$ is a commutative semiring with identity $\ung_{R\RL}$ and no zero divisors, it will be called a \rl integral domain.

Furthermore, if every element of $R\RL$ is a unit (right and left invertible), $R\RL$ is a {\bf division semiring\/}.

And, a \rl {\bf semifield\/} is a commutative division semiring.

\item A \rl {\bf adele semiring\/} is the product of the primary completions of the \rl semifield.

\item Let $R\RL$ be a \rl semiring.  A \rl {\boldmath \bf $R\RL$-semimodule\/} is an additive abelian \rl semigroup $M\RL$ together with a function $M_R\times R\RL\to M_R$ (resp. $R\LR\times M_L\to M_L$~) such that:
\Bean
\item $(a_R+b_R)\ r\RL=a_R\ r\RL+b_R\ r\RL$ where $a_R\ r\RL=(a\cdot r)_R\in M_R$ (resp. $r\LR\ (a_L+b_L)=r\LR\ a_L+r\LR\ b_L$ where $r\LR\ a_L=(r\cdot a)_L\in M_L$~), \quad $\forall\ r\RL\in R\RL$~, $a\RL,b\RL\in M\RL$~.

\item $a_R\ (r\RL+s\RL)=a_R\ r\RL+a_R\ s\RL$ (resp. $(r\LR+s\LR)\ a_L=r\LR\ a_L+s\LR\ a_L$~), \quad $\forall\ s\RL\in R\RL$~.

\item $(a_R\ s_R)\ r_R=a_R\ (s_R\ r_R)$ (resp.  $r_L\ (s_L\ a_L)=(r_L\ s_L)\ a_L$~).
\Ee

If $R\RL$ has an identity element $\ung\RL$ such that $a_R\ \ung_R=a_R$ (resp. $\ung_L\ a_L=a_L$~), $M\RL$ is a \rl unitary $R\RL$ semimodule.

If $R\RL$ is a \rl division  semiring, then the unitary \rl $R\RL$-semimodule is a \rl {\bf vector semispace\/}.

\item If $R\RL$ is a commutative semiring with identity, a {\boldmath\bf $R\RL$-semialgebra $\As\RL$} is a semiring $\As\RL$ such that:
\Bean
\item $(\As\RL,+)$ is a unitary \rl $R\RL$-semimodule.
\item $(a_R\ b_R)\ r\RL=a_R\ (b_R\ r\RL)=b_R\ (a_R\ r\RL)=(a\ b\ r)_R\in \As_R$
(resp. $r\LR\ (a_L\ b_L)=(r\LR\ a_L)\ b_L=a_L\ (r\LR\ \ b_L)= (a\ b\ r)_L\in \As_L$~).
\Ee

If $\As\RL$ is a division semiring, then $\As\RL$ is called a {\bf division semialgebra\/}.
\Ei
\vskip 11pt 

The generation of the global algebraic extension (semi)fields and of the corresponding completions considered in this paper presents some analogy with the construction of local $p$-adic number fields which will be recalled in the following section.
\vskip 11pt 

\subsection{Classical notions about local fields}  

The field $K$~, which is a finite extension of $\QQ_p$~, is a $p$-adic field.  Let $\Os_K$ denote its ring of integers, $\wp_K$ the unique maximal ideal of $\Os_K$~, $k(\upsilon_K)=k(\wp_K)=\Os_K/\wp_K$ its residue field and $\widetilde \omega_K$ a uniformiser in $\Os_K$~.  Let $\upsilon_K:K^*\to \ZZ$ be the unique valuation so that the absolute value on $K$ is defined by $|\ |_K=|\ |_{\upsilon_K}$ with $|x|_K=(\#k(\upsilon_K))^{-\upsilon_K(x)}$ for $x\in K^*$~.
\vskip 11pt
The number of elements in $k(\wp_K)$ is $q=p^f$ where $f_{\upsilon_K}=[k(\upsilon_K):\FF_p]$ is the residue degree over $\QQ_p$~.  The ideal $\wp_K\Os_K$ of $\Os_K$ has the form $\wp_K^{e_{\upsilon_K}}=\tilde\omega_K^{e_{v_K}}$ where $e_{\upsilon_K}$ is the ramification degree of $K$ over $\QQ_p$~.  

Then, we have $[K:\QQ_P]=e_{\upsilon_K}\cdot f_{\upsilon_K}$ so that $e_{\upsilon_K}=[K:\QQ_P]/f_{\upsilon_K}$~.
\vskip 11pt
The maximal unramified extension of $K$ is denoted $K^{nr}$ and its completion is $\widehat K^{nr}$~.  
The inertia subgroup $I_k$ is such that:
\[  \Gal(K^{ac}/K)/I_K\;\overset{\sim}{\longrightarrow}\;\Gal(K^{nr}/K)\;\overset{\sim}{\longrightarrow}\;\Gal (k(\upsilon_K)^{ac}/k(\upsilon_K))\;.\]
\vskip 11pt

The local Weil group $W_K\subset \Gal (K^{ac}/K)$ is the inverse image of $\Frob^{\ZZ}_{k(\upsilon_K)}\linebreak \subset \Gal (k(\upsilon_K)^{ac}/k(\upsilon_K))$ in $\Gal(K^{ac}/K)$~.
\vskip 11pt 

\subsection{Infinite places of a global number field}  

Let $k$ denote a global number field of characteristic 0 and let $k[x]$ 
be a polynomial ring composed of a family of pairs of polynomials $\{P_{T_\mu} (x),P_{T_\mu} (-x)\}$~, $1<\mu <\infty $~.  Then, according to section 1.1, the algebraic extension $\widetilde F$ of $k$~, assumed to be generally closed, is the
 symmetric splitting field $\widetilde F=\widetilde F_R\cup \widetilde F_L$ composed of a right extension semifield $\widetilde F_R$ 
and of a left extension semifield $\widetilde F_L$ in one-to-one correspondence.  If the algebraic extension field of 
$k$ is complex, then $\widetilde F_L$ (resp. $\widetilde F_R$~) is composed of the set of complex (resp. conjugate complex) 
simple roots of the polynomial ring $k[x]$~, while, if the algebraic extension field is real, then it will 
be noted $\widetilde F^+=\widetilde F^+_R\cup \widetilde F^+_L$ where $\widetilde F_L^+$ (resp. $\widetilde F_R^+$~) is the \lr algebraic extension semifield composed 
of the set of positive (resp. symmetric negative) simple roots. 
 \vskip 11pt

The left and right equivalence classes of the local completions of $\widetilde F^{(+)}_L$ and of $\widetilde F^{(+)}_R$ are the left and right 
real (resp.  complex) infinite places of $\widetilde F^{(+)}_L$ and of $\widetilde F^{(+)}_R$~: they are equal in number and noted, in the real case:
\[ v=\{ v_1,\cdots,v_n,\cdots,v_s\} \quad \text{and} \quad 
\o v=\{\o v_1,\cdots,\o v_n,\cdots,\o v_s\}\]  and, in the complex case:
\[\omega =\{\omega _1,\cdots,\omega _n,\cdots,\omega _s\} \quad \text{and} \quad 
\o \omega =\{\o \omega _1,\cdots,\o \omega _n,\cdots,\o\omega _s\}\;\]

With reference to the classical $p$-adic treatment of local fields, ``pseudo-unramified'' infinite real places will be 
characterized algebraically by their global class residue degrees given by:
\[ [F^{+,nr}_{v_n}:k]=f_{v_n}=n \quad \text{and} \quad
[F^{+,nr}_{\o v_n}:k]=f_{\o v_n}=n \;, \quad n\in\nit\;, \]
where $F^{+,nr}_{v_n}$ (resp. $F^{+,nr}_{\o v_n}$~), denoting a \lr pseudo-unramified local completion of $\widetilde F^+_L$ 
(resp. $\widetilde F^+_R$~), is a \lr $k$-semimodule whose dimension is given by the global class residue degree $f_{v_n}=f_{\o v_n}=n$~.
\vskip 11pt

The integer $n$ is the order of the archimedean valuation $v$ (resp. $\o v$~) defined as a homomorphism of $F^{+,nr}_{v_n}$ (resp. $F^{+,nr}_{\o v_n}$~) into the group of positive real numbers in such a way that, for all $x_1\in F^{+,nr}_{v_1}$ (resp. $\o x_1\in F^{+,nr}_{\o v_1}$~) and $x_n\in F^{+,nr}_{v_n}$ (resp. $\o x_n\in F^{+,nr}_{\o v_n}$~), we have:
\begin{equation}
v(n\ x_1)\ge v(x_n) \qquad \text{(resp.} \quad 
\o v(n\ \o x_1)\ge \o v(\o x_n) \ )\tag*{\cite{Koc}}\end{equation}
\vskip 11pt

Similarly, (infinite) complex places will be defined by their global class residue degrees:
\[ [F^{nr}_{\omega _n}:k] = f_{\omega _n}=n
\quad \text{and} \quad
[F^{nr}_{\o\omega _n}:k] = f_{\o\omega _n}=n\;, \quad n\in\nit\;, \]
where $F^{nr}_{\omega _n}$ (resp. $F^{nr}_{\o\omega _n}$~), denoting a \lr pseudo-unramified local completion of $\widetilde F_L$ 
(resp. $\widetilde F_R$~), is a \lr $k$-semimodule having as dimension the global class residue degree $f_{\omega _n}=f_{\o \omega _n}=n$~.
\vskip 11pt

The corresponding pseudo-ramified completions are assumed to be generated from irreducible $k$-semimodules 
$F^+_{v_n^1}$ (resp. $F^+_{\o v_n^1}$~) of rank $N$ in the real case, $1\le n\le s$~, $N\in\NN$~, and 
from ``irreducible'' $k$-semimodules $F_{\omega ^1_n}$ (resp. $F_{\o \omega ^1_n}$~) of rank $N\cdot m^{(n)}$ in the complex 
case, where $m ^{(n)}=\sup(m_n+1)$ denotes the multiplicity of the $n$-th real completion covering its complex equivalent as it will be seen in the following.

So, the ranks of the real pseudo-ramified completions $F^+_{v_n}$ (resp. $F^+_{\o v_n}$~) will be given by integers modulo $N$  according to:
\[\begin{array}[t]{lllll}
[F^+_{v_n}:k] & =*+n\cdot N &\qquad & \text{(resp.} \quad [F^+_{\o v_n}:k] & =*+n\cdot N\\
& \simeq n\cdot N &&& \simeq n\cdot N\ )\end{array}
\]
where $*$ denotes an integer inferior to $N$~.

In the complex case, the ranks of the pseudo-ramified completions $F_{\omega _n}$ (resp. $F_{\o\omega _n}$~) will also be given by the integers modulo $N$ according to:
\[\begin{array}[t]{lllll}
[F_{\omega _n}:k] & =*+n\cdot N\cdot m^{(n)} &\qquad & \text{(resp.} \quad [F_{\o \omega _n}:k] & =*+n\cdot N\cdot m^{(n)}\\
& \simeq n\cdot N\cdot m^{(n)} &&& \simeq n\cdot N\cdot m^{(n)}\ )\end{array}\]
in such a way that each complex completion $F_{\omega _n}$ (resp. $F_{\o\omega _n}$~) be covered by the set of $m^{(n)}$ real completions $F^+_{v_{n,m_n}}$ (resp. $F^+_{\o v_{n,m_n}}$~).

Indeed, as a place is an equivalence class of completions, we have to consider at each real place $v_n$ (resp. $\o v_n$~) a set of real completions $\{F^+_{v_n,m_n}\}_{m_n}$~, $m_n\in \NN$ (resp. $\{F^+_{\o v_n,m_n}\}_{m_n}$~), equivalent to $F^+_{v_n}$ (resp. $F^+_{\o v_n}$~) (with $m_n=0$~), and characterized by the same rank as $F^+_{v_n}$ (resp. $F^+_{\o v_n}$~).

Similarly, at each complex place $\omega _n$ (resp. $\o \omega _n$~), a set of complex completions\linebreak
$\{F_{\omega _{n,m_{\omega _n}}}\}_{m_{\omega _n}}$~, $m_{\omega _n}\in \NN$ (resp. $\{F_{\o \omega _ {n,m_{\omega _n}}}\}_{m_{\omega _n}}$~), equivalent to 
$F_{\omega _n}$ (resp. $F_{\o \omega _n}$~) and characterized by the same rank as $F_{\omega _n}$ 
(resp. $F_{\o\omega _n}$~) has to be considered.

Let $F_\omega =\{F_{\omega _1},\cdots,F_{\omega _{n,m_{\omega _n}}},\cdots,F_{\omega _{s,n_s}}\}$
(resp. $F_{\o\omega} =\{F_{\o\omega _1},\cdots,F_{\o\omega _{n,m_{\o\omega _n}}},\cdots,F_{\omega _{s,n_s}}\}$~) denote the set of complex pseudo-ramified completions of the number semifield $\widetilde F_L$ (resp. $\widetilde F_R$~) at the set of complex places $\omega $ (resp. $\o\omega $~) and let
$F^+_v=\{F^+_{v_1},\cdots,F^+_{v_{n,m_n}},\cdots,F^+_{v_{s,m_s}}\}$ (resp.
$F^+_{\o v}=\{F^+_{\o v_1},\cdots,F^+_{\o v_{n,m_n}},\cdots,F^+_{\o v_{s,m_s}}\}$~) be the corresponding real pseudo-ramified completions of the number semifield $\widetilde F^+_L$ (resp. $\widetilde F^+_L$~) at the set of real places $v$ (resp. $\o v$~).

Then, the direct sum of the complex pseudo-ramified completions is given by:
\[F_{\omega _\oplus}=\bigoplus_n \bigoplus_{m_{\omega _n}} F_{\omega _{n,m_{\omega _n}}}\qquad
\text{(resp.} \quad
F_{\o\omega _\oplus}=\bigoplus_n \bigoplus_{m_{\omega _n}}  F_{\o\omega _{n,m_{\omega _n}}}\ )\]
while the direct sum of the real pseudo-ramified completions is given by:
\[F^+_{v _\oplus}=\bigoplus_n \bigoplus_{m_{ _n}} F^+_{v _{n,m_{n}}}\qquad
\text{(resp.} \quad
F^+_{\o v _\oplus}=\bigoplus_n \bigoplus_{m_{ _n}}  F^+_{\o v _{n,m_{n}}}\ )\]

On the other hand, if
 the global class residue degree is a prime $p$ or an integer modulo $p$~, then 
$F^{+,nr}_{v_p}$ and $F^{+,nr}_{v_{p+i}}$ can refer respectively to a completion of an (unramified) $p$-adic (semi)field and to a completion of an 
extension of this $p$-adic (semi)field which can be identified to a finite extension of $\qit_p$~; indeed, 
$F^{+,nr}_{v_p}$ and $F^{+,nr}_{v_{p+i}}$ are pseudo-unramified local (semi)fields at $p$ elements which correspond to 
$p$ Galois automorphisms. More concretely, $F^{+,nr}_{v_p}$ and $F^{+,nr}_{v_{p+i}}$ are characterized by their 
global class residue degrees given by: \[ [F^{+,nr}_{v_p}:k] =f_{v_p}=p 
\quad \text{and} \quad
[F^{+,nr}_{v_{p+i}}:k] =f_{v_{p+i}}=k\ p+i'=p+i\;, \]
\mbox{} \hfill $ k\in\nit\;, \quad 0\le i'\le p-1\;, \quad 1\le i\le \infty$~. \qquad  \vskip 11pt

Remark that $F^{+,nr}_{v_p}$ and $F^{+,nr}_{v_{p+i}}$ are (unramified) $p$-adic fields if:
\Be
\item $f_{v_p}=p$ and $f_{v_{p+i}}=kp+i'=p^\alpha $~, $\alpha \in\NN$~, respectively.
\item the number of nonunits of $F^{+,nr}_{v_p}$ and of $F^{+,nr}_{v_{p+i}}$ is a power of $p$~.
\item $F^{+,nr}_{v_p}$ and $F^{+,nr}_{v_{p+i}}$ are considered as completions of $\QQ$ in the $p$-adic metric.
\Ee
\vskip 11pt

In the pseudo-ramified case, we should have: 
\[ [F^{+}_{v_{p+i}}:k] \simeq f_{v_{p+i}}\cdot N=(k\ p+i')\cdot N=(p+i)\cdot N\]
where $F^+_{v_{p+i}}$ is a  local (semi)field at $p$ elements corresponding to $(p+i)\cdot N$ 
 Galois automorphisms from a global point of view.
$F^+_{v_{p+i}}$ is a $p$-adic field if:
\Be
\item $f_{v_{p+i}}\cdot N\cdot\#(Nu)=p^\beta $~, $\beta \in \NN$~, where $\#(Nu)$ is the number of nonunits of $F^+_{v_{p+i}}$~.
\item $F^+_{v_{p+i}}$ is defined as a completion of $\QQ$ with respect to the $p$-adic metric.
\Ee
\vskip 11pt


If the number of left and right real places of the completions of $\widetilde F^+_L$ and $\widetilde F^+_R$ are equal respectively to the number of left and right complex places of the completions of $\widetilde F_L$ and of $\widetilde F_R$~, then $\widetilde F=\widetilde F_R\cup \widetilde F_L$ can be interpreted as a CM number field given by $\widetilde F=\widetilde F^+\cdot \widetilde E$ where $\widetilde F^+$ is the real number field and where $\widetilde E$ is an imaginary quadratic field.\vskip 11pt

Then a finite extension $K$ of $\qit_p$ can be identified to the completion of $\widetilde F^+$ in a place $v'_1$ above $p$~. And, as $p$ is decomposed in $E$ into two places $\wp$ and $\o \wp$~, the places of $\widetilde F$ dividing $\wp$ are divided into a set of left complex places $\omega '=\{\omega '_1,\cdots,\omega '_n,\cdots,\omega '_s\}$ above $\wp$ belonging to the completions of $\widetilde F_L$ and into a set of right conjugate places $\o\omega '=\{\o\omega '_1,\cdots,\o\omega '_n,\cdots,\o\omega '_s\}$ above $\o \wp$ belonging to the completions of $\widetilde F_R$~.\vskip 11pt

\subsection{Galois subgroups and global inertia subgroups}

Let $\Gal( \widetilde F^{+,nr}_{v_n}\big/k)$ (resp. $\Gal(\widetilde F^{+,nr}_{\o v_n}\big/k)$~) be the  Galois subgroup of the pseudo-unramified extension $\widetilde F^{+,nr}_{v_n}$ (resp. $\widetilde F^{+,nr}_{\o v_n}$~) of $k$ 
corresponding to the pseudo-unramified completion $F^{+,nr}_{v_n}$ (resp. $F^{+,nr}_{\o v_n}$~),
 and let $\Gal(\widetilde F^+_{v_n}\big/k)$ (resp. $\Gal(\widetilde F^+_{\o v_n}\big/k)$~) be the Galois subgroup of the pseudo-ramified extension $\widetilde F^+_{v_n}$ (resp. $\widetilde F^+_{\o v_n}$~) of $k$ corresponding to the pseudo-ramified completion $F^+_{v_n}$ (resp. $F^+_{\o v_n}$~).
\vskip 11pt

Then, the global inertia subgroup $I_{F^+_{v_n}}$ (resp. $I_{F^+_{\o v_n}}$~) of
$\Gal(\widetilde F^+_{v_n}\big/k)$ (resp. $\Gal(\widetilde F^+_{\o v_n}\big/k)$~), being the Galois subgroup of the irreducible extension $\widetilde F^+_{v^1_n}$ (resp. $F^+_{\o v^1_n}$~), can be defined by:
\begin{align*}
I_{F^+_{v_n}}
&= \Gal( \widetilde F^+_{v_n} \big/k)\Big/ \Gal(\widetilde F^{+,nr}_{v_n}\big/k)\\
\text{(resp.} \quad
I_{F^+_{\o v_n}}
&= \Gal(\widetilde F^+_{\o v_n}\big/k)\Big/ \Gal(\widetilde F^{+,nr}_{\o v_n}\big/k)\;)\end{align*}
leading to the exact sequence:
\begin{alignat*}{7}
1 \quad &\To \quad 
I_{F^+_{v_n}} & \quad 
\To \quad 
& 
 \Gal(\widetilde F^+_{v_n}\big/k) & \quad 
\To \quad 
& \Gal(\widetilde F^{+,nr}_{v_n}\big/k) & \quad \To \quad 
& 1\\
\text{(resp.} \quad
1 \quad 
&\To \quad 
I_{F^+_{\o v_n}} & \quad 
\To \quad 
& 
 \Gal(\widetilde F^+_{\o v_n}\big/k) & \quad 
\To \quad 
& \Gal(\widetilde F^{+,nr}_{\o v_n}\big/k) & \quad 
\To \quad 
 & 1\;).\end{alignat*}\vskip 11pt

The global inertia subgroup $I_{F^+_{v_n}} $ (resp. $I_{F^+_{\o v_n}}$~) of order $N$ can then be viewed as the normal subgroup or inner automorphisms of Galois with respect to the Galois subgroup $\Gal(\widetilde F^+_{v_n}\big/k) $ (resp. $\Gal(\widetilde F^+_{\o v_n}\big/k) $~) which can be considered as a subgroup of modular automorphisms of Galois \cite{Pie1}.\vskip 11pt

Similar conclusions are obtained for the Galois subgroups $\Gal(\widetilde F_{\omega _n}\big/k)$ (resp.\linebreak $\Gal(\widetilde F_{\o\omega _n}\big/k)$~) of the pseudo-ramified extension $\widetilde F_{\omega _n}$ (resp.$\widetilde F_{\o\omega _n}$~) of $k$~.
\vskip 11pt

As we are concerned with Galois class fields, all the  ``(pseudo-)ramification"" orders are supposed to be equal to the same positive integer $N$
in the real case and to the positive integers $N\cdot m^{(n)}$ in the complex case.
\vskip 11pt

\subsection{Global Weil semigroups}  

Let $ \widetilde F^+_{v_n,m_n}$ (resp. $\widetilde F^+_{\o v_n,m_n}$~) denote a real pseudo-ramified extension referring to the infinite place $v_n$ (resp. $\o v_n$~) and let
$\widetilde F^{+,nr}_{v_n,m_n}$ (resp. $\widetilde F^{+,nr}_{\o v_n,m_n}$~) denote the corresponding pseudo-unramified extension.

Then, we set:
\begin{align*}
\Gal(\widetilde F^+_{v}\big/k)
&= \txt\bigoplus\limits_{n,m_n} \Gal(\widetilde F^+_{v_n,m_n}\big/k)\\[11pt]
\text{(resp.} \quad
\Gal(\widetilde F^+_{\o v}\big/k)
&= \txt\bigoplus\limits_{n,m_n} \Gal(\widetilde F^+_{\o v_n,m_n}\big/k)\;)\end{align*}
as well as 
\begin{align*}
\Gal(\widetilde F^{+,nr}_{v}\big/k)
&= \txt\bigoplus\limits_{n,m_n} \Gal(\widetilde F^{+,nr}_{v_n,m_n}\big/k)\\[11pt]
\text{(resp.} \quad
\Gal(\widetilde F^{+,nr}_{\o v}\big/k)
&= \txt\bigoplus\limits_{n,m_n} \Gal(\widetilde F^{+,nr}_{\o v_n,m_n}\big/k)\;).\end{align*}
\vskip 11pt

As in the $p$-adic case, the Weil group is the Galois subgroup of the elements inducing on the residue field an integer power of a Frobenius element, we shall assume that, in the characteristic zero case, the Weil group will be the Galois subgroup of the pseudo-ramified extensions characterized by extension degrees $d=0\mod N$~.

In this respect, let $\dot{\widetilde F}^+_{v_n}$ (resp. $\dot{\widetilde F}^+_{\o v_n}$~) denote a pseudo-ramified Galois extension characterized by the degree
\[ [\dot{\widetilde F}^+_{v_n}:k] \equiv [\dot{\widetilde F}^+_{\o v_n}:k]=n\cdot N\]
in such a way that the sum of the  closed global Weil sub(semi)groups  is given by:
\begin{align*}
\WW_{F^+_v} &= \txt\bigoplus\limits_{n,m_n} \Gal(\dot{\widetilde F}^{+}_{v_{n,m_n}}\big/k)
\\
\text{(resp.} \quad
\WW_{F^+_{\o v}} &= \txt\bigoplus\limits_{n,m_n} \Gal(\dot{\widetilde F}^{+}_{\o v_{n,m_n}}\big/k)
\;).
\end{align*}
And, the product of  $W_{F^+_{\o v}}$ and of $W_{F^+_{v}}$ gives:
\[
\WW_{F^{+}_{\o v}}\times \WW_{F^{+}_{v}} =
\txt\bigoplus\limits_{n,m_n} \L(\Gal(\dot{\widetilde F}^{+}_{\o v_{n,m_n}}\big/k) \times 
\Gal(\dot{\widetilde F}^{+}_{v_{n,m_n}}\big/k)\R)
\;.\]\vskip 11pt

\subsection{Representations of algebraic bilinear semigroups}  

Let $B_L\big/\widetilde F^+_v$ denote a left division semialgebra of dimension $r$ over the semifield $\widetilde F_{v}^+ = \{\widetilde F^+_{v_1},\cdots,\widetilde F^+_{v_{s,m_s}}\}$ and let $B_R\big/\widetilde F^+_{\o v} $ denote the corresponding right division semialgebra of the same dimension $r$ over the semifield $\widetilde F_{\o v}^+ = \{\widetilde F^+_{\o v_1},\cdots,\widetilde F^+_{\o v_{s,m_s}}\}$ in such a way that: 
\Bi \item $B_R$ is the opposite semialgebra of $B_L$~; 
\item  the center of $B_L$ is in one-to-one correspondence with the center of $B_R$~. \Ei\vskip 11pt

If it is assumed that the division semialgebra $B_L$ is isomorphic to the matrix algebra 
$T_r(\widetilde F^+_v)$ of Borel upper triangular matrices over $\widetilde F^+_v$ and that the opposite division semialgebra $B_R$ is isomorphic to the matrix algebra $T^t_r(\widetilde F^+_{\o v})$ of Borel lower triangular matrices over $\widetilde F^+_{\o v}$~, then we have: 
\[ B_R\otimes B_L\simeq T^t_r(\widetilde F^+_{\o v})\times T_r(\widetilde F^+_v)
\equiv\GL_r(\widetilde F^+_{\o v}\times \widetilde F^+_v)\;.\]
So, the tensor product $B_R\otimes B_L$ of the division semialgebras $B_R$ and $B_L$ is  isomorphic to the product of the group $T^t_r(\widetilde F^+_{\o v})$ of lower triangular matrices by the 
group $T_r(\widetilde F^+_v)$ of upper triangular matrices. 
Such a product $T^t_r(\widetilde F^+_{\o v})\times T_r(\widetilde F^+_v)$ is denoted $\GL_r(\widetilde F^+_{\o v}\times \widetilde F^+_v)$ and corresponds to a general bilinear semigroup having a representation in the tensor product $ M_R\otimes  M_L$ of a right $T^t_r(\widetilde F^+_{\o v})$-semimodule $ M_R$ by a left $T_r(\widetilde F^+_v)$-semimodule $ M_L$~.\vskip 11pt

So, the general bilinear semigroup
\[ \GL_r(\widetilde F^+_{\o v}\times \widetilde F^+_v)= T^t_r(\widetilde F^+_{\o v})\times T_r(\widetilde F^+_v)\;,\]
has a representation in a $\GL_r(\widetilde F^+_{\o v}\times \widetilde F^+_v)$-bisemimodule 
$ M_R\otimes  M_L$ of dimension $r^2$~.\vskip 11pt

\subsection{Proposition}  

{\em Let $B_L\big/\widetilde F^+_v$ (resp. $B_R\big/\widetilde F^+_{\o v}$~) be a \lr division semialgebra of dimension $r$ 
over the extension semifield $\widetilde F^+_v$ (resp. $\widetilde F^+_{\o v}$~) such that
\[
B_L \simeq T_r(\widetilde F^+_v)\;,\qquad B_R \simeq T^t_r(\widetilde F^+_{\o v})\;.\] Then, we have \[ B_R\otimes B_L \simeq \GL_r(\widetilde F^+_{\o v}\times \widetilde F^+_v)\] implying that the general bilinear semigroup $\GL_r(\widetilde F^+_{\o v}\times \widetilde F^+_v)$ has the natural bilinear  Gauss decomposition:
\begin{align*}
g_r(\widetilde F^+_{\o v}\times \widetilde F^+_v) &= t^t_r(\widetilde F^+_{\o v})\times t_r(\widetilde F^+_v)\\
&= [u^t_r(\widetilde F^+_{\o v})\times u_r(\widetilde F^+_v)]\times [d_r(\widetilde F^+_{\o v})\times d_r(\widetilde F^+_v)]\end{align*}
for any regular matrix 
\Bi 
\item $g_r(\widetilde F^+_{\o v}\times \widetilde F^+_v)
\in\GL_r(\widetilde F^+_{\o v}\times \widetilde F^+_v)$~; 

\item $u_r(\widetilde F^+_{\o v})\in UT_r(\widetilde F^+_v)$ 
where $UT_r(\widetilde F^+_v)$ is the group of upper unitriangular  matrices; 

\item $d_r(\widetilde F^+_v)\in D_r(\widetilde F^+_v)$ where $D_r(\widetilde F^+_v)$ is the group of diagonal matrices. \Ei}\vskip 11pt

\bpr The classical case
\[ B^e\simeq \End_{\widetilde F^+}(B)\simeq \End_{\widetilde F^+}((\widetilde F^+)^r)\simeq \Ms_r(\widetilde F^+)\]
where:
\Bi
\item $B^e=B\times_FB^{\rm op}$ with $B$ a central simple $F$-algebra and with $B^{\rm op}$ the opposite algebra,
\item $\End_{\widetilde F^+}(B)$ is isomorphic to $B^e$ if $B$ is an Azumaya algebra (then, $B$ is a projective $F$-module),
\item $\widetilde F^+$ is a ring in such a way that $\widetilde F^+= \widetilde F^+_{\o v}\cup \widetilde F^+_v$~,
\item $\Ms_r(\widetilde F^+)$ is the ring of $r\times r$ matrices with entries in $\widetilde F^+$~,
\Ei
becomes in the considered bilinear case:
\[ B^e\equiv B_R\otimes B_L\simeq \End_{\widetilde F^+_{\o v}\times \widetilde F^+_v}(B_R\times B_L)\simeq \End_{\widetilde F^+_{\o v}\times \widetilde F^+_v} ((\widetilde F^+_{\o v}\times \widetilde F^+_v)^r)\simeq \GL_r(\widetilde F^+_{\o v}\times \widetilde F^+_v)\]
where:
\Bi
\item $B_L$ (resp. $B_R$~) is a (division) central simple $\widetilde F^+_v$ (resp. $\widetilde F^+_{\o v}$~)-semialgebra of dimension $r$ over $\widetilde F^+_v$ (resp. $\widetilde F^+_{\o v}$~);
\item $\End_{\widetilde F^+_{\o v}\times \widetilde F^+_v}(B_R\times B_L) \simeq 
\End_{\widetilde F^+_{\o v}\times \widetilde F^+_v}(B_L)
\simeq 
\End_{\widetilde F^+_{\o v}\times \widetilde F^+_v}(B_R)$~.\epr
\Ei
\vskip 11pt

\subsection{Proposition}

{\em Let $F^+=F^+_R\cup F^+_L\simeq \widetilde F^+_{\o v}\cup \widetilde F^+_v$ be a compact real algebraically closed field defined from the symmetric algebraically closed semifields $F^+_R\simeq \widetilde F^+_{\o v}$ and $F^+_L\simeq \widetilde F^+_v$~.

Let $\GL_r(F^+)$ be the group of invertible $(r\times r)$  matrices with entries in $F^+$~.  The algebraic general linear group $\GL_r(F^+)$ has for representation space a vectorial space $V$ of dimension $r^2$ isomorphic to $(F^+)^{r^2}$~.

On the other hand, let $\GL_r(F^+_R\times F^+_L)$ be the algebraic bilinear semigroup, introduced in sections 1.7 and 1.8, and having for representation space the $\GL_r(F^+_R\times F^+_L)$-bisemimodule $M_R\otimes M_L$ of dimension $r^2$~.

Then, we can state the one-to-one correspondence:
\begin{align*}
I_{\GL}: \quad \GL_r(F^+) & \begin{CD}@>\sim>>\end{CD} \GL_r(F^+_R\times F^+_L)\;,\\
V &\begin{CD}@>\sim>>\end{CD}M_R\otimes M_L\;, \end{align*}
between the algebraic linear group $\GL_r(F^+)$ and the algebraic bilinear semigroup $\GL_r(F^+_R\times F^+_L)$ sending the vectorial space $V$ into the $\GL_r(F^+_R\times F^+_L)$-bisemimodule $M_R\otimes M_L$~.
}
\vskip 11pt

{\parindent=0pt
\bpr $M_L$ (resp. $M_R$~) is a $r$-dimensional vector semispace over $F^+_L$ (resp. $F^+_R$~).
Every basis of $M_L$ (resp. $M_R$~) defines an isomorphism of the group $T_r(M_L)$ (resp. $T_r^t(M_R)$~) 
of the automorphisms of $M_L$ (resp. $M_R$~) with $T_r(F^+_L)$ (resp. $T_r^t(F^+_R)$~) $\subset \GL_r(F^+_R\times F^+_L)$~.

And, every basis of $(M_R\otimes M_L)$ defines an isomorphism of the group $\GL_r(M_R\otimes M_L)$ of the automorphisms of $(M_R\otimes M_L)$ with $\GL_r(F^+_R\times F^+_L)$~.

Now, $\GL_r(F^+_R\times F^+_L)$ has the bilinear Gauss decomposition:
\[ g_r(F^+_R\times F^+_L) =(u^t_r(F^+_R)\times u_r(F^+_L))\times (d_r(F^+_R)\times d_r(F^+_L))\]
for any matrix $g_r(F^+_R\times F^+_L) \in\GL_r(F^+_R\times F^+_L)$~,\\
while $\GL_R(F^+)$ has the Gauss decomposition:
\[  g_r(F^+)=(u^t_r(F^+)\times u_r(F^+))\times d_r(F^+)\]
for $g_r(F^+)\in\GL_r(F^+)$~.

If we take into account the maps \eqref{eq:*}
\begin{equation}\label{eq:*}
\begin{cases}
u^t_r(F^+)&\To u^t_r(F^+_R)\;, \\
u_r(F^+)&\To u_r(F^+_L)\;, \\
d_r(F^+)&\To d_r(F^+_R\times F^+_L)\;, \end{cases}\tag{$\star$}\end{equation}
then,
\Be
\item $V\simeq M_R\otimes M_L$ in such a way that the basis of the $n^2$-dimensional vector space $V$ corresponds to the basis of the $n^2$-dimensional $\GL_r(F^+_R\times F^+_L)$-bisemimodule $M_R\otimes M_L$~;
\item the linear algebraic group $\GL_r(F^+)$ is covered by the bilinear algebraic semigroup $\GL_r(F^+_R\times F^+_L)$~.
\epr
\Ee}
\vskip 11pt 

\subsection{Proposition}

{\em Let $F^+=F^+_R\cup F^+_L$ be a symmetric algebraic extension field.

If the Gauss decompositions of the linear algebraic group $\GL_r(F^+)$ and of the bilinear algebraic semigroup $\GL_r(F^+_R\times F^+_L)$ correspond under the conditions \eqref{eq:*} of proposition 1.9, then:
\Be
\item $\GL_r(F^+) \begin{CD}@>\sim>>\end{CD} \GL_r(F^+_R\times F^+_L)$
expressing that $\GL_r(F^+) $ is covered by $\GL_r(F^+_R\times F^+_L)$~;
\item the $r^2$-dimensional representation space $V=\Repsp(\GL_r(F^+))$ of $\GL_r(F^+)$ coincides with the  $r^2$-dimensional representation space $M_R\otimes M_L=\Repsp(\GL_r(F^+_R\times F^+_L))$ of $\GL_r(F^+_R\times F^+_L)$~.
\Ee}
\vskip 11pt

We are then led to a bilinear version of the Wedderburn theorem \cite{K-M-R-T}:
\vskip 11pt

\subsection{Proposition}  

{\em Let $B_L$ and $B_R$ be the division semialgebras respectively over the semifields $F^+_L\simeq \widetilde F^+_v$ and $F^+_R\simeq \widetilde F^+_{\o v}$~.  
The following conditions are then equivalent:

\Bean

\item $B_R\otimes B_L$ {is central simple} of dimension {$r ^2$~}.

\item $B_R\otimes B_L\simeq \GL_r(F^+_R\times F^+_L)$~,

$B_L\otimes B_R\simeq \GL_r(F^+_L\times F^+_R)$~.

\item The canonical map \quad $\begin{cases}B_R\otimes B_L\to \End_{F^+_R\times F^+_L}(B_L)\;,\\ B_L\otimes B_R\to 
\End_{F^+_L\times F^+_R}(B_R)\;,\end{cases}$
which relates respectively to  \quad $\begin{cases} b_R\otimes b_L\\ b_L\otimes b_R\end{cases}$
the map  \quad $\begin{cases} x^2\to b_Rx^2b_L\;,\\ x^2\to b_Lx^2b_R\;,\end{cases}$ \quad  is an isomorphism.

\item $\begin{cases} B_R/F^+_R\simeq T^t_r(F^+_R)\\
B_L/F^+_L\simeq T_r(F^+_L)\end{cases} $
implies \quad $B_R/F^+_R\otimes B_L/F^+_L\simeq T^t_r(F^+_R)\times T_r(F^+_L)$~.\Ee}
\vskip 11pt

\subsection{Proposition}  

{\em \Be
\item Every involution of the first kind on the central simple $F^+_L$-semialgebra $B_L$ (resp. $F^+_R$-semialgebra $B_R$~) is defined by a mapping $\sigma_{L\to R}:B_L/F^+_L\to B_R/F^+_R$ (resp. $\sigma_{R\to L}:B_R/F^+_R\to B_L/F^+_L$~) such that: $t_r(F^+_L)\to t^t_r(F^+_R)$ (resp. $t_r^t(F^+_R)\to t_r(F^+_L)$~) transposes every upper (resp. lower) Borel upper triangular matrix $t_r(F^+_L)\in T_r(F^+_L)$ (resp.
$t_r^t(F^+_R)\in T_r^t(F^+_R)$~).
\item Every involution of the second kind $\sigma_{R\times L\to L\times R}$ on the tensor product of the central simple $F^+_R$-semialgebra $B_R$ by the central simple $F^+_L$-semialgebra $B_L$ is given by the exchange involution:
\[ \sigma_{R\times L\to L\times R}:\quad B_R/F^+_R\otimes B_L/F^+_L \To B_L/F^+_L\otimes B_R/F^+_R\]
such that
\[t_r^t(F^+_R)\times t_r(F^+_L)\To t_r(F^+_L)\times t_r^t(F^+_R)\]
resulting in a double involution of the first kind:
\[\sigma_{R\times L\to L\times R}=\sigma_{R\to L}\circ\sigma_{L\to R}\;.\]\Ee}
\vskip 11pt
\bpr  This is a direct consequence of Propositions 1.9 and 1.12.  However, let us prove it more explicitly inspired by the book of \cite{K-M-R-T}.
\\
Let $\gamma:M_R\times M_L\to F^+_R\times F^+_L$ be a bilinear form on the product of the semifields $F^+_R$ and $F^+_L$ and let
\begin{align*}
\widehat \gamma : \quad & M_L\To M_R = \Hom(M_L,F^+_R)\;, \\
\widehat \gamma^{-1} : \quad & M_R\To M_L = \Hom(M_R,F^+_L)\;.\end{align*}
For any $f_L\in \End_{F^+_L}(M_L)$ (resp. $f_R\in \End_{F^+_R}(M_R)$~), we define $\sigma_\gamma(f_L)\in \End_{F^+_L}(M_L)$ by $\sigma_\gamma(f_L)=\widehat \gamma^{-1}\circ f_R\circ \widehat\gamma$ (resp. $\sigma_{\gamma^{-1}}(f_R)\in \End_{F^+_R}(M_R)$ by $\sigma_{\gamma^{-1}}(f_R)=\widehat \gamma\circ f_L\circ \widehat\gamma^{-1}$~).  It is then clear that any involution of this first kind $\sigma_{L\to R}$ (resp. $\sigma_{R\to L}$~) on a division semialgebra $B_L$ (resp. $B_R$~) can be explained as follows:
\begin{align*}
(B_L,\sigma_{L\to R})\simeq \End_{B_L}(M_L) &\simeq (T_r(F^+_L),\sigma_\gamma)\;, \\
(B_R,\sigma_{R\to L})\simeq \End_{B_R}(M_R) &\simeq (T_r^t(F^+_R),\sigma_{\gamma^{-1}})\;.\tag*{\eop}\end{align*}
\vskip 11pt

 \section{Modular representations of general bilinear semigroups in terms of global elliptic bisemimodules}  Until the end of this paper, all developments will refer to the real and complex dimensions which are inferior or equal to 2. \vskip 11pt

 Let us summarize our terminology:
 
 \Bi
 \item 
  We consider the left and right symmetric completions $F_\omega$ and $F_{\o \omega }$ decomposing respectively according to the set of places $\omega =\{\omega _1,\cdots,\omega _n,\cdots,\omega _s\}$ and\linebreak $\o\omega =\{\o\omega _1,\cdots,\o\omega _n,\cdots,\o\omega _s\}$~) of the complex algebraic extension semifields $ \widetilde F_L$ and $\widetilde F_R$ of a global number field of characteristic zero $k$~. If the algebraic extension semifields $\widetilde F^+_L$ and  $\widetilde F^+_R$ are real, then the corresponding left and right symmetric completions $F^+_v$ and $F^+_{\o v}$ decompose according to the set of real places $v=\{ v_1,\cdots,v_n,\cdots,v_s\}$ and $\o v=\{ \o v,\cdots,\o v_n,\cdots,\o v_s\}$~.\vskip 11pt

\item In the real case, let $B_{F^+_v}$ (resp. $B_{F^+_{\o v}}$~)be a left (resp. right) division semialgebra  over the completions $F^+_v$ (resp. $F^+_{\o v}$~).  Assuming the isomorphisms:
\[ B_{F^+_v} \simeq T_2(F^+_v) \qquad \text{and} \qquad 
B_{F^+_{\o v}}\simeq T^t_2(F^+_{\o v}) \;, \]
we have 
\[ B_{F^+_{\o v}}\otimes B_{F^+_v}\simeq T^t_2(F^+_{\o v})\times T_2(F^+_v)\equiv
\GL_2(F^+_{\o v}\times F^+_v)\]
where $\GL_2(F^+_{\o v}\times F^+_v)$ is a general bilinear semigroup, with entries in the product of the 
semifield completions $F^+_{\o v}\times F^+_v$~, generated from the product, right by left, 
of triangular matrix (semi)groups of order 2 which leads to a bilinear Gauss decomposition as 
introduced in proposition 1.9. 

\item $GL_2(F^+_{\o v}\times F^+_v)$ has a representation in the $B_{F^+_{\o v}}\otimes B_{F^+_v}$-bisemimodule $M_{F^+_{\o v}}\otimes M_{F^+_v}$ 
composed of the set $\{M_{F^+_{\o v_{n,m_n}}} \otimes
M_{F^+_{v_{n,m_n}}}\}^s_{n=1,m_n}$ of subbisemimodules $M_{F^+_{\o v_{n,m_n}}}\otimes
M_{F^+_{v_{n,m_n}}}$ following the places $\o v_n$ and $v_n$ considered with their $m_n$ representatives.

\item Similarly, if $F^+_{v_\oplus}=\bigoplus_{n,m_n} F^+_{v_{n,m_n}}$ (resp.
$F^+_{\o v_\oplus}=\bigoplus_{n,m_n} F^+_{\o v_{n,m_n}}$~) denotes the sum of the real pseudo-ramified completions, the algebraic bilinear semigroup $\GL_2(F^+_{\o v_\oplus}\times F^+_{v_\oplus})$~, with entries in the product of $F^+_{\o v_\oplus}$ by $F^+_{v_\oplus}$~, has a representation in the $B_{F^+_{\o v_\oplus}}\otimes B_{F^+_{v_\oplus}}$-bisemimodule $M_{F^+_{\o v_\oplus}}\otimes M_{F^+_{v_\oplus}}$ which decomposes according to:
\[  M_{F^+_{\o v_\oplus}}\otimes M_{F^+_{v_\oplus}}
= \bigoplus_{n,m_n}\L( M_{F^+_{\o v_{n,m_n}}}\otimes M_{F^+_{v_{n,m_n}}}\R)\]
in such a way that:
\[B_{F^+_{\o v_\oplus}}\otimes B_{F^+_{v_\oplus}}
\simeq T^t_2( F^+_{\o v_\oplus}) \times T_2(F^+_{v_\oplus})\equiv \GL_2
( F^+_{\o v_\oplus}\times F^+_{v_\oplus})\;.\]

\item The representation space $(M_{F^+_{\o v}}\otimes M_{F^+_{v}})$ of $\GL_2
(F^+_{\o v}\times F^+_{v})$ is homomorphic to the representation space
$(M_{F^+_{\o v_\oplus}}\otimes M_{F^+_{v_\oplus}})$ of $\GL_2
(F^+_{\o v_\oplus}\times F^+_{v_\oplus})$~:
\[\begin{CD}
\GL_2
(F^+_{\o v}\times F^+_{v}) @>>> \GL_2
(F^+_{\o v_\oplus}\times F^+_{v_\oplus})\\
@VVV @VVV\\
(M_{F^+_{\o v}}\otimes F^+_{v}) @>>> (M_{F^+_{\o v_\oplus}}\otimes F^+_{v_\oplus})\end{CD}\]
\Ei
\vskip 11pt

\subsection{Definition: Left and right real (pseudo-)ramified lattices} 

 In this respect, a \lr maximal order $ \Lambda ^{(2)}_{v}$ (resp. $ \Lambda ^{(2)}_{\o v}$~) over $\ZZ\big/N\ \ZZ$  in the \lr division semialgebra $B_{F^+_v}$ (resp. $B_{F^+_{\o v}}$~) can be introduced by: \Bean \item the isomorphisms
\[ \Lambda ^{(2)}_{v}\simeq T_2(\ZZ\big/N\ \ZZ) \qquad \text{and} \qquad
\Lambda ^{(2)}_{\o v}\simeq T^t_2(\ZZ\big/N\ \ZZ)\]
where $T_2(\ZZ\big/N\ \ZZ)$ is the subgroup of upper triangular matrices over the integers modulo N as it will be developed in section 2.3.
 \item its representation in a \lr lattice decomposing according to the set of places $v$ (resp. $\o v$~) of $F^+_v$ (resp. $F^+_{\o v}$~) following:
\[ \Lambda ^{(2)}_{v_\oplus} = \txt\bigoplus\limits_{n,m_n} \Lambda ^{(2)}_{v_{n,m_n}}\qquad
\text{(resp.} \quad 
\Lambda ^{(2)}_{\o v_\oplus} = \txt\bigoplus\limits_{n,m_n} \Lambda ^{(2)}_{\o v_{n,m_n}} \;)\]
where $\Lambda ^{(2)}_{v_{n,m_n}}$ (resp. $\Lambda ^{(2)}_{\o v_{n,m_n}}$~) is a \lr (Hecke) pseudo-ramified sublattice associated with the $n$-th place $v_n$ (~$m_n$-th representative) and corresponding to the $T_2(F^+_{v_n})$-subsemimodule $M_{F^+_{v_{n,m_n}}}$ (resp. $T^t_2(F^+_{\o v_n})$-subsemimodule $M_{F^+_{\o v_{n,m_n}}}$~).\vskip 11pt
\Ee

\subsection{Real (pseudo-)ramified bilattices} 

 The tensor product $\Lambda ^{(2)}_{\o v}\otimes \Lambda ^{(2)}_{v}$ of the right (pseudo-)\linebreak ramified lattice
$\Lambda ^{(2)}_{\o v}$ by the left (pseudo-)ramified lattice $\Lambda ^{(2)}_{v}$ constitutes the representation space of the (bilinear) arithmetic group $\GL_2(\ZZ\big/N\ \ZZ)^2\equiv T^t_2(\ZZ\big/N\ \ZZ)\times T_2(\ZZ\big/N\ \ZZ)$~. Indeed we have 
that:
\[ \Repsp(\GL_2(\ZZ\big/N\ \ZZ)^2)= \Lambda ^{(2)}_{\o v}\otimes \Lambda ^{(2)}_{v}\]
is homomorphic to $\Lambda ^{(2)}_{\o v_\oplus}\otimes \Lambda ^{(2)}_{v_\oplus}$ which
decomposes into the direct sum of bisublattices $\Lambda ^{(2)}_{\o v_{n,m_n}} \otimes
\Lambda ^{(2)}_{v_{n,m_n}}$~, according to the set of biplaces $\o v\times v$ of $F^+_{\o v}\times F^+_v$~:
\[\Lambda ^{(2)}_{\o v_\oplus} \otimes
\Lambda ^{(2)}_{v_\oplus}=\txt\bigoplus\limits_n\bigoplus\limits_{m_n} 
(\Lambda ^{(2)}_{\o v_{n,m_n}} \otimes
\Lambda ^{(2)}_{v_{n,m_n}})\;.\]
On the other hand, the (bilinear) arithmetic group $\GL_2(\ZZ\big/N\ \ZZ)^2$ has the following bilinear Gauss decomposition:
\[ g_2(\ZZ\big/N\ \ZZ)^2 = [d_2(\ZZ\big/N\ \ZZ)\times d_2(\ZZ\big/N\ \ZZ)] \times [u^t_2(\ZZ\big/N\ \ZZ)\times u_2(\ZZ\big/N\ \ZZ)]\] for every element $g_2(\ZZ\big/N\ \ZZ)^2\in \GL_2(\ZZ\big/N\ \ZZ)^2$ where \Bi \item $d_2(\ZZ\big/N\ \ZZ)$ is a diagonal matrix of order 2, also called a split Cartan subgroup element; \item $u_2(\ZZ\big/N\ \ZZ)$ is a two-dimensional upper unitriangular matrix; \item $u ^t_2(\ZZ\big/N\ \ZZ)$ is a two-dimensional lower unitriangular matrix. \Ei\vskip 11pt

\subsection{Complex (pseudo-)ramified bilattices}

The complex case can be handled similarly as it was done for the real case.

\Bean
\item So, let $B_{F_\omega }$ (resp. $B_{F_{\o\omega} }$~) be a \lr division semialgebra over the complex completions ${F_\omega }$ (resp. ${F_{\o\omega }}$~).

Setting $B_{F_\omega } \simeq T_2(F_\omega )$ and $ B_{F_{\o\omega }} \simeq T_2(F_{\o\omega })$~, the complex bilinear semigroup $\GL_2(F_{\o\omega }\times F_\omega )$ with entries in the product 
$F_{\o\omega }\times F_\omega $ of the complex semifield completions can be introduced by:
\[ B_{F_{\o\omega }} \otimes B_{F_{\omega }} \simeq T^t_2(F_{\o\omega })\times T_2(F_{\omega })
\equiv \GL_2 (F_{\o\omega }  \times F_{\omega })\;.\]
The  $\GL_2 (F_{\o\omega }  \times F_{\omega })$-bisemimodule $M_{F_{\o\omega }}\otimes M_{F_\omega }$ is homomorphic to $M_{F_{\o\omega _\oplus}}\otimes M_{F_{\omega _\oplus}}$ which
decomposes into:
\[ M_{F_{\o\omega_\oplus }}\otimes M_{F_{\omega _\oplus}}=\txt\bigoplus\limits_{n,m_{\omega _n}} 
(M_{F_{\o\omega _{n,m_{\omega _n}}}}\otimes M_{F_{\omega_{n,m_{\omega _n}}} })\;,\]
following the complex places $\omega _n\in \omega $ (resp. $\o\omega _n\in\o\omega $~) 
with representatives $m_{\omega _n}$~, constitutes a representation of 
$\GL_2(F_{\o\omega} \times F_\omega )$~.

\item This can be justified if a \lr lattice $\Lambda ^{(2)}_{\omega }$ (resp. $\Lambda ^{(2)}_{\o\omega }$~) over $\ZZ\big/N\ \ZZ$ in the \lr division semialgebra $B_{F_\omega }$ (resp. $B_{F_{\o\omega} }$~) is introduced by the isomorphisms:
\[ \Lambda ^{(2)}_{\omega }\simeq T_2(\ZZ\big/N\ \ZZ)
\qquad \text{(resp.} \quad
\Lambda ^{(2)}_{\o\omega }\simeq T^t_2(\ZZ\big/N\ \ZZ)\ )\]
where $ T_2(\ZZ\big/N\ \ZZ)$ is the subgroup of upper triangular matrices over the integers modulo $N$ according to section 1.4.

This \lr lattice $\Lambda ^{(2)}_{\omega }$ (resp. $\Lambda ^{(2)}_{\o\omega }$~) 
is homomorphic to $\Lambda ^{(2)}_{\omega_\oplus }$ (resp. $\Lambda ^{(2)}_{\o\omega _\oplus}$~) which decomposes following the set of complex places $\omega $ (resp. $\o\omega $~) according to:
\[ \Lambda ^{(2)}_{\omega _\oplus} =\txt\bigoplus\limits_{n,m_{\omega _n}}\Lambda ^{(2)}_{\omega_{n,m_{\omega _n}} } \qquad
\text{(resp.} \quad 
\Lambda ^{(2)}_{\o\omega_\oplus } =\txt\bigoplus\limits_{n,m_{\omega _n}}\Lambda ^{(2)}_{\o\omega_{n,m_{\omega _n}} } \ ).\]

\item The tensor product $\Lambda ^{(2)}_{\o\omega } \otimes\Lambda ^{(2)}_{\omega } $ 
constitutes the representation space of $\GL_2 (\ZZ\big/N\ \ZZ)^2\equiv T^t_2 (\ZZ\big/N\ \ZZ) \times T_2(\ZZ\big/N\ \ZZ)$ in such a  way that:
$\Lambda ^{(2)}_{\o\omega_\oplus } \otimes\Lambda ^{(2)}_{\omega_\oplus } $  decomposes into the bisublattices according to:
\[ \Lambda ^{(2)}_{\o\omega _\oplus} \otimes\Lambda ^{(2)}_{\omega _\oplus} 
=\txt\bigoplus\limits_{n,m_{\omega _n}}(\Lambda ^{(2)}_{\o\omega_{n,m_{\omega _n}} } \otimes
 \Lambda ^{(2)}_{\omega_{n,m_{\omega _n}} } ).\]
\Ee
\vskip 11pt



\setcounter{subsection}{3}  
\subsection{Borel-Serre compactification type of the (pseudo-)ramified\protect\linebreak  lattice bisemispace}

Let us introduce the (pseudo-)ramified complex lattice bisemispace $X_{S_{R\times L}}$ as the quotient semigroup:
\[ X_{S_{R\times L}} = \GL_2(\widetilde F_{R}\times \widetilde F_L )\big/\GL_2(\ZZ\big/N\ \ZZ)^2
\approx M_{F_{\o\omega }}\otimes M_{F_\omega }\]
where the complex algebraic bilinear semigroup $\GL_2(\widetilde F_{R }\times \widetilde F_L )$ is taken over the product $\widetilde F_{R }\times \widetilde F_L$ of the symmetric (algebraically closed) splitting  semifields $ \widetilde F_R$ and $\widetilde F_L$~.\vskip 11pt

$X_{S_{R\times L}}$ has a representation in a $\GL_2(F_{\o \omega }\times F_\omega )$-bisemimodule $M_{F_{\o \omega }}\otimes M_{F_\omega }$ in the sense that:
\[ X\SRTL \approx M_{F_{\o\omega }}\otimes M_{F_\omega }\;.\]
The Borel-Serre toroidal compactification of $X_{S_{R\times L}}$ can be considered as the toroidal projective isomorphism of compactification given by:
\[ \gamma ^c_{R\times L}:\quad X_{S_{R\times L}}\To \o X_{S_{R\times L}}\]
where
\[ \o X_{S_{R\times L}}= \GL_2(F^T_{R }\times F^T_L)\big/\GL_2(\ZZ\big/N\ \ZZ)^2
\approx M_{F^T_{\o\omega }}\otimes M_{F^T_\omega }\] such that:
\Bi
\item $X_{S_{R\times L}}$  may be viewed as the interior of $\o X_{S_{R\times L}}$ in the sense that $\gamma ^c_{R\times L}$ is an inclusion isomorphism: \quad $X_{S_{R\times L}}\hookrightarrow 
\o X_{S_{R\times L}}$~;
\item $F^T_R $ and $F^T_L$ are    ``toroidal''  symmetric algebraic(ally closed) semifields;
\item $M_{F^T_{\o\omega }}\otimes M_{F^T_\omega }$ is the toroidal equivalent of $M_{F_{\o\omega }}\otimes M_{F_\omega }$~.
\Ei\vskip 11pt


\subsection{Proposition} 

{\em If the multiplicity of the complex places is equal to one, the boundary $\partial\o X\SRTL $ of $\o X\SRTL$ is introduced as resulting from the surjective morphism:
\[ \gamma ^\delta \RTL: \quad \o X\SRTL\to \partial\o X\SRTL\]
sending the complex toroidal (pseudo-)ramified lattice bisemispace $\o X\SRTL$ into its boundary $\partial\o X\SRTL$ given by
\[\partial \o X\SRTL=\GL_2( F^{+,T}_{R} \times F^{+,T}_{L} )\big/\GL_2(\ZZ\big/N\ \ZZ)^2
\approx M_{F^{+,T}_{\o v}}\otimes M_{F^{+,T}_v}\;, \]
where 
\Bi
\item $F^{+,T}_R$ and $F^{+,T}_L$ are ``toroidal'' symmetric algebraic(ally closed) semifields;
\item $M_{F^{+,T}_{\o v}} \otimes M_{F^{+,T}_{v}}$ is the toroidal equivalent of
$M_{F^{+}_{\o v}}\otimes M_{F^{+}_{v}}$~;
\item
$F^{+,T}_{v} $ (resp. $F^{+,T}_{\o v} $~) are the toroidal completions associated with 
$F^+_v$ (resp. $F^+_{\o v}$~) at the set of real places $v$ (resp. $\o v$~) {\em \cite{Pie2}}. 
\Ei
}
\vskip 11pt

{\parindent=0pt
\bpr \Bean\item According to section 2.3, the complex toroidal $\GL_2({F^T_{\o \omega }} \times {F^T_{\omega }})$-bisemimodule $M_{F^T_{\o \omega }} \otimes M_{F^T_{\omega }} $ is the representation space of $\GL_2(F^T_{\o\omega }\times F^+_\omega )$ in such a way that:
\[ \Repsp(\GL_2( F^T_{\o \omega_\oplus } \times F^T_{\omega _\oplus}))
= M_{F^T_{\o \omega_\oplus }} \otimes M_{F^T_{\omega _\oplus}} 
=\txt\bigoplus\limits_{n,m_{\omega _n}} (M_{F^T_{\o \omega_{n,m_{\omega _n}} }} \otimes 
M_{F^T_{\omega_{n,m_{\omega _n} }} })\]
where $F^T_{\o\omega }$ and $F^T_\omega $ denote the sets of toroidal completions corresponding respectively to $F_{\o\omega }$ and $F_\omega $~.

Similarly, the real $B_{F^{+,T}_{\o v}} \otimes B_{F^{+,T}_{v}} 
$-bisemimodule $M_{F^{+,T}_{\o v}} \otimes M_{F^{+,T}_{v}} $ is the representation space of
$\GL_2(F^{,T}_{\o v}\times F^{+,+}_v)$ in such a way that:
\[ \Repsp(\GL_2( F^{+,T}_{\o v_\oplus} \times F^{+,T}_{v_\oplus}))
= M_{F^{+,T}_{\o v_\oplus}} \otimes M_{F^{+,T}_{v_\oplus}} 
=\txt\bigoplus\limits_{n,m_n} (M_{F^{,T}_{\o v_{n,m_n} }} \otimes M_{F^{+,T}_{v_{n,m_n }}} )\;.\]

\item At the  condition that $m^{(\omega _n)}=1$~, $\forall\ n$~, the  surjective morphism: $\gamma ^\delta \RTL:\o X\SRTL\to\partial \o X\SRTL$~, defining the boundary of the Borel-Serre compactification, corresponds to a covering of $\Repsp(\GL_2^{\rm (res)}(F^T_{\o \omega }\times F^T_\omega ))$ by 
$\Repsp(\GL_2(F^{+,T}_{\o v}\times F^{+,T}_v))$ in such a way that:
\Be
\item the representation space of the ``restricted'' algebraic bilinear semigroup\linebreak $\GL_2^{\rm (res)}(F^T_{\o\omega }\times F^T_\omega )$ be given by:
\[\Repsp (\GL_2^{\rm (res)}(F^T_{\o\omega }\times F^T_\omega ))=M^{\rm res}_{F^T_{\o\omega }}\otimes M^{\rm res}_{F^T_\omega }\]
homomorphic to $M^{\rm res}_{F^T_{\o\omega_\oplus }}\otimes M^{\rm res}_{F^T_{\omega _\oplus}} =\bigoplus_n(M_{F^T_{\omega _n}}\otimes M_{F^T_{\omega _n}})$~;
\item 
each complex subbisemimodule
\[ (M_{F^T_{\o \omega_n }} \otimes M_{F^T_{\omega_n}} )\;, \qquad 1\le n\le s\le\infty \;, \quad m_{\omega _n}=1\;, \]
be covered by the $m^{(n)}$ real subbisemimodules $\{
(M_{F^{+,T}_{\o v_{n,m_n} }} \otimes M_{F^{+,T}_{v_{n,m_n}}} )\}^{m^{(n)}}_{m_n=1}$~.
\Ee
\item
So, the complex ``bi''points of the boundary of $\o X\SRTL$ are in one-to-one correspondence with the real 
``bi''points of $\partial \o X\SRTL$~.\epr
\Ee}
\vskip 11pt

\subsection{Corollary}

{\em The surjective morphism
\[ \gamma ^\delta \RL: \quad \o X\SRL\to \partial\o X\SRL\]
is a morphism of inclusion in the sense that:
\Bean
\item the real lattice $\Lambda ^{(2)}_v$ (resp. $\Lambda ^{(2)}_{\o v}$~) is included into the corresponding restricted complex lattice  $\Lambda ^{(2)(\rm res)}_\omega =\sum_n\Lambda ^{(2)}_{\omega _n}$ (resp. $\Lambda ^{(2)(\rm res)}_{\o \omega} =\sum_n\Lambda ^{(2)}_{\o\omega _n}$~):
\[ \Lambda ^{(2)}_v \subset\Lambda ^{(2)(\rm res)}_{\omega } \qquad
\text{(resp.} \quad 
\Lambda ^{(2)}_{\o v}\subset\Lambda ^{(2)(\rm res)}_{\o\omega } \ ).\]

\item the real lattice $\Lambda ^{(2)}_v =\bigoplus\limits_{n,m_n}\Lambda ^{(2)}_{v_{n,m_n}}$ is commensurable with the restricted complex lattice $\Lambda ^{(2)(\rm res)}_{\omega } =\sum\limits_n
\Lambda ^{(2)}_{\omega_n } $ if and only if:
\[\Lambda ^{(2)}_{\omega_n } =\bigoplus\limits_{m_n}
\Lambda ^{(2)}_{v_{n,m_n} } \;, \quad 1\le n\le s\;.\]\Ee
}
\vskip 11pt

{\parindent=0pt
\bpr In correspondence with the covering of the representation space 
$M^{\rm res}_{F^T_{\o\omega }}\otimes M^{\rm res}_{F^T_\omega }$ of the restricted algebraic bilinear semigroup 
$\GL_2^{\rm (res)}(F^T_{\o\omega }\times F^T_\omega )$ by the representation space 
$M_{F^{+,T}_{\o v}}\otimes M_{F^{+,T}_v}$ of the real bilinear semigroup 
$\GL_2(F^{+,T}_{\o v}\times F^{+,T}_v)$ as given in proposition 2.5, we have that 
$\Lambda ^{(2)(\rm res)}_{\o\omega }\otimes \Lambda ^{(2)(\rm res)}_\omega $ is covered by 
$\Lambda ^{(2)}_{\o v}\otimes \Lambda ^{(2)}_v$~.

It then results that $\Lambda ^{(2)}_{\omega _n}=\bigoplus_{m_n} \Lambda ^{(2)}_{v_{n,m_n}}$ (resp.
$\Lambda ^{(2)}_{\o\omega _n}=\bigoplus_{m_n} \Lambda ^{(2)}_{\o, v_{n,m_n}}$~).
\epr
}
\vskip 11pt


\subsection{Bilinear parabolic subsemigroup}

Let $P_2(F^T_{[\omega _1]})$ (resp. $P_2(F^T_{[\o\omega _1]})$~) be a minimal parabolic (locally compact) subgroup over irreducible toroidal completions of $F_\omega $ (resp. $F_{\o\omega }$~), i.e. restricted to a global class residue degree $f_{\omega _1}=1$ (resp. $f_{\o\omega _1}=1$~). Let then
\[ P_2(F^T_{[\o\omega _1]})\times F^T_{[\omega _1]})\equiv
P_2(F^T_{[\o\omega _1]})\times P_2(F^T_{[\omega _1]})\]
denote a bilinear complex parabolic semigroup: it is the smallest (pseudo-)ramified normal bilinear subsemigroup of $\GL_2(F^T_{\o \omega }\times F^T_\omega )$ constituting a representation of the product $I_{F_{\o \omega _n}}\times I_{F_{\omega _n}}$~, $1\le n\le s$~, of the global inertia subgroups $I_{F_{\o \omega _n}}$ and $ I_{F_{\omega _n}}$~. 
Remark that $\GL_2(F^T_{\o \omega }\times F^T_\omega )$ acts by conjugation on
$P_2(F^T_{[\o\omega _1]}\times F^T_{[\omega _1]})$~.
Then, a double coset decomposition of $\GL_2(F^T_{R}\times F^T_L)$ leads to the following compactified bisemispace:
\[ \o S^{P_2}_{\GL_{2_{\ZZ_N}}}= P_2(F^T_{[\o\omega _1]}\times F^T_{[\omega _1]}) \setminus
\GL_2(F^T_{R}\times F^T_L ) \big/\GL_2(\ZZ\big/N\ \ZZ)^2\]
such that the modular conjugacy classes of $\GL_2(F^T_{\o\omega }\times F^T_\omega )$ (with respect to
$(I_{F_{\o \omega _n}}\times I_{F_{\omega _n}})$~) correspond to the cosets of $ P_2(F^T_{[\o\omega _1]}\times F^T_{[\omega _1]}) \setminus
\GL_2(F^T_{R}\times F^T_L ) $~.\vskip 11pt

The toroidal isomorphism of compactification 
$\gamma ^c_{R\times L}: X_{S_{R\times L}}\to \o X_{S_{R\times L}}$ with boundary 
$\partial \o X_{S_{R\times L}}$ leads to define the boundary
$ \partial \o S^{P_2}_{\GL_{2_{\ZZ_N}}}$ of $ \o S^{P_2}_{\GL_{2_{\ZZ_N}}}$ by:
\[ \partial \o S^{P_2}_{\GL_{2_{\ZZ_N}}}
= P_2(F^{+,T}_{[\o v_1]}\times F^{+,T}_{[v_1]})\setminus 
\GL_2(F^{+,T}_{R}\times F^{+,T}_L)\big/ \GL_2(\ZZ\big/N\ \ZZ)^2\]
where $P_2(F^{+,T}_{[\o v_1]}\times F^{+,T}_{[v_1]})$ is the bilinear real parabolic semigroup over the product $F^{+,T}_{[\o v_1]}\times F^{+,T}_{[v_1]}$ of irreducible toroidal completions of $F^+_{\o v}\times F^+_v$~.\\
$ \partial \o S^{P_2}_{\GL_{2_{\ZZ_N}}}$ is the equivalent of a Shimura (bisemi)variety following \cite{Pie2}.\vskip 11pt

It was shown in \cite{Pie1} that the representation space of $(P_2(F^{+,T}_{[v_1]}))$ (resp. $(P_2(F^{+,T}_{[\o v_1]}))$~) in a $P_2(F^{+,T}_{[v_1]})$-subsemimodule $M^I_{[v_1]}$
(resp. $P_2(F^{+,T}_{[\o v_1]})$-subsemimodule $M^I_{[\o v_1]}$~) is a \lr quantum of quantum field theory.\vskip 11pt

\subsection{\boldmath Analytic development of the $\GL_2(F^T_{\o\omega}\times F^T_\omega)$-bisemimodule $M_{F^T_{\o\omega}}\otimes M_{F^T_\omega}$}

It appears from section 2.3 that the complex bilinear algebraic semigroup $\GL_2(F_{\o\omega }\times F_\omega )$ has a representation space given by the $B_{F_{\o\omega }}\otimes B_{F_\omega }$-bisemimodule $M_{F_{\o\omega_\oplus}}\otimes M_{F_{\omega_\oplus} }$ which decomposes according to the equivalent representatives of the places of $F_{\o\omega }$ and $F_\omega $ as follows:
\[ M_{F_{\o\omega _{\oplus}}}\otimes M_{F_{\omega _{\oplus}}}=\bigoplus_n\bigoplus_{m_{\omega _n}}
( M_{F_{\o\omega_{n,m_{\omega _n}} }}\otimes M_{F_{\omega_{m,m_m}} })\]
where the set $\{ M_{F_{\omega_{n,m_{\omega _n}} }}\}_{F_{\omega _{n,m_{\omega _n}}}}$ (resp. $ \{M_{F_{\o\omega_{n,m_{\omega _n}}} }\}_{F_{\o\omega _{n,m_{\omega _n}}}}$~) of subsemimodules forms a tower of conjugacy class representatives of $T_2(F_\omega )\subset \GL_2(F_{\o\omega }\times F_\omega )$
(resp. $T^t_2(F_{\o\omega} )\subset \GL_2(F_{\o\omega }\times F_\omega )$~) characterized by increasing ranks, which are increasing integers modulo $N$ as developed in section 1.5.
\vskip 11pt

Remark also that the decomposition of $(M_{F_{\o\omega _\oplus}}\otimes M_{F_{\omega_\oplus} }$) into subbisemimodules also results from the action of the product of Hecke operators $(T_{q_R}\otimes T_{q_L})$ as it will be seen in the following.
 \vskip 11pt

So, the $B_{F_{\o\omega }}\otimes B_{F_\omega }$-bisemimodule $M_{F_{\o\omega }}\otimes M_{F_\omega }$ decomposes into a double symmetric tower of conjugacy class representatives corresponding each other respectively in the upper and in the lower half space.\vskip 11pt

The toroidal projective isomorphism of compactification sends the $B_{F_{\o\omega }}\otimes B_{F_\omega }$-bisemi\-module $M_{F_{\o\omega _\oplus}}\otimes M_{F_{\omega _\oplus}}$  into the
$B_{F^T_{\o\omega }}\otimes B_{F^T_\omega }$-bisemimodule
$M_{F^T_{\o \omega _\oplus}}\otimes M_{F^T_{\omega_\oplus} }$ which decomposes following the ``~$s$~'' places of
$F^T_{\o \omega }$ and $F^T_\omega $ according to:
\[ M_{F^T_{\o \omega_\oplus}}\otimes M_{F^T_{\omega_\oplus}}
= \txt\bigoplus\limits_{n=1}^s \bigoplus\limits_{m_{\omega _n}}\L( M_{F^T_{\o\omega_{n,m_{\omega _n}}}} \otimes 
M_{F^T_{\omega_{n,m_{\omega _n}}}} \R)\;.\]
With respect to proposition 2.5, we take into account the restricted $\GL_2(F^T_{\o\omega }\times F^T_\omega )$-bisemimodule
\[ M^{\rm res}_{F^T_{\o\omega_\oplus }}\otimes M^{\rm res}_{F^T_{\omega _\oplus}}= 
\txt\bigoplus\limits^s_{n=1} \L( M_{F^T_{\o\omega_n }}
\otimes M_{F^T_{\omega _n}}\R)\;, \quad m^{(\omega _n)}=1 \;,\]
where the \lr subsemimodule $M_{F^T_{\omega _n}}$ (resp. 
$M_{F^T_{\o\omega _n}}$~) is a \lr two-dimensional semitorus 
$T^2_L[n]$ (resp. $T^2_R[n]$~) localized in the upper (resp. lower) half plane.

\vskip 11pt

Indeed, let $\vec z=\sum\limits^2_{\alpha =1}z_\alpha \ \vec e_\alpha $ be a vector of $\cit ^1$
and, more precisely, of $F^T_\omega $~, 
and fix $z=\sum\limits^2_{\alpha =1}z_\alpha \ |\vec e_\alpha |$~. Then, every \lr 2-dimensional 
real semitorus has the analytic development:
\begin{align*}
T^2_L[n] &\simeq \lambda(2,n)\ e^{2\pi inz}\\[11pt]
\text{(resp.} \quad
T^2_R[n] &\simeq \lambda (2,n)\ e^{-2\pi inz}\; ),\end{align*}
where 
\Bi
\item $\lambda(2,n)$~, introduced in section 2.12,  can be considered as a Hecke character;
\item  $\lambda (2,n)\ e^{2\pi inz}$ corresponds to the product of two (semi)circles localized in perpendicular planes \cite{Pie2}.\Ei
\vskip 11pt

More concretely, we have that the (semi)torus
\[ T^2_L[n] \approx S^1_{d_1}[n]\times S^1_{d_2}[n]\]
is diffeomorphic to the product of two circles $S^1_{d_1}[n]$ and $S^1_{d_2}[n]$ localized in perpendicular planes.

This gives rise to the following decomposition:
\begin{align*}
T^2_L[n]
&\approx \lambda (2,n)\ e^{2\pi inz}\;, \quad z=x_{d_1}+i\ y_{d_2}\;, \quad x_{d_1}\in\rit\;, \quad y_{d_2}\in \rit\;,\\
&\simeq S^1_{d_1}[n]\times S^1_{d_2}[n]\\
&= r_{S^1_{d_1}}\ e^{2\pi inx_{d_1}} \times r_{S^1_{d_2}}\ e^{2\pi in(iy_{d_2})}
\;,\end{align*}
where $r_{S^1_{d_1}}$ and $r_{S^1_{d_2}}$ are radii verifying
\[ \lambda (2,n)\simeq r_{S^1_{d_1}}\times r_{S^1_{d_2}}\;,\]
in such a way that $ e^{2\pi in(iy_{d_2})}$~, localized in a plane perpendicular to
$e^{2\pi inx_{d_1}}$ and defined over $i\rit$~, is effectively the equation of a circle.

This is justified by the fact that a rotation of 90$^0$ of the circle
$S^1_{d_2}[n]$ over $i\rit$ transforms it into the circle
$S^1_{d_{2_\perp}}[n]$ over $\rit$ localized in the same plane as the circle
$S^1_{d_1}[n]$ according to:
\begin{align*}
{\rm rot}(90^0) : \quad S^1_{d_2}[n]
&\To S^1_{d_{2_\perp}}[n]\\
r_{S^1_{d_2}}\ e^{2\pi in(iy_{d_2})}\Big/ i\rit
&\To r_{S^1_{d_2}}\ e^{2\pi in(y_{d_2})}\Big/ \rit\;,\end{align*}
where:
\Bi
\item $S^1_{d_2}[n]= r_{S^1_{d_2}} [\cos(2\pi iny_{d_2}]+i\sin (2\pi iny_{d_2})]$~,
\item $S^1_{d_{2_\perp}}[n]= r_{S^1_{d_2}} [\cos(2\pi ny_{d_2}]+i\sin (2\pi ny_{d_2})]$~.
\Ei

\vskip 11pt

Then, $M^{\rm res}_{F^T_{\omega_\oplus}}$ (resp. $M^{\rm res}_{F^T_{\o\omega_\oplus}}$~) has the analytic development:
\begin{align*}
\EIS_L(2,n)
&\simeq \txt\bigoplus\limits_{n=1}^s  \lambda(2,n)\ e^{2\pi inz}\;, 
\quad s\le \infty\;, \\[11pt]
\text{(resp.}\quad
\EIS_R(2,n)
&\simeq \txt\bigoplus\limits_{n=1}^s  \lambda(2,n)\ e^{-2\pi inz}\;, 
\quad s\le \infty\;)\end{align*}
where $\EIS_L(2,n)$ (resp. $\EIS_R(2,n)$~) is the Fourier development of the equivalent of a normalized \lr Eisenstein series of weight $k=2$ restricted to the upper (resp. lower) half plane \cite{Pie2} and verifying $\lambda (2,n)=\sigma ^{\rm res}_{k-1}(n)\approx n\cdot N$ according to proposition 2.15 while classically
$\sigma _{k-1}(n)=\sum_{d\mid n}d^{k-1}$ \cite{Ser2}.\vskip 11pt

So, $\EIS_L(2,n)$ (resp. $\EIS_R(2,n)$~) corresponds to a \lr generalized cusp form of weight 
$k=2$ as it will be introduced in the following sections;
it can be decomposed into a tower of semitori $T^2_L[n]$ (resp. $T^2_R[n]$~), characterized by increasing ranks and representing analytically the conjugacy class representatives of
$T_2(F^T_\omega )$ (resp. $T^t_2(F^T_{\o\omega} )$~) $\subset \GL_2(F^T_{\o\omega }\times F^T_\omega )$~).\vskip 11pt

\subsection{\boldmath Analytic development of the $\GL_2(F^{+,T}_{\o v}\times F^{+,T}_v)$-bisemimodule 
$M_{F^{+,T}_{\o v}}\otimes M_{F^{+,T}_v}$}

According to proposition 2.5, the boundary $\partial \o X\SRL$ of  $\o X\SRL$ is given by
$\partial \o X\SRL=\linebreak \GL_2(F^{+,T}_{R}\times F^{+,T}_L)\big/\GL_2(\ZZ\big/N\ \ZZ)^2$~.  So, the toroidal compactification of the $B_{F^+_{\o v}}\otimes B_{F^+_v}$-bisemimodule 
$M_{F^+_{\o v}}\otimes M_{F^+_v}$ is given by the 
$B_{F^{+,T}_{\o v}}\otimes B_{F^{+,T}_v}$-bisemimodule
$M_{F^{+,T}_{\o v}}\otimes B_{F^{+,T}_v}$ constituting the representation space of the real algebraic bilinear semigroup\linebreak $\GL_2({F^{+,T}_{\o v}}\times {F^{+,T}_v})$~.\vskip 11pt

$M_{F^{+,T}_{\o v_\oplus}}\otimes M_{F^{+,T}_{v_\oplus}}$ can be developed as a direct sum of subbisemimodules 
following the ``~$s$~'' places of $F^{+,T}_{v}$ and $ F^{+,T}_{\o v}$ according to:
\[
M_{F^{+,T}_{\o v_\oplus}}\otimes M_{F^{+,T}_{v_\oplus}}
=\txt\bigoplus\limits_{n=1}^s\bigoplus\limits_{m_n} 
(M_{F^{+,T}_{\o v_{n,m_n}}}\otimes M_{F^{+,T}_{v_{n,m_n}}})\]
where $M_{F^{+,T}_{v_{n,m_n}}}$ (resp. $ M_{F^{+,T}_{\o v_{n,m_n}}}$~) is a \lr real 
one-dimensional semitorus localized in the upper (resp. lower) half plane.\vskip 11pt

So, $M_{F^{+,T}_{v}}$ (resp. $M_{F^{+,T}_{\o v}}$~) decomposes according to a tower of conjugacy
class representatives of $T_2(F^{+,T}_v) \subset \GL_2({F^{+,T}_{\o v}}\times {F^{+,T}_v})$
(resp. $T^t_2(F^{+,T}_{\o v})$~), characterized by increasing ranks, which are increasing integers modulo $N$~, and localized respectively in the upper (resp. lower) half space.
\vskip 11pt

As every \lr 1-dimensional real semitorus has the analytical development
\begin{align*}
T^1_L[n,m_n]
&\simeq r(1,n,m_n)\ e^{2\pi inx}\;, \quad x\in F^{+,T}_v\;, \\
\text{(resp.}\quad
T^1_R[n,m_n]
&\simeq r(1,n,m_n)\ e^{-2\pi inx}\;),\end{align*}
where $r(1,n,m_n)=( \lambda _+(n^2_N,m^2_N)-\lambda _-(n^2_N,m^2_N))\big/2$ according to proposition 2.15 and section 2.20,
$M_{F^{+,T}_{v}}$ (resp. $ M_{F^{+,T}_{\o v}}$~) will be developed analytically 
according to
\begin{align*}
\ELLIP_L(1,n,m_n)
&\simeq \txt\bigoplus\limits_{n=1}^s \bigoplus\limits_{m_n}r(1,n,m_n)\ e^{2\pi inx}\;, 
\quad s\le \infty\;, \\[11pt]
\text{(resp.}\quad
\ELLIP_R(1,n,m_n)
&\simeq \txt\bigoplus\limits_{n=1}^s \bigoplus\limits_{m_n}r(1,n,m_n)\ e^{-2\pi inx}\;, \quad s\le \infty\;)
\end{align*}
as it will be seen in section 2.18.\vskip 11pt

Consequently, the analytical development of the representation space of $T_2(F^{+,T}_v)$ (resp. $T^t_2(F^{+,T}_{\o v})$~)
$\subset \GL_2(F^{+,T}_{\o v}\times F^{+,T}_v)$ is given by the Fourier series $\ELLIP_L(1,n,m_n)$ (resp. $\ELLIP_R(1,n,m_n)$~) in such a way that each term of $\ELLIP_L(1,n,m_n)$ (resp.\linebreak $\ELLIP_R(1,n,m_n)$~) be the analytical representation of a conjugacy class representative of $T_2(F^{+,T}_v)$ (resp. $T^t_2(F^{+,T}_{\o v})$~).  By this way, the Fourier series receive an algebraic and geometric interpretation.
\vskip 11pt

\subsection[Irreducible representations entering into the two-dimensional bilinear Langlands correspondences]{Irreducible representations entering into the Langlands two-dimensional bilinear correspondences}  

\Bean \item According to section 2.8, a representation space $\Repsp(\GL_2(F^T_{\o\omega}\times F^T_\omega))$ of the general bilinear semigroup $GL_2(F^T_{\o \omega}\times F^T_\omega)$ is given by the bisemimodule $(M^{\rm res}_{F^{T}_{\o \omega_\oplus }}\otimes M^{\rm res}_{F^{T}_{\omega_\oplus} })$ whose analytic representation is given by $\EIS_R(2,n)\otimes \EIS_L(2,n)$~. As the product, right by left, of the global Weil semigroups $W_{F_{\o\omega}}\times W_{F_\omega}$ corresponds to the product, right by left, of the subgroups of Galois modular automorphisms of $\GL_2(F_{\o\omega}\times F_\omega)$ following section 1.7, we see that the sum of the products, right by left, of the equivalence classes of the irreducible 2-dimensional Weil-Deligne representation $\Irr\Rep^{(2)}(W_{F_{\o\omega}}\times W_{F_\omega})$ of the bilinear global Weil semigroup $(W_{F_{\o\omega}}\times W_{F_\omega})$ is given by
$\Repsp(\GL_2({F_{\o\omega}}\times {F_\omega}))$~. On the other hand, $\EIS_R(2,n)\otimes \EIS_L(2,n)$ constitutes the sum of the products, right by left, of the equivalence classes of the irreducible cuspidal representation $\Irr\Cusp(\GL_2(F^T_{\o\omega}\times F^T_\omega))$ of $\GL_2(F^T_{\o\omega}\times F^T_\omega)$~.\vskip 11pt

\item Similarly, the representation space $\Repsp(\GL_2(F^{+,T}_{\o v}\times F^{+,T}_v))$ of 
$GL_2(F^{+,T}_{\o v}\times F^{+,T}_v)$ is  the bisemimodule $(M_{F^{+,T}_{\o v }}\otimes M_{F^{+,T}_v })$ whose analytic representation is given by $\ELLIP_R(1,n,m_n)\otimes \ELLIP_L(1,n,m_n)$~.   Let $W_{F_{\o v}}\times W_{F_v}$ be the product, right by left, of the global Weil semigroups interpreted as the product, right by left, of the 
subgroups of Galois modular automorphisms of $\GL_2(F^+_{\o v}\times F^+_v)$~. Then, the sum of the products, right by left, of the equivalence classes of the representation $\Irr\Rep^{(1)}(W_{F_{\o v}}\times W_{F_v})$ of the bilinear global Weil group $(W_{F_{\o v}}\times W_{F_v})$ is given by $\Repsp(\GL_2({F^+_{\o v}}\times {F^+_v}))$~. On the other hand, $\ELLIP_R(1,n,m_n)\otimes \ELLIP_L(1,n,m_n)$ constitutes the sum of the products, right by left, of the equivalence classes of the irreducible cuspidal representation $\Irr\ELLIP( \GL_2(F^{+,T}_{\o v}\times F^{+,T}_v) )$
of $\GL_2(F^{+,T}_{\o v}\times F^{+,T}_v) $~.\Ee \vskip 11pt  So, we can state the proposition~:\vskip 11pt

\subsection{Proposition}  {\em \Bean \item Over the sum of the products $F_{\o\omega}\times F_\omega$ of the completions of the complex number semifields $ \dot{\widetilde F} _R$ and $\dot{\widetilde F} _L$~, there is a Langlands bilinear global correspondence given by the bijection:
\begin{align*}
 \sigma _{F_{\o\omega}\times F_\omega} : \quad
\Irr\Rep^{(2)}(W_{F_{\o\omega}}\times W_{F_\omega}) &\To
\Irr\Cusp(\GL_2({F^T_{\o\omega}\times F^T_\omega} ))\;;\\
\Repsp(\GL_2(F_{\o \omega }\times F_\omega ))  & \To
\EIS_R(2,n)\times \EIS_L(2,n)\;.\end{align*}
 \item Over the sum of the products $F_{\o v}\times F_v$ of the completions of the real number semifields $\dot{\widetilde F} ^+_R$ and $\dot{\widetilde F} ^+_L$~, there is a Langlands bilinear global correspondence given by the bijection:
\begin{align*}
 \sigma _{F_{\o v}\times F_v} : \quad
\Irr\Rep^{(1)}(W_{F_{\o v}}\times W_{F_v}) &\To
\Irr\ELLIP(\GL_2({F^{+,T}_{\o v}\times F^{+,T}_v} ))\;;\\
\Repsp(\GL_2(F^+_{\o v}\times F^+_v))  & \To
\ELLIP_R(1,n,m_n)\times \ELLIP_L(1,n,m_n)\;.\end{align*}
\Ee}\vskip 11pt

{\parindent=0pt
\bpr Taking into account the section 2.10, the Langlands global correspondence $\sigma _{F_{\o\omega }\times F_\omega }$ implies that, to the double symmetric tower of conjugacy class representatives of $\GL_2(F_{\o\omega }\times F_\omega )$ corresponds the double symmetric tower of the cuspidal representations of these conjugacy class representatives consisting in terms of the Fourier developments of the considered cusp forms.

More concretely, let $g^{(2)} _{R\times L} [n]$ be a conjugacy class representative of $\GL_2(F_{\o\omega }\times F_\omega )$~, i.e. a $B_{F_{\o\omega_n }}\otimes B_{F_{\omega _n}}$-subbisemimodule $ M_{F_{\o\omega _n}} \otimes
M_{F_{\omega _n}} \subset M_{F_{\o\omega }} \otimes M_{F_{\omega} }$ and let
$\eis_{R\times L}(2,n)
= \lambda (2,n)\cdot e^{-2\pi inz } \otimes 
\lambda (2,n)\cdot e^{2\pi inz }\in \EIS_R(2,n)\otimes \EIS_L(2,n)$ denote its analytic cuspidal
(toroidal) counterpart.

Then, we have the following set of bijections:
\begin{align*}
g^{(2)} _{R\times L}[1] \quad &  \begin{CD}@>\sim>>\end{CD} \quad \eis _{R\times L}(2,1)\;,\\
\vdots \qquad & \qquad \qquad \qquad \vdots \\
g^{(2)}_{R\times L}[n] \quad &  \begin{CD}@>\sim>>\end{CD} \quad \eis _{R\times L}(2,n)\;,\\
 \vdots \qquad & \qquad \qquad \qquad \vdots \\
g^{(2)} _{R\times L}[s] \quad &  \begin{CD}@>\sim>>\end{CD} \quad \eis _{R\times L}(2,s)\;,\\
\noalign{implying the bijection}
\Repsp(\GL_2(F_{\o\omega_\oplus }\times F_{\omega_\oplus} )) & \begin{CD}@>\sim>>\end{CD} \quad 
\EIS_{R\times L}(2,n) 
=\EIS_{R}(2,n)\times \EIS_{L}(2,n)\;,\end{align*}
where:
\Bi
\item $\Repsp(\GL_2(F_{\o\omega _\oplus}\times F_{\omega_\oplus} )) =\txt\bigoplus\limits^s_{n=1}  \Repsp(g^{(2)}_{R\times L}[n])$~;
\item $\EIS_{R\times L}(2,n) =\txt\bigoplus\limits^s_{n=1}  \eis_{R\times L}(2,n)$~.
\Ei

The second Langlands global correspondence $\sigma _{F_{\o v}\times F_v}$ can be handled similarly.\epr
\vskip 11pt}


Note that this way of introducing two-dimensional Langlands correspondences by global class field concepts corresponds to a new approach of this problem which was developed locally more particularly by G. Henniart \cite{Hen}, H. Carayol \cite{Car}, B. Conrad, F. Diamond and R. Taylor \cite{C-D-T}, S. Gelbart \cite{Ge} and M. Harris and R. Taylor \cite{H-T}.\vskip 11pt


\subsection{Representations of products of Hecke operators}

The ring of endomorphisms of the $GL_2(F_{\o \omega }\times F_\omega )$-bisemimodule 
$(M_{F_{\o \omega }}\otimes M_{F_\omega })$~, decomposing it into the set of subbisemimodules 
$(M_{F_{\o \omega _{n,m_{\omega _n}}}}\otimes M_{F_{\omega _{n,m_{\omega _n }}}})$ following the bisublattices 
$(\Lambda ^{(2)}_{\o \omega _{n,m_{\omega _n}}} \otimes
\Lambda ^{(2)}_{\omega _{n,m_{\omega _n}}})$~, is generated over $\ZZ\big/N\ \ZZ$ by the products 
$(T_{q_R}\otimes T_{q_L})$ of Hecke operators $T_{q_R}$ and $T_{q_L}$ for $q\nmid N$ and by the products $(U_{q_R}\otimes U_{q_L})$ of Hecke operators $U_{q_R}$ and $U_{q_L}$ for $q\mid N$ \cite{M-W}: it is noted $T_H(N)_R\otimes T_H(N)_L$~.\vskip 11pt

Note that the ring of endomorphisms of the real $\GL_2(F^+_{\o v}\times F^+_v)$-bisemimodule $(M_{F^+_{\o v}}\times M_{F^+_v})$ can also be defined classically by products, right by left, of Hecke operators, if we take into account the proposition 2.5 and corollary 2.6.
\vskip 11pt


The coset representative of $U_{q_L}$~, referring to the upper half plane, can be chosen to be upper triangular and given by the  integral matrix $\left(\begin{smallmatrix} 1 & b_N \\ 0 & q_N\end{smallmatrix}\right)$ of the congruence subgroup $\Gamma_L(N)$ in $\GL_2(\ZZ)$ such that:
\Bi
\item \quad $q_N= *+q\cdot N$~;

\quad $b_N= *+b\cdot N$~,

\quad where $*$ denotes an integer inferior to $N$~.

Indeed, $q_N=q\cdot N$ involves that $q\mid q_N$ and also that $q\mid N$ if $N=r q$~, $r\in\NN$~.

On the other hand, $q\nmid N$ is equivalent to the condition $q\nmid q_N$ since then,\linebreak $q_N=q\cdot N+s$~.

\item \quad $q$ is the cardinality of the infinite place $v_q$~;
\item \quad $b$ refers to the multiplicity of the Hecke sublattice of level $q$~.
\Ei
\vskip 11pt

Similarly, the coset representative of $U_{q_R}$~, referring to the lower half plane, can be chosen to be lower triangular and given by the integral matrix 
$\left(\begin{smallmatrix} 1 & 0 \\b_N &  q_N\end{smallmatrix}\right)$ of the congruence subgroup $\Gamma_R(N)$ in $\GL_2(\ZZ)$~.

For general $n$~, we would have respectively integral matrices
$\left(\begin{smallmatrix} a_N & b_N \\ 0 & d_N\end{smallmatrix}\right)$ and
$\left(\begin{smallmatrix}  a_N & 0 \\ b_N& d_N\end{smallmatrix}\right)$ of determinant $a_N\cdot d_N=*+n\cdot N$~.
\vskip 11pt

In fact, the congruence subgroup $\Gamma _L(N)$ denotes a general inverse image \cite{Hus} of the congruence subgroups $\Gamma (N)$~, $\Gamma _1(N)$ or $\Gamma _0(N)$ in $\SL_2(\ZZ\big/ N\ \ZZ)$ \cite{M-W}, \cite{La}, denoting the group of matrices $\L(\begin{smallmatrix} a & b\\ c & d\end{smallmatrix}\R)$ with the suitable congruences and acting as group(s) of symmetries of the hyperbolic upper half plane by the rule $z\to (a\ z+b)\big/(c\ z +d)$~, $z\in \CC$ being a point of order $N$~, $\Im (z)>0$~.

The congruence subgroup $\Gamma _R(N)=\Gamma _L(N)^t$ denotes the general inverse image of the group of matrices $\L(\begin{smallmatrix} a & c\\ b & d\end{smallmatrix}\R)$ with the suitable congruences acting as group of symmetries of the hyperbolic lower half plane by the rule $z^*\to (a\ z^*+b)\big/(c\ z^*+d)$~, $z^*\in \CC$~, $\Im (z^*)<0$~.
\vskip 11pt


Considering that the group of matrices
\[ u_2(b_N) = \BM 1&b_N \\ 0&1\EM \qquad \text{and} \qquad u_2(b_N)^t=\BM 1&0\\ b_N& 1\EM\;, \]
elements of the nilpotent group relative to the $v_q$-th infinite place, generates $\FF_q$~, the following coset representative
\[ g_2(q^2_N,b_N)=\left[
\BM 1&b_N\\ 0&1\EM
\BM 1&0\\ b_N&1\EM\right]
\BM 1&0\\ 0&q^2_N\EM\]
will be adopted for $U_{q_R}\otimes U_{q_L}$~, where $\left( \begin{smallmatrix} 1&0\\ 0&q^2_N\end{smallmatrix}\right)$ is the element of the split Cartan subgroup referring to a pseudo-ramified quadratic infinite place $v_{q^2}\equiv \o v_q\times v_q$~.
\vskip 11pt
$g_2(q^2_N,b_N) $ has a Gauss decomposition form in diagonal and unipotent parts.  The unipotent part  of $g_2(q^2_N,b_N)$ 
is $u_2(b_N)\cdot u_2(b_N)^t$ (and not  $u_2(b_N)^t\cdot u_2(b_N)$~), since we are dealing with (bi)linear functionals according to the Riesz lemma, and
corresponds to the element of the decomposition group associated with the split Cartan subgroup element $\alpha_{q^2}=\left(\begin{smallmatrix}1&0\\ 0&q^2_N\end{smallmatrix}\right)$~.  \vskip 11pt

\subsection{Hecke eigenvalues and decomposition group eigenvalues coincide}
The {\bf decomposition group\/} $D_{q_N^2}$ of the square of $q_N$ associated with the split Cartan subgroup element $\alpha _{q_N^2}=
\left( \begin{smallmatrix} 1&0\\ 0&q^2_N\end{smallmatrix}\right)$ is given by 
$\left\{D_{q_N^2;b_N}=u_2(b_N)\cdot u_2(b_N)^t\right\}_{b_N}$~.  Indeed, the semisimple form of $D_{q_N^2;b_N}$ is unimodular since 
$\det ( \lambda^+_{D_{q_N^2;b_N}} \cdot \lambda^-_{D_{q_N^2;b_N}} )=1$ where $\lambda^{\pm}_{D_{q_N^2;b_N}}$ are the two eigenvalues of $D_{q_N^2;b_N}$ which thus maps $q^2_N$ into itself.
\vskip 11pt

Let $D_{q_N^2;b_N}$ be the element of the decomposition group acting on the split Cartan subgroup element $\alpha _{q_N^2}$~.  Let $\lambda _+(q_N^2,b_N^2)$ and $\lambda _-(q_N^2,b_N^2)$ be the eigenvalues of $(D_{q_N^2;b_N^2}\cdot \alpha _{q_N^2})$ so that the determinant of the semisimple form of $(D_{q_N^2,b_N^2}\cdot \alpha _{q_N^2})$ be given by $\det (D_{q_N^2;b_N^2}\cdot \alpha_{q_N^2})_{ss}=\lambda _+(q_N^2,b_N^2)\cdot\lambda _-(q_N^2,b_N^2)$~.  Then, the decomposition group $D_{q_N^2;b_N^2}$ is such that:
\begin{align*}
 D_{q_N^2;b_N^2}:\quad \det (\alpha _{q_N^2})&\To \det (D_{q_N^2;b_N^2}\cdot \alpha_{q_N^2})_{ss}\\
q^2_N &\To q^2_N\end{align*}
where $(\ )_{ss}$ denotes the corresponding semisimple form.
\vskip 11pt

\subsection{\boldmath Definition: The Frobenius element $\Frob(q^2)$}
In the pseudo-unramified case, i.e. when $N=1$~, $D_{q^2,b^2}=\Frob(q^2)$ which gives:
\begin{align*}
 D_{q^2,b^2}: \quad \det (\alpha _{q^2})&\To \det (\alpha  _{q^2}\cdot\Frob (q^2))\\
q^2 &\To q^2\end{align*}
mapping the square of the global residue field class degree $f_{{v_q}}=q$ into itself.
\vskip 11pt

So, {\bf the Hecke eigenvalues and the Frobenius eigenvalues also coincide\/} in the pseudo-unramified case.

As a result of the coincidence of the Hecke eigenvalues with the decomposition group eigenvalues or with the Frobenius eigenvalues,
 we can then state the following proposition and corollary:
\vskip 11pt

\subsection{Proposition}
{\em There is an explicit irreducible semisimple (pseudo-)ramified representation, associated with a weight two cusp form
\[ \rho_{\lambda_{\pm}}: \quad \Gal ( \widetilde F^{+}_{\o v}/k)\times \Gal (\widetilde F^{+}_v/k)\To\GL_2(T_H(N)_R\otimes T_H(N)_L)\;,\]
 having eigenvalues:
\[\lambda_{\pm}(q_N^2,b_N^2) =\F{(1+b^2_N+q^2_N)\pm [(1+b^2_N+q^2_N)^2-4q^2_N]^{\half}}2\]
verifying
\begin{align*}
\tr \rho_{\lambda_{\pm}} &= 1+ b^2_N+q^2_N\;,\\
\det \rho_{\lambda_{\pm}} &= \lambda_+(q_N^2,b_N^2)\cdot \lambda_-(q^2_N,b^2_N)=q^2_N\;.
\end{align*}}
\vskip 11pt

\bpr Indeed, we have a representation $\rho_{\lambda_{\pm}}$ whose coset representatives are given by the matrices:
\[ g_2(q^2_N,b_N) = \left[
\BM 1 & b_N \\ 0&1\EM \BM 1 & 0 \\ b_N & 1\EM\right] \BM 1&0 \\ 0&q^2_N\EM
= \BM 1+b^2_N & b_Nq^2_N\\ b_N & q^2_N\EM\]
having eigenvalues $\lambda_{\pm} (q_N^2,b_N^2)$~.
\\
It is then easy to check that $\tr\rho_{\lambda_{\pm}}=1+b^2_N+q^2_N$ and that $\det\rho_{\lambda_{\pm}}=q^2_N$~.\epr
\vskip 11pt

\subsection{Corollary}
{\em Let $\rho^{nr}_{\lambda_{\pm}}$ be the corresponding irreducible semisimple (pseudo-)unramified representation, associated with a weight two cusp form:
\[\rho^{nr}_{\lambda_{\pm}}: \quad \Gal (\widetilde F^{+,nr}_{\o v}/k)\times \Gal (\widetilde F^{+,nr}_v/k)\To\GL_2(T_{H_R}\otimes T_{H_L})\;,\]
where
\[ T_{H_R}\otimes T_{H_L}=\left. T_H(N)_R\otimes T_H(N)_L\right|_{N=1}\;, \]
having eigenvalues:
\[\lambda^{nr}_{\pm}(q^2,b^2) =\F{(1+b^2+q^2)\pm [(1+b^2+q^2)^2-4q^2]^{\half}}2\]
verifying
\begin{align*}
\tr\rho^{nr}_{\lambda_{\pm}}(\Frob q^2) &= 1+ b^2+q^2\;,\\
\det \rho^{nr}_{\lambda_{\pm}}(\Frob q^2) &= \lambda^{nr}_+(q^2,b^2)\cdot \lambda^{nr}_-(q^2,b^2)=q^2\;.
\end{align*}}
Then, $\rho^{nr}_{\lambda_{\pm}}$ has a characteristic polynomial having the form \[X^2-\tr \rho^{nr}_{\lambda_{\pm}}(\Frob q^2)X + q^2=0\]
where $X$ is an indeterminate \cite{La-Tr}, \cite{Shi}, \cite{Ser1}, \cite{Clo}.
\vskip 11pt


\subsection{On the relevance of considering bialgebras for modular representations}

Assume classically that $f$ is a normalized eigenform associated with the congruence subgroup $\Gamma _1(N)$ in $SL_2(\ZZ)$ of weight $k\ge 2$~.  Then, $T_n$ is the Hecke operator verifying $T_nf=c(n,f)f$~, for each integer $n$~, where $c(n,f)$ is an algebraic integer.  Let $K_f$ be the number field generated over $\QQ$ by the $\{c(n,f)\}$ and $\theta $ its ring of integers.
\vskip 11pt
The normalized eigenforms $f$~, expanded in formal power series $f=\sum\limits_n a_nq^n$~, are non-zero cusp forms of the space $S(N)$ and are eigenvectors for all the $T_n$~, satisfying $a_1=1$ and $a_n=c(n,f)$~: so, Fourier coefficients of $f$ and eigenvalues of $T_n$ coincide.  The space $S(N)$ of the $f=\sum\limits_nc(n,f)q^n$ is thus a  [semi-]algebra \cite{Pie3} of cusp forms.
\vskip 11pt

If $H$ denotes the Poincare upper half plane in $\CC$~, the eigenforms $f$~, elements of the  [semi-]algebra $S(N)$~, 
are holomorphic in $H$ and defined in $\{\Im (z)>0\}$ with respect to the variable $z\in \CC$ 
of $q=e^{2\pi iz}$~.  The dual  [semi-]algebra of the [semi-]algebra $S(N)$~, relabelled $S_L(N)$~, is defined as the right [semi-]algebra 
$S_R(N)$ whose elements are normalized eigenforms $f_R$ associated with the congruence subgroup 
$\Gamma _1(N)^t$ (of the transposed matrices of $\Gamma_1(N)$~) and defined in the Poincare lower half plane 
$H^*$~.  (~``$L$~'' is for left and refers to the upper half plane while ``~$R$~'' is for right and refers to 
the lower half plane).
\vskip 11pt

These eigenforms $f_R$ of $S_R(N)$~, expanded in Fourier series $f_R=\sum\limits_na_{n,R}q^n_R=$\linebreak $\sum\limits_na_{n,R}e^{-2\pi inz}$~, are holomorphic in $H^*$ and defined in $\{\Im (z)<0\}$~.
The $f_R$ are eigenfunctions of Hecke operators $T_{n,R}$~, defined with respect to the congruence subgroup $\Gamma _1(N)^t$ and verifying
\[ T_{n,R}f_R=c(n,f_R)f_R \qquad \text{with} \quad a_{n,R}=c(n,f_R)\;.\]
\vskip 11pt
The bi[semi]algebra $S^e_L(N)$ associated with the [semi-]algebra $S_L(N)$ is given by $S_L^e(N)=S_R(N)\otimes_\theta S_L(N)$ and is of special importance owing to the natural homomorphism $\psi :S_L^e(N)\to \End_\theta (S_L(N))$ where $\End_\theta (S_L(N))$ is isomorphic to the Hecke algebra $\Hs_L$ generated by the Hecke operators $T_{n,L}(N)$~\cite{F-D}.
\vskip 11pt
Considering the isomorphism
\[\omega :\quad \End_\theta (S_R(N))\otimes \End_\theta (S_L(N))\To \End_\theta (S_R(N)\otimes S_L(N))\]
where $\End_\theta (S_R(N))\otimes \End_\theta (S_L(N))$ is isomorphic to the product $\Hs_R\otimes \Hs_L$ of the Hecke algebras $\Hs_R$ and $\Hs_L$~, generated by the Hecke operators $T_{n,R}$ and $T_{n,L}$~, while $\End_\theta (S_R(N)\otimes S_L(N))$ is isomorphic to the Hecke bialgebra $\Hs_{R\otimes L}$~, generated by tensor products $T_{n,R}\otimes T_{n,L}$ of Hecke operators $T_{n,R}$ and $T_{n,L}$ acting on $S^e_L(N)$~, we are led to the following commutative diagram:
\[ \begin{array}{ccc}
S^e_L(N) & \To & \End_\theta (S_L(N))\\
&\searrow & \Big\uparrow\\
& & \End_\theta (S^e_L(N))\end{array}\]
It then becomes clear that $\End_\theta (S_L(N))$ will be worked out successfully by taking into account the endomorphisms $\End_\theta (S^e_L(N))$ of the bialgebra $S^e_L(N)$~, i.e. by considering the Hecke bialgebra $\Hs_{R\otimes L}$ of tensor products $T_{n,R}\otimes T_{n,L}$ of right Hecke operators $T_{n,R}$ by left Hecke operators $T_{n,L}$ acting on tensor products of cusp forms $f_R\otimes f_L=\sum\limits_n a_{n,R}q^n_R\otimes \sum\limits_n a_{n,L}q^n_L\in S_R(N)\otimes S_L(N)$~.
\vskip 11pt

$f_R\otimes_D f_L$ is called a modular biform and is naturally defined in the present context on the complex algebraic bilinear semigroup $\GL_2(F^T_{\o\omega }\times F^T_\omega )$~.

\noindent $f_R\otimes _D f_L\in C^\infty(\GL_2(F^T_{R }\times F^T_L )\big/\GL_2((\ZZ\big/N\cdot \ZZ)^2$~.
\vskip 11pt

\subsection{The space of global elliptic semimodules}  
According to section 2.7, the compactified bisemispace
\[ \o S^{P_2}_{\GL_{2_{\ZZ_N}}}=P_2( F^T_{[\o \omega _1]}\times F^T_{[\omega _1]})
\setminus \GL_2( F^T_{R }\times F^T_{L}) \Big/ \GL_2(\ZZ\big/N\cdot \ZZ)^2\]
decomposes into pairs of two-dimensional semitori following the ``modular'' conjugacy classes of 
$\GL_2( F^T_{\o \omega }\times F^T_{\omega})$ which are in one-to-one correspondence with the places of $ F^T_{\o \omega }$ or $F^T_{\omega}$ counted with their multiplicities.\vskip 11pt

Remark that there is a one-to-one correspondence between $\o S^{P_2}_{\GL_{2_{\ZZ_N}\mid L}}$~, restricted to the upper half space (left case), and given by the semispace
\[ \o S^{P_2}_{\GL_{2_{\ZZ_N}\mid L}}
= P_2(F^T_{[\omega_1]})\setminus T_2(F^T_L)\Big/ T_2(\ZZ\big/N\cdot \ZZ)\;, \]
and the jacobian
\[ J(N)_\CC=H^0(X(N)_\CC,\Omega)\Big/ H_1(X(N)_\CC,\ZZ)\]
of the Riemann surface $X(N)_\CC$ corresponding to the group $\Gamma(N)$ associated with$\Gamma_L(N)$ (see section 2.12).

Indeed, the subgroup $T_2(\ZZ\big/N\ \ZZ)$ is a representation of Hecke lattice operators cutting $T_2(F^T_L)$ following its modular conjugacy classes which are in one-to-one correspondence with the sublattices or periods of the jacobian  $J(N)_\CC$ \cite{Hin}, \cite{Bos}.  In order to have a natural model over $\QQ$ of $X(N)_\CC$ leading to $J(N)_\QQ={\rm Pic}(X(N)_\QQ)$~, 
 the boundary
\[\partial \o S^{P_2}_{\GL_{2_{\ZZ_N}}}=P_2( F^{+,T}_{[\o v _1]}\times F^{,T}_{[v _1]})
\setminus \GL_2( F^{+,T}_{R}\times F^{+,T}_{L}) \Big/ \GL_2(\ZZ\big/N\ \ZZ)^2\]
of $\o S^{P_2}_{\GL_{2_{\ZZ_N}}}$ is assumed to be decomposed into pairs of packets of 
one-dimensional semitori following the biplaces of $F^{+,T}_{[\o v]}\times F^{+,T}_{[v]}$ such that each packet of $1D$-semitori covers the corresponding complex representative of the considered conjugacy class of $T_2(F^{T}_\omega )$ or of $T^t_2(F^{T}_{\o \omega })$~.\vskip 11pt

Let the number of infinite complex places be equal to 
the number of real infinite places.  Then, every $2D$-semitorus 
$T^2_L[n,m_n=1]\in \o S^{P_2}_{\GL_{2_{\ZZ_{N\mid L}}}}$ 
(resp. $T^2_R[n,m_n=1]\in \o S^{P_2}_{\GL_{2_{\ZZ_{N\mid R}}}}$~) 
can be decomposed into a packet of $1D$-semitori $T^1_L[n,m_n]$ (resp. $T^1_R[n,m_n]$~) 
referring to the $n$-th place of $F^{+,T}_{[v]}$ (resp. $F^{+,T}_{[\o v]}$~) 
in one-to-one correspondence with the $n$-th modular conjugacy class of $\GL_2(F^{+,T}_{[\o v]}\times F^{+,T}_{[v]})$ corresponding to the $n$-th coset of the boundary $\partial\o S^{P_2}_{\GL_{2_{\ZZ_N}}}$ of $\o S^{P_2}_{\GL_{2_{\ZZ_N}}}$~.\vskip 11pt

So, every $2D$-semitorus of $\o S^{P_2}_{\GL_{2_{\ZZ_N}}}$ can be decomposed into a packet $\{T^1 _L[n,m_n]\}_{m_n}$ of $1D$-semitori isomorphic to a set of affine curves of
the $T_2(F^+_{v_n})$-subsemimodule $M_{F^+_{v_n}}$ and characterized each one by an integer ``~$n$~'' with respect to the representation of the smallest pseudo-ramified normal subsemigroup of class ``~$1$~'' of $T_2(F^+_v)$ which is a $P_2(F^+_{[v_1]})$-subsemimodule $M^I_{[v_1]}$ also called a left quantum. The integer $n$ corresponds to the global class residue degree
\[ [F^{+,nr}_{v_n}:k]=f_{v_n}=n\]
in such a way that a pseudo-ramified one-dimensional semitorus of class $n$ has a rank given by $[F^+_{v_n}:k]=n\cdot N$ where $N$ is the rank of the quantum $M^I_{[v_1]}$~.\vskip 11pt

The set of continuous complex-valued functions over the conjugacy class representatives of the
$T_2(F^{+,T}_{v})$-semimodule $M_{F^{+,T}_{v}}$ (resp. $T^t_2(F^{+,T}_{\o v})$-semimodule 
$M_{F^{+,T}_{\o v}}$~) is a semisheaf of rings whose set of sections is noted $A_L\equiv \Gamma (M_{F^{+,T}_{v}})$ (resp. $A_R\equiv \Gamma (M_{F^{+,T}_{\o v}})$~).\vskip 11pt

This set of sections $A_L$ (resp. $A_R$~) has the structure of a semiring $C(M^I_{[v_1]})$ (resp. $C(M^I_{[\o v_1]})$~) of continuous complex valued functions over \lr quanta $M^I_{[v^1_n]}$ (resp. $M^I_{[\o v_n^1]}$~) due to the injective maps:
\begin{align*}
 m ^I_L(n) : \quad & M^I_{[v_1]} \To M_{F^{+,T}_{v_n}}\;, \quad 1\le n\le s\;, \\
\text{(resp.} \quad 
 m ^I_R(n) : \quad & M^I_{[\o v_1]} \To M_{F^{+,T}_{\o v_n}} \;, \quad 1\le n\le s\;),
\end{align*} where $M_{F^{+,T}_{v_n}}$ (resp. $M_{F^{+,T}_{\o v_n}}$~) is a left- (resp. right)-$T_2(F^{+,T}_{v_n})$-subsemimodule $\in M_{F^{+,T}_{v}}$  (resp. $T^t_2(F^{+,T}_{\o v_n})$-subsemimodule $\in M_{F^{+,T}_{\o v}}$~).\vskip 11pt

For every section $ (s_L)_{n,b} \in A_L$ (resp. $ (s_R)_{n,b}\in A_R$~), let $\End_{F^{+,T}_{v}}(A_L)$ (resp.\linebreak $\End_{F^{+,T}_{\o v}}(A_R)$~) be the Frobenius endomorphism of $A_L$ (resp. $A_R$~) and let 
\begin{align*}
q\Big/\QQ_L &\To q^n\big/\QQ_L\in \End_{F^{+,T}_{v_n}} (s_L)_{n,b} \\
\text{(resp.} \quad q\Big/\QQ_R &\To q^n\big/\QQ_R\in 
\End_{F^{+,T}_{\o v_n}} (s_R)_{n,b} \;)
\end{align*}
be the Frobenius substitution where $q^n\big/\QQ_L =e^{2\pi inx}$ (resp. $q^n\big/\QQ_R =e^{-2\pi inx}$~); $x\in F^{+,T}_v$ being a point of order $N$~. 
\vskip 11pt

A global elliptic \lr-$A_L$-semimodule $\phi_L(s_L)$ (resp. $A_R$-semimodule $\phi_R(s_R)$~) over $\QQ$ is a ring homomorphism:
\begin{align*}
\phi_L : \quad A_L &\To \End_{F^{+,T}_{v}}(A_L)\\
\text{(resp.} \quad
\phi_R : \quad A_R &\To \End_{F^{+,T}_{\o v}}(A_R)\;)\end{align*}
given, for the sections $s_L\in A_L$ (resp. $s_R$ in $A_R$~), by:
\begin{align*}
\phi _L(s_L) &= \sum_{n=1}^s \sum_b \phi (s_L)_{n,b}\ q^n\Big/ \QQ_L\;, 
\qquad s\le \infty\;,
\\
\text{(resp.} \quad
\phi _R(s_R) &= \sum_{n=1}^s \sum_b \phi (s_R)_{n,b}\ q^n\Big/ \QQ_R\;, 
\qquad s\le \infty\;), \end{align*} where
 \Bi 
\item $\sum\limits_n$ runs over the ~classes ``~$n$~'' of the set of pairs of $1D$-semitori isomorphic to the set of pairs of the affine curves of $\GL_2(F^+_{R}\times F^+_L)\big/\GL_2(\ZZ\big/N\ \ZZ)^2$~;

\item $\sum\limits_b$ runs over the representatives of the pairs of the $1D$-semitori of class ``~$n$~''. These representatives correspond to the number of (semi)orbits or ideals of the decomposition group $D_{n^2_N;b_N} $ according to definition 2.13;
\item $\phi(s_L)_{n,b}$ (resp. $\phi(s_R)_{n,b}$~) is a coefficient of the Fourier series $\phi_L(s_L)$ (resp. $\phi_R(s_R)$~). \Ei\vskip 11pt

Remark that the name of global elliptic \lr $A_L$-semimodule $\phi _L(s_L)$ (resp. $A_R$-semimodule $\phi _R(s_R)$~), also noted $\ELLIP_L(1,n,m_n)$ (resp. $\ELLIP_R(1,n,m_n)$~) in section 2.11 for these (truncated) Fourier series is justified by the facts that: 
\Bean\item a subset of $1D$-semitori of level $p$~, where $p$ is a prime number, can be identified with the orbit space of an elliptic curve $E(\FF_p)$ under the action of the decomposition group $D_{p^2_N}$ as it will be seen in chapter 3; \item these global elliptic (semi)modules have some analogy of construction with respect to the local elliptic modules (over function fields) introduced by V. Drinfeld \cite{Drin} (see also \cite{And}). \Ee\vskip 11pt

\subsection{Proposition}  {\em Let $S_L(f)$ (resp. $S_R(f)$~) be the left [semi-]algebra (resp. right [semi-]algebra) of modular forms $f_L(z) = \sum\limits_n a_{n,L}\ q^n_L$ (resp. $f_R(z) = \sum\limits_n a_{n,R}\ q^n_R$~) $q_L=e^{2\pi iz}$
(resp. $q_R=e^{-2\pi iz}$~), $z\in\cit$~,
\Bi \item being normalized eigenforms of Hecke operators associated with the congruence subgroup $\Gamma _L(N)$ (resp. $\Gamma _R(N)$~) introduced in section 2.12; \item characterized by a weight two and a level $N$~. \Ei 
On the other hand, let $S_L(\phi )$ (resp. $S_R(\phi )$~) be the left [semi-]algebra (resp. right [semi-]algebra) of global elliptic $A_L$-semimodules  $\phi _L(s_L) $ (resp. $\phi _R(s_R) $~) in the sense that $f_L(z)\simeq \phi _L(s_L)$ (resp. $f_R(z)\simeq \phi _R(s_R)$~). Then, we have the following inclusions of left [semi-]algebras (resp. right [semi-]algebras):
\begin{align*}
S_L(\phi) &\hookrightarrow S_L(f)\\
\text{(resp.} \quad 
S_R(\phi) &\hookrightarrow S_R(f)\;).\end{align*} }\vskip 11pt

\paragraph{Proof:} according to section 2.8, the (truncated) Eisenstein series $\EIS_L(2,n,m_n=1)$
(resp. $\EIS_R(2,r,m_n=1)$~), which is a (pseudo)modular form $f_L(z)$ (resp. $f_R(z)$~) \cite{Ma3}, constitutes the analytic development of the $T_2(F^T_\omega )$-semimodule $M^{\rm res}_{F^T_\omega }$
(resp. $T^t_2(F^T_{\o\omega} )$-semimodule $M^{\rm res}_{F^T_{\o\omega} }$~) according to the global Langlands correspondence a) of proposition 2.11. On the other hand, the global elliptic $A_L$-semimodule $\phi_L(s_L)$ (resp. $A_R$-semimodule $\phi_R(s_R)$~) constitutes the analytic development of the $T_2(F^{+,T}_v)$-semimodule $M_{F^{+,T}_v}$ (resp.
$T_2(F^{+,T}_{\o v})$-semimodule $M_{F^{+,T}_{\o v}}$~) according to section 2.10 and the global correspondence of Langlands b) of proposition 2.11.\\ 
Then, if we assume that (conditions of proposition 2.5):
\Be
\item the complex places of $F^T_\omega$ have simple multiplicities given by $m_n=1$~;
\item the number of complex places of $F^T_\omega$ is equal to the number of real places of $F^{+,T}_v$~;
\item $\partial \o S^{P_2}_{\GL_{2_{\ZZ_N}}}\hookrightarrow \o S^{P_2}_{\GL_{2_{\ZZ_N}}}$ (see section 2.18), implying that complex irreducible completions coincide with real irreducible completions,
\Ee
 it is then clear that
\Bi \item $f_L(z)\simeq \phi_L(s_L)$ \qquad (resp. \quad $f_R(z)\simeq \phi_R(s_R)$~); \item $S_L(\phi)\hookrightarrow S_L(f)$ \qquad (resp. \quad $S_R(\phi)\hookrightarrow S_R(f)$~).\epr\Ei\vskip 11pt  


\subsection{Supercuspidal representations in terms of global elliptic\protect\linebreak bisemimodules}
Let $\otimes_D$ denote a diagonal tensor product, i.e. a tensor product whose only diagonal terms with respect to a basis $\{e_{n,b}\otimes e_{n,b}\}$ are different from zero.
Let $\lambda _+(q_N^2,b_N^2)$ and $\lambda _-(q_N^2,b_N^2)$ be the eigenvalues of the products of Hecke operators having as coset representatives $g_2(q^2_N,b_N)$~.

Assume the existence of a global  elliptic $A_{R-L}$-bisemimodule $\phi _R(s_R)\otimes_D\phi _L(s_L)
=$\linebreak $\sum\limits_{n,b}\phi (s_R)_{n,b}q^n_{/\QQ_R}\otimes_D \sum\limits_{n,b}\phi (s_L)q^n_{/\QQ_L}$ whose coefficients $\phi (s_{L})_{n,b}$ (resp. $\phi (s_{R})_{n,b}$~)are given by $\phi (s_{L})_{n,b}=(\lambda _{+}(n_N^2,b_N^2))$ (resp. $\phi (s_{R})_{n,b}=(\lambda _{-}(n_N^2,b_N^2))$~) and which verifies $f_R(z)\otimes _Df_L(z)\simeq \phi _R(s_R)\otimes_D\phi _L(s_L)$~, i.e. the ``inclusion'', in the sense of 2.19, of the elliptic $A_{R-L}$-bisemimodule into a diagonal tensor product of weight two cusp forms.
Consider the isomorphism:
\begin{align*}
 i_R\otimes i_L:\quad &\phi _R(s_R)\otimes_D \phi _L(s_L)=\sum_{n,b}(\lambda _+(n_N^2,b_N^2))q^n_{/\QQ_R}\otimes_D \sum_{n,b}(\lambda_-(n_N^2,b_N^2))q^n_{/\QQ_L} \\
& \qquad \quad\To \widetilde\phi _R(s_R)\otimes_D \widetilde\phi _L(s_L)=\sum_{n,b}(r(n_N^2,b_N^2))q^n_{/\QQ_R}\otimes_D \sum_{n,b}(r(n_N^2,b_N^2))q^n_{/\QQ_L}
\end{align*}
sending the eigenvalues $(\lambda _{\pm}(n_N^2,b_N^2))$ to $(r(n_N^2,b_N^2))$ defined by $(r(n_N^2,b_N^2))=\linebreak ((\lambda _+(n_N^2,b_N^2))+(\lambda _-(n_N^2,b_N^2))/2 = \tr \rho _{\lambda _\pm}(n^2_N,b^2_N)/2$ given in proposition 2.15.

Then, $ \widetilde \phi _R(s_R)\otimes_D \widetilde  \phi _L(s_L)$ decomposes into a sum of (tensor) products of irreducible semitoric curves localized in the upper and in the lower half space corresponding to each other by pairs of same class $n$ and same value of $b$ in such a way that each pair of semitoric curves be characterized by a radius $(r(n_N^2,b_N^2))$ and a center at the origin.  The isomorphism $i_R\otimes i_L$ has then to be interpreted as a translation of the toric curves since $(\lambda _{\pm}(n_N^2,b_N^2))$ may be viewed as a pair $\{{\rm cent}(n_N^2,b^2)),(r(n_N^2,b_N^2))\}$ where $({\rm cent}(n_N^2,b_N^2))= 
((\lambda _+(n_N^2,b_N^2))-(\lambda _-(n_N^2,b_N^2))/2$ is the image of the translated center under $i_{R,L}$~.
It then results that the eigenvalues $\lambda _+(n _N^2,b_N^2)$ and $\lambda _-(n_N^2,b^2_N)$ are equivalent.
\vskip 11pt

A diagonal tensor product of weight two cusp forms thus has for representation\linebreak $\widetilde\phi_R(s_R)\otimes_D\widetilde\phi (s_L)$ which will be rewritten according to:
\[ \widetilde\phi_R(s_R)\otimes_D\widetilde\phi (s_L) \;
{\mathrel\mid\joinrel\mathrel\Rightarrow}\;\Pi_R\otimes_D \Pi_L=\sum_n (m_R(n)\Pi_R(n)) \otimes_D (m_L(n)\Pi_L(n)) \]
where $m_{L}(n)$ (resp. $m_{R}(n)$~) denotes the multiplicity of the irreducible curve $\Pi_{L}(n)$ (resp. $\Pi_{R}(n)$~) of global level $n$~, i.e. associated with the global \lr place $v_n$ (resp. $\o v_n$~).
\vskip 11pt

\subsection{Proposition}

 {\em Let $GL_2(F^+_{\o v}\times F^+_v)$ be the general bilinear algebraic semigroup  over $F^+_{\o v}\times F^+_v$~.  Then, $\GL_2(F^+_{\o v}\times F^+_v)$ has for irreducible supercuspidal representation the Grothendieck group $\groth (\GL_2(F^+_{\o v}\times F^+_v))$ defined by
\[\groth (GL_2(F^+_{\o v}\times F^+_v))=\Pi_R\otimes_D \Pi_L=\sum_n (m_R(n)\Pi_R(n)) \otimes_D (m_L(n)\Pi_L(n)) \;.\]}
\vskip 11pt
\bpr Indeed $\Pi_R\otimes_D\Pi_L$ constitutes an irreducible supercuspidal representation of\linebreak $GL_2(F^+_{\o v}\times F^+_v)$ such that the sum of the representations of the conjugacy classes of $\GL_2(F^+_{\o v}\times F^+_v)$ are in bijection with the sum of the irreducible supercuspidal subrepresentations of $\Pi_R\otimes_D\Pi_L$ over the quadratic places $(\o v_n\times v_n)$~.\epr
\vskip 11pt

\subsection{Injective nilpotent morphisms}

Assume that the bi[semi]algebra $S_R(\phi)\otimes_D S_L(\phi)$ of elliptic $A_{R-L}$-bisemimodules is included into the bi[semi]algebra $S_R(f)\otimes_D S_L(f)$ of weight two cusp forms.  

Let $\Pi ^0_{R\times L}:GL_2(F^+_{\o v}\times F^+_v)\to \phi^0 _R(s_R)\otimes_D \phi^0 _L(s_L)$ denote a simple supercuspidal representation of elliptic type of $GL_2(F^+_{\o v}\times F^+_v)$ such that
\[ \phi ^0_R(s_R)\otimes_D \phi ^0_L(s_L)=\sum_n(\lambda (n_N^2,b_N^2=0)q^n_{/\QQ_R}\otimes_D \sum_n( \lambda (n_N^2,b_N^2=0)) q^n_{/\QQ_L}\]
be defined with $\lambda (n_N^2,b_N^2=0)=n_N=n\cdot N$ as the square root of the eigenvalue $n^2_N$ of the coset representative
\begin{align*}
g_2(n^2_N,0) &= \left[\BM 1&0\\0&1\EM \BM 1&0\\ 0&1\EM \right] \BM 1&0\\ 0&n^2_N\EM\\[11pt]
&= I\cdot \alpha _{n_N^2}\;.\end{align*}
Suppose that $w^0_{v_{n_{R\times L}}} 
$ denotes the action of the decomposition group $D_{n^2_N}$ on the irreducible supercuspidal subrepresentation $\Pi_R(n)\otimes_D \Pi_L(n)$~.
\vskip 11pt

$\omega ^0_{v_{n_{R\times L}}}$ then corresponds to the injective (nilpotent) morphism:
\[ \omega ^0_{v_{n_{R\times L}}}: \quad
\Pi_R(n)\otimes _D\Pi_L(n) \To m_R(n)\ \Pi_R(n)\otimes_D m_L(n)\ \Pi_L(n)\]
which generates the multiplicity of the supercuspidal subrepresentation
\[ \Pi_R(n)\otimes_D \Pi_L(n)\equiv \lambda (n^2_N,b^2_N=0)\cdot q^n\big/Q_R \otimes_D
\lambda (n^2_N,b^2_N=0)\cdot q^n\big/Q_L\;.\]\vskip 11pt

We thus have a family $\omega ^0_{v_{{R\times L}}}$ of injective morphisms at all quadratic places $v_{R\times L}\equiv \o v\times v$ such that:
\[\Pi^0_{R\times L}(\omega ^0_{v_{{R\times L}}}): \quad \GL_2(F^+_{\o v}\times F^+_v) \To \phi_R(s_R)\otimes_D\phi_L(s_L)\]
where $\phi_R(s_R)\otimes_D \phi_L(s_L)$ is the global elliptic $A_{R-L}$-bisemimodule:
\[\phi_R(s_R)\otimes_D \phi_L(s_L)= \sum_{n,b}\lambda_+(n^2_N,b^2_N)\ q^n_{/
\QQ_R}
\otimes_D \sum_{n,b}\lambda_-(n^2_N,b^2_N)\ q^n_{/ \QQ_L}\;.\]\vskip 11pt

\subsection{Proposition}  {\em There exists a family $\omega ^0_{v_{{R\times L}}}$ of injective morphisms associated with the 
action of the set $\{D_{n^2_N}\}_{n^2}$ of decomposition groups such that:
\[ \omega ^0_{v_{{R\times L}}}:\quad \Pi^0_{R\times L}\To \Pi^0_{R\times L}(\omega^0_{v_{R\times L}})\]
implies the commutative diagram:

\setlength{\unitlength}{1mm}
\centerline{%
\begin{picture}(115,30)(0,0)
\put(29,8){$\mbox{}^{\Pi^0_{R\times L}(\omega^0_{v_{R\times L}})}$}
\put(45,15){\vector(3,-1){18}}
\put(0,0){\makebox(115,30)[b]{$\begin{array}{ccc}
\GL_2(F^+_{\o v}\times F^+_v) &
\begin{CD} @>{\quad\Pi^0_{R\times L}\quad}>> \end{CD}&
\phi^0_R(s_R)\otimes _D\phi^0_L(s_L) \\
&& \Big\downarrow \mbox{}_{\omega ^0_{v_{R\times L}}} \\
&&
\phi_R(s_R)\otimes_D\phi_L(s_L) \end{array}$}}
\end{picture}}
}
\vskip 11pt
 \subsection{Corollary}  {\em The family of injective morphisms $\omega^0_{v_{R\times L}}$ corresponds to the action of 
the product $\Ws_R\times \Ws_L$ of the Weyl groups acting on the supercuspidal representation of 
elliptic type $\phi ^0_R(s_R)\otimes _D\phi ^0_L(s_L)$~.
}\vskip 11pt

\paragraph{Proof:} indeed, $\phi ^0_R(s_R)\otimes _D\phi ^0_L(s_L)$ consists of the sum of products of pairs of maximal semitori $T^1_R[n]\times T^1_L[n]$ at all quadratic places $v_{R\times L}\equiv \o v\times v$~.\\
And, on the other hand, the action of the set $\{D_{n^2_N}\}_{n^2}$ of the decomposition groups, associated with the family of injective morphisms $\omega ^0_{v_{R\times L}}$~, corresponds to the action of the product $\Ws_R\times \Ws_L$ of the Weyl groups since, for every product $\Ws_{\o v_n}\times \Ws_{v_n}\in\Ws_R\times \Ws_L$ of Weyl subgroups restricted to the quadratic place $\o v_n\times v_n$~, we have:
\[\Ws_{\o v_n}\times \Ws_{v_n}=\F{N_R}{T^1_R[n]}\times \F{N_L}{T^1_L[n]}\]
where $N_L$ (resp. $N_R$~) is the normalizer of $T^1_L[n]$ (resp. $T^1_R[n]$~) in $\Pi_L(n)$ (resp. $\Pi_R(n)$~).\epr\vskip 11pt

\section{Applications: The treatment of some conjectures}

Chapter 2 deals with the representation of a complex algebraic bilinear semigroup $\GL_2(F_{\o\omega }\times F_\omega )$ into a $B_{F_{\o\omega }}\otimes B_{F_\omega }$-bisemimodule 
$M_{F_{\o\omega _\oplus}}\otimes M_{F_{\omega _\oplus}}$ decomposing into  direct sum of bisubsemimodules
$M_{F_{\o\omega_{n,m_{\omega _n}} }}\otimes M_{F_{\omega_{n,m_{\omega _n}} }}$ characterized by increasing ranks.
\vskip 11pt

A toroidal compactification of the Borel-Serre type of 
$M_{F_{\o\omega_\oplus }}\otimes M_{F_{\omega_\oplus} }$~:
\Be
\item maps each subsemimodule $M _{F_{\omega _{n,m_{\omega _n}}}}$ (resp. $M _{F_{\o\omega _{n,m_{\omega _n}}}}$~) into a two-dimensional semitorus;
\item is such that its boundary is the restricted complex 
$\GL_2({F^T_{\o\omega}}\times F^T_{\omega})$-bisemimodule
\[ M^{\rm res}_{F^T_{\o\omega_\oplus}} \otimes M^{\rm res}_{F^T_{\omega _\oplus}} 
= \txt\bigoplus\limits^s_{n=1}\L( M^{\rm res}_{F^T_{\o\omega_n }} \otimes M^{\rm res}_{F^T_{\omega _n}} \R)\]
covered by the real 
$\GL_2({F^{+,T}_{\o v}}\times F^{+,T}_{v})$-bisemimodule
\[ M_{F^{+,T}_{\o v_\oplus}} \otimes M_{F^{+,T}_{v_\oplus}} 
= \txt\bigoplus\limits_{n,m_n}\L( M_{F^{+,T}_{\o v_{n,m_n} }} \otimes M_{F^{+,T}_{v_{n,m_n}}} \R)\]
at the conditions of proposition 2.5.
\Ee
\vskip 11pt

The analytic development of the $T_2(F^T_\omega )$-semimodule $M^{\rm res}_{F^T_{\omega _\oplus}} $ (resp. $T^t_2(F^T_{\o\omega} )$-semi\-module $M^{\rm res}_{F^T_{\o\omega_\oplus }} $~) is the equivalent of the Eisenstein series $\EIS_L(2,n)$ (resp. $\EIS_R(2,n)$~) which is a modular form $f_L(z)$ (resp. $f_R(z)$~), restricted to the upper (resp. lower) half plane.
\vskip 11pt

So, a modular form on the representation space of the complex bilinear algebraic semigroup $\GL_2(F^T_{\o\omega }\times F^T_\omega )$ is a modular biform $f_R(z) \otimes f_L(z)\in C^\infty(\GL_2(F^T_{R }\times F^T_L)\big/$ \linebreak $\GL_2(\ZZ\big/N\ \ZZ)^2)$~.
\vskip 11pt

The analytic development of the $T_2(F^{+,T}_v)$-semimodule  $M_{F^{+,T}_{v_\oplus}}$ (resp.  
$T^t_2(F^{+,T}_{\o v})$-semi\-module  $M_{F^{+,T}_{\o v_\oplus}}$~) is given by the global elliptic $A_L$-semimodule $\phi _L(s_L)\linebreak =\sum\limits_{n,m_n}\lambda _-(n^2_N,m^2_n)\ e^{2\pi inx}$ (resp.
$A_R$-semimodule $\phi _R(s_R)= \sum\limits_{n,m_n}\lambda _+(n^2_N,m^2_n)\ e^{-2\pi inx}$~), in such a way that the $m_n$ representatives of the $n$-th class of $\phi _L(s_L)$ (resp. $\phi _R(s_R)$~), which are semicircles, cover the $n$-th representative of $\EIS_L(2,n)$
(resp. $\EIS_R(2,n)$~), which is a $T^2_L[n]$ (resp. $T^2_R[n]$~) semitorus.
\vskip 11pt

So, the modular representation of $f_L(z)$ (resp. $f_R(z)$~) can be given by a set of $n$~, $1\le n\le \infty $~, two-dimensional semitori $T^2_L[n]$ (resp. $T^2_R[n]$~), restricted to the upper (resp. lower) half plane and covered each one by $m_n$ semicircles of the $n$-th class of the global elliptic semimodule $\phi _L(s_L)$ (resp. $\phi _R(s_R)$~).
\vskip 11pt

This kind of modular representation will be used in this chapter in order to analyze the three following conjectures.
\vskip 11pt


\subsection{The Shimura-Taniyama-Weil conjecture}

The Shimura-Taniyama-Weil conjecture deals with the modular representation of an elliptic curve.  Let us recall the main tools of this modular representation.
\vskip 11pt

Let $T_H(N)_R\otimes T_H(N)_L$ be the ring of products of Hecke operators $T_{q_R}\otimes T_{q_L}$ and $U_{q_R}\otimes U_{q_L}$~.  Assume that a global elliptic $A_{R-L}$-bisemimodule $\phi _R(s_R)\otimes_D\phi _L(s_L)=\linebreak \sum\limits_{n,m}(\lambda _+(n_N^2,m_N^2))q^n_{/\QQ_R}\otimes_D\sum\limits_{n,m}(\lambda _-(n_N^2,m_N^2))q^n_{/\QQ_L}$ corresponds to each normalized eigenform\linebreak $f_R(z)\otimes_Df_L(z)=\sum\limits_n 
a_{n,R}q^n_R \otimes_D \sum\limits_n 
a_{n,L}q^n_L$~, $q_{R,L}=e^{\mp 2\pi iz}$~, $z\in\CC
$~, of a product of Hecke operators.  $\phi _R(s_R)\otimes _D\phi _L(s_L)$ can then be viewed as an automorphic representation of $f_R(z)\otimes_D f_L(z)$ composed of a double tower of one-dimensional irreducible curves $E_f(n_N,m_N)_{R,L}$~: 1) being each one a semitorus $T^1_{R,L}[n,m_n]$ of class $n$ with respect to a quantum $M^I_{[v_1]}$ of class 1 (see section 2.18);
2) having centers and radii given by the pairs $\{{\rm cent}(n_N^2,m_n^2),(r(n_N^2,m_n^2))\}$~.
\vskip 11pt


\subsubsection{{\bbf Euclidian uniformization of the elliptic curve $E(\CC)$}}

Let $E(\CC)$ be an elliptic curve over $\CC$ given by the equation $Y^2=4X^3+AX+B$ (~$A,B\in\CC$~) arising in connection with the nonlinear differential equation $(\wp ')^2=4\wp^3+A\wp+B$ where the $\wp$-function is the Weierstrass $\wp$-function which is periodic and related to a lattice in $\CC$~.

Assume that this lattice in $\CC$ is precisely the lattice $\Lambda ^2_{\omega }$ (resp. $\Lambda ^2_{\o\omega }$~) in the $B_{F_{\omega }}$-semimodule $M_{F_{\omega }}$ (resp. $B_{F_{\o\omega }}$-semimodule $M_{F_{\o\omega }}$~) as developed in section 2.3.

So, an Euclidian uniformization of our elliptic curve $E(\CC)$ of the type introduced by B. Mazur in \cite{Ma1} would be obtained by considering the surjective mapping:
\[ \Es_{\GL_2\to E(\CC)}: \quad \GL_2(F_{R} \times F_{L})\big/ \GL_2((\ZZ\big/N\cdot \ZZ)^2)\To E(\CC)\]
of the quotient of the complex bilinear algebraic semigroup $\GL_2(F_{R} \times F_{L })$ by the subgroup $\GL_2((\ZZ\big/N\cdot \ZZ)^2)$~, constituting a representation of the bilattice $\Lambda ^2_{\o\omega }  \otimes\Lambda ^2_{\omega }$~.  This mapping $\Es_{\GL_2\to E(\CC)}$ identifies the complex points of the elliptic curve $E(\CC)$ with certain ``bipoints'' of the orbit space  $\GL_2(F_{R} \times F_{L })\big/ \GL_2((\ZZ\big/N\cdot \ZZ)^2)$ with respect to $\Lambda ^2_{\o\omega }  \otimes\Lambda ^2_{\omega }$~.
\vskip 11pt

\subsubsection{{\bbf Hyperbolic uniformization of the elliptic curve $E(\CC)$}}

The quotient bilinear semigroup $\GL_2(F_{R} \times F_{L })\big/ \GL_2((\ZZ\big/N\cdot \ZZ)^2)$ is a bisemispace whose representation is the $\GL_2(F_{\o\omega} \times F_{\omega })$-bisemimodule  $M_{F_{\o\omega }}\otimes M_{F_{\omega }}$~. 

As it was developed in chapter 2, a modular representation of 
$M_{F_{\o\omega_\oplus}}\otimes M_{F_{\omega_\oplus }}$ can be only obtained if we consider the restricted
 $\GL_2(F^T_{\o\omega} \times F^T_{\omega })$-bisemimodule 
$M^{\rm res}_{F^T_{\o\omega_\oplus }}\otimes M^{\rm res}_{F^T_{\omega_\oplus}}
=\bigoplus\limits^s_{n=1}\L( M^{\rm res}_{F^T_{\o\omega_n }}
\otimes M^{\rm res}_{F^T_{\omega_n }}\R)$~, 
$1\le n\le s\le \infty $~, having a cuspidal representation given by the product, right by left, $\EIS_R(2,n)\otimes \EIS_L(2,n)$ of the equivalent of Eisenstein series (see section 2.8).

Let $f_R(z)\otimes_D f_L(z)\simeq \EIS_R(2,n)\otimes_D \EIS_L(2,n)$  be the normalized eigenform of the above mentioned products of Hecke operators: it is a product, right by left, of cuspidal forms expanded in formal power series.

Thus, $f_R(z)\otimes_D f_L(z)$ constitutes a cuspidal representation of the restricted
$\GL_2(F^T_{\o\omega} \times F^T_{\omega }))$-bisemimodule 
$M^{\rm res}_{F^T_{\o\omega }}\otimes_D M^{\rm res}_{F^T_{\omega}}$~.

A hyperbolic uniformization of our elliptic curve $E(\CC)$ will be reached by considering the surjective mapping:
\[ \Hs_{\GL^{\rm res}_{2_\CC}\to E(\CC)}: \quad 
M^{\rm res}_{F^T_{\o\omega }}\otimes_D M^{\rm res}_{F^T_{\omega}}\To E(\CC)\]
identifying the complex points of $E(\CC)$ with certain complex bipoints of
$M^{\rm res}_{F^T_{\o\omega }}\otimes_D M^{\rm res}_{F^T_{\omega}}$~.

So, the hyperbolic uniformization of $E(\CC)$ is especially obtained by getting rid of the
 multiples of the subbisemimodules $M_{F^T_{\o\omega_n }}\otimes M_{F^T_{\omega_n}}$~, $1\le n\le s\le \infty $~, i.e. by considering:
\begin{align*}
M^{\rm res}_{F^T_{\o\omega _\oplus}}\otimes_D M^{\rm res}_{F^T_{\omega_\oplus}}
&= \sum^s_{n=1}\sum_{m_{\omega _n}} \L(M_{F^T_{\o\omega_{n,m_{\omega _n}} }}\otimes_D M_{F^T_{\omega_{n,m_{\omega _n}}}}\R)
-\sum^s_{n=1}\sum_{m_{\omega _n}>1} \L(M_{F^T_{\o\omega_{n,m_{\omega _n}} }}\otimes_D M_{F^T_{\omega_{n,m_{\omega _n}}}}\R)\\[11pt]
&= \sum^s_{n=1} \L( M^{\rm res}_{F^T_{\o\omega_{n} }}\otimes_D 
M^{\rm res}_{F^T_{\omega_{n}}} \R)\;.\end{align*}

As $f_R(z)\otimes_D f_L(z)$ is a cuspidal representation of 
$M^{\rm res}_{F^T_{\o\omega}}\otimes_D M^{\rm res}_{F^T_{\omega}} $~, it constitutes a modular representation of the elliptic curve by the surjective mapping:
\[ \Hs^{\rm cusp}_{\GL^{\rm res}_{2_\CC}\to E(\CC)}: \quad 
f_R(z)\otimes_D f_L(z)\To E(\CC)\]
if we take into account the mapping $\Hs_{\GL^{\rm res}_{2_\CC}\to E(\CC)}$~.

This hyperbolic uniformization \cite{Rib2} is realized in the upper and in the lower half planes which are periodic respectively with the congruence subgroups $\Gamma _L(N)$ and $\Gamma _R(N)$ (representations respectively of the Hecke operators $U_{q_L}$ and $U_{q_R}$~) being restricted subgroups of $\Gamma (N)$ (or $\Gamma _1(N)$ or $\Gamma _0(N)$~) and
$\Gamma ^t(N)$ (or $\Gamma ^t_1(N)$ or $\Gamma ^t_0(N)$~) respectively (see section 2.12).
\vskip 11pt

\subsubsection{{\bbf Hyperbolic uniformization of arithmetic type of the elliptic curve $E(\QQ)$}}

The next step consists in envisaging the modular representation of the elliptic curve 
$E(\QQ)$~, covering its complex equivalent $E(\CC)$~, by means of the global elliptic $A_{R-L}$-bisemimodule $\phi _R(s_R)\otimes_D \phi _L(s_L)$~. So, we have to consider the diagram:
\[\begin{CD}
f_R(z)_{\rm res} \otimes_D f_L(z)_{\rm res} 
@>{\Hs^{\rm cusp\ (res)}_{\GL_{2_\CC}^{\rm res}\to\CC}}>> E(\CC)\\
@A{\Ms_{\phi \to f}}AA @AA{\Ms_{E(\QQ)\to E(\CC)}}A\\
\phi _R(s_R)_{\rm res} \otimes_D \phi _R(s_L)_{\rm res} 
@>{\Hs^{\rm cusp\ (res)}_{\GL_{2_\CC}^{\rm res}\to\QQ}}>> E(\QQ)\end{CD}\]
in such a way that the covering $\Ms_{E(\QQ)\to E(\CC)}$ of the elliptic curve $E(\CC)$ over $\CC$  by the elliptic curve $E(\QQ)$ over $\QQ$ can be reached throughout the covering $\Ms_{\phi \to f}$ of the restricted cuspidal biform $f_R(z)_{\rm res}\otimes_D f_L(z)_{\rm res}$ by the restricted global elliptic bisemimodule $\phi _R(s_R)_{\rm res} \otimes_D \phi _R(s_L)_{\rm res} $~, as it will be developed in the following.

As the mapping $\Ms_{\phi \to f}$ was justified in chapter 2, the covering mapping
$\Ms_{E(\QQ)\to E(\CC)}$ of the elliptic curve $E(\CC)$ by the elliptic curve 
$E(\QQ)$ will be justified if the surjective mapping 
$\Hs^{\rm cusp(res)}_{\GL_{2_\RR}^{\rm res}\to E(\QQ)}$ of the elliptic curve $E(\QQ)$ 
by the global elliptic bisemimodule 
$\phi _R(s_R)_{\rm res} \otimes_D \phi _R(s_L)_{\rm res} $ can be justified.  
This surjective mapping $\Hs^{\rm cusp(res)}_{\GL_{2_\RR}^{\rm res}\to E(\QQ)}$ corresponds precisely to the hyperbolic uniformization of arithmetic type of the elliptic curve $E(\QQ)$  because it corresponds to a modular representation of the elliptic curve $E(\QQ)$ by 
$(\phi _R(s_R)_{\rm res} \otimes_D \phi _R(s_L)_{\rm res} )$~.

Now, describing the elliptic curve  $E(\QQ)$ globally over $\QQ$ is equivalent to consider 
the set $\{E(\ZZ\big/n\ \ZZ)\}_{n\in\NN}$ of the elliptic curves, and, especially, the set 
$\{E(\FF_p)\}_p$ of elliptic curves over $\FF_p$~, $p$ being a prime number (because our
 equation of $E(\QQ)$~, having integer coefficients, can be reduced modulo $p$~) 
\cite{Tat}, \cite{Ma2}.

As mentioned above, $\phi _R(s_R) \otimes_D \phi _R(s_L)$ is the direct sum of the $n$ sets of products, right by left, 
\[ \{E_f(n_N,m_n)_R) \otimes  E_f(n_N,m_n)_L\}_{m_n}\equiv
\{T^1_R(n,m_n]\otimes T^1_M[n,m_n]\}_{m_n}\]
of one-dimensional irreducible curves $E_f(n_N,m_N)$ which are semitori $T^1[n,m_n]$~.

So, at each integer $n$~, we have a subset of $m_n$~, $m_n\in\NN$~, products of semitori $T^1_L[n,m_n]$~, localized in the upper half plane, by their equivalents 
$T^1_R[n,m_n]$ in the lower half plane in such a way that
$T^1_R[n,m_n]$ be, for instance, projected on $T^1_L[n,m_n]$~.

These $m_n$-semitori $T^1_L[n,m_n]$ and $T^1_R[n,m_n]$~,  generated under the action of the decomposition group $D_{n^2_N}$ according to section 2.13, are isomorphic.

Taking into account that:
\Bean
\item $E(\FF_p)\simeq (\ZZ\big/ p\ \ZZ)^2$~,
\item $T^1_L(p,m_p]$ (resp. $T^1_R[p,m_p]$~) is a cyclic semigroup whose order is a multiple of $p$~, 
\Ee
it is reasonable to associate with an elliptic curve $E(\FF_p)$ (having good reduction at $p$~), the subset $\{E_f(p_N,m_p)_R\otimes E_f(p_N,m_p)_L\}_{m_p}$ of $m_p$ products of semitori,  $p$ being a prime  integer, in such a way that there exists a natural surjective mapping:
\[ \{E_f(p_N,m_p)_R\otimes E_f(p_N,m_p)_L\}_{m_p}\To E(\FF_p)\]
which identifies the elliptic curve $E(\FF_p)$ with the orbit space of the $p$-th class of products, right by left, of semitori of the global elliptic bisemimodule $\phi _R(s_R)\otimes_D \phi _L(s_L)$ under the action of the decomposition group $D_{p^2_N}$~.
\vskip 11pt

\subsubsection{Proposition}

{\em If we have the natural surjective mapping:
\[\{E_f(p_N,m_p)_R) \otimes  E_f(p_N,m_p)_L\}_{m_p}\To E(\FF_p)\;, \]
then the elliptic curve $E(\FF_p)$ will have $m_p$ generating points given by:
\begin{align*}
\#E(\FF_p)
&= (|\det \rho _{\lambda (p^2_N,m^2_p)}- \tr \rho _{\lambda (p^2_N,m^2_p)}+1|)^\half\\
&= m_p\end{align*}
where $\det \rho _{\lambda (p^2_N,m^2_p)}$ and $\tr\rho _{\lambda (p^2_N,m^2_p)}$ are given in proposition 2.15.}
\vskip 11pt

\bpr
\Be
\item The sum of the $m_p$ products, right by left, of semitori is given, according to section 2.20, by:
\[ \txt \bigoplus\limits_{m_p} 
\L(E_f(p_N,m_p)_R) \otimes  E_f(p_N,m_p)_L\R)
= \bigoplus\limits_{m_p} \L( \lambda _+(p^2_N,m^2_p)\ e^{-2\pi ipx}\otimes_D
\lambda _+(p^2_N,m^2_p)\ e^{2\pi ipx}\R)\; ,\]
 where the coefficients $\lambda _\pm(p^2_N,m^2_p)$ are eigenvalues of products of Hecke operators given by the coset representatives $g_2(p^2_N,b_N)$ in sections 2.12, 2.13, and in proposition 2.15.

\item Recall that the number of points on an elliptic curve $E(\FF_p)$ is traditionally given by:
\[ \#E(\FF_p)=p+1-a_p\]
where $a_p$ results from $T_p(f)=a_p\ f$ if
\Bi
\item $f$ is an eigenfunction of the coset representative $T_p$ of the Hecke operator;
\item $a_p$ is the sum of the eigenvalues of $T_p$.
\Ei
A similar formula can be introduced by taking into account the present context.

\item As the number of Galois automorphisms on a semitoric curve $E_f(p_N,m_p)$ is equal to its rank $r_{E_f(p_N,m_p)}=p_N=p\cdot N$~, it is reasonable to associate with the orbit space of the subset $\{E_f(p_N,m_p)_R\otimes E_f(p_N,m_p)_L\}_{m_p}$ the elliptic curve $E(\FF_p)$ over $\FF_p$ and to think that the number of generating points on $E(\FF_p)$ could be given by the rank
$r^0_{E^2_f(p_N,m_p)}=p^2_N$ of this orbit space \cite{Ma1} viewed as a product, right by left, of semicircles.  Now, this rank $r^0_{E^2_f(p_N,m_p)}$ is precisely given by
\begin{align*}
\det \rho _{\lambda (p^2_N,m^2_p)} &= \lambda _+(p^2_N,m^2_p)\centerdot \lambda _-(p^2_N,m^2_p)\\
&= p^2_N\end{align*}
according to proposition 2.15.

So, $\#E(\FF_p)$ could be given by:
\[\#E(\FF_p)= (|\det \rho _{\lambda (p^2_N,nm^2_p)} |)^\half\;.\]

\item But, the orbit space of these toric curves $E_f(p_N,m_p)$ must be centered at the origin in order to have an isogeny \cite{Kna}.

Then, we have that:
\begin{align*}
\#E(\FF_p)
&= \L( |\det \rho _{\lambda (p^2_N,m^2_p)} -\tr \rho _{\lambda (p^2_N,m^2_p)} +1|\R)^\half\\
&= \L(|p^2_N-(1+m^2_p+p^2_N)+1|\R)^\half = m_p\end{align*}
where:
\Bi
\item $\tr \rho _{\lambda (p^2_N,m^2_p)} =1+m^2_p+p^2_N$~, according to proposition 2.15, contributes to the above mentioned isogeny;
\item the term ``~$+1$~'' comes from the point at infinity.\epr
\Ei
\Ee

\subsubsection{{\bbf Modular representation of the elliptic curve $E(\QQ)$}}

The problem which arises now is to precise the set $\{E(\FF_p)\}_p$ of elliptic curves $E(\FF_p)$ over $\FF_p$ (~$p$ varying) which are locally equivalent to our elliptic curve $E(\QQ)$ to which a modular representation must correspond.  So, we have to introduce the Hecke $L$-series of the elliptic curve $E(\QQ)$ and, more generally, the Hecke $L$-series associated to the cusp forms $f_L(z)$ and $f_R(z)$~.
\vskip 11pt


\subsubsection{\boldmath Definition: Hecke $L$-series}

Let $\lambda _+(q_N^2,m_N^2)$ and $\lambda _-(q_N^2,m_N^2)$ be the eigenvalues of the Hecke operators $U_{q_R}\otimes U_{q_L}(q\mid N)$ and $T_{q_R}\otimes T_{q_L}(q\nmid N)$ (see section 2.3).

We can then define the Hecke $L$-series:\\[6pt]
\indent $ L_R(s_-)=\sum_{n}\lambda _\mp(n_N^2,m_N^2)n^{-s_-}\;, \qquad s_-\in\CC \quad \text{with\ } \{\Im s_-<0\}\;,  $\\[6pt]
in the lower half plane, and\\[6pt]
\indent $ L_L(s_+)=\sum_{n}\lambda _\mp(n_N^2,m_N^2)n^{-s_+}\;, \qquad s_+\in\CC \quad \text{with\ } \{\Im s_+>0\}$\\[6pt]
in the upper half plane, with respect to
\[\lambda_\mp(n^2_N,m^2_N)=\F{(1+m^2_N+n^2_N)\mp[(1+m^2_N+n^2_N)^2-4n^2_N]^{\half}}2\]
(see section 2.13 and proposition 2.15),
and associated respectively with a pseudo-ramified \rl cusp form
\[ f_R(z)=\sum_{n}\lambda_\mp(n^2_N,m^2_N)\ e^{2\pi inz}\quad
\text{(resp.} \quad f_L(z)=\sum_{n}\lambda_\mp(n^2_N,m^2_N)\ e^{-2\pi inz}\;);\quad z\in\CC\;.\]
\vskip 11pt

\subsubsection{Proposition}
{\em 
Let $L_R(s_-)$ and $L_L(s_+)$ be the Hecke $L$-series defined respectively in the lower and in the upper half plane.  The $L^{\rm deg}_{R-L}(\Re(s))$ Hecke $L$-series can be defined from their product (degenerate case):
\[ L^{\rm deg}_{R-L}(\Re(s))=L_R(s_-)\cdot L_L(s_+)=\sum_{n}\lambda^2_{\pm}(n_N^2,m_N^2)n^{-2x}\;, \qquad x\in \Re(s_{\pm})\]
which corresponds to the classical Eulerian development {\em \cite{La}}:
\[L^{\rm deg}_{R-L}(\Re(s))=\prod_{q\mid N} (1-\lambda^2 _{\mp}(q_N^2)q^{-2x})^{-1}\prod_{q\nmid N}
(1-\lambda^2 _{\mp}(q_N^2)\ \varepsilon (q)^2q^{2k-2x-2})^{-1}\]
where \Bi
\item $\varepsilon (q):(\ZZ\ / \ N\ZZ)^*\to \CC^*$ is the Dirichlet character
\item $\lambda ^2_\mp(q^2_N)=\sum\limits_m\lambda ^2_\mp(q^2_N,m^2_N)$~.
\Ei}
\vskip 11pt

{\parindent=0pt 
\bpr  \Bean
\item As the Hecke $L$-series $L_{R-L}(\Re(s))$ is defined on a space of semisimple type, the off-diagonal terms of the product of $L_R(s_-)$ by $L_L(s_+)$ must be zero.

\item Taking into account that the $\lambda_{-}(q^2_N,m^2_N)$ are eigenvalues of $(D_{q^2_N;m^2_N})\cdot\alpha _{q^2_N}$ (see section 2.13) and considering the statement of proposition 2.15, it is clear that $L^{\rm deg}_{R-L}(\Re(s))$ is also an Artin $L$-function.
\epr\Ee}
 \vskip 11pt

\subsubsection{Corollary}

{\em The Hecke $L$-series $L_{R-L}(s)=L_R(s_-)\cdot L_L(s_+)$ have also the following Eulerian development (nondegenerate case):
\begin{align*}
L_{R-L}(s)&=\prod_{q\mid N} (1-\lambda_\mp(q^2)q^{-s_-})^{-1}
\prod_{q\nmid N} (1-\lambda _\mp(q_N^2)\ \varepsilon (q) q^{k-s_--1})^{-1}\\
&\quad \cdot\prod_{q\mid N} (1-\lambda_\mp(q_N^2)q^{-s_+})^{-1}
\prod_{q\nmid N } (1-\lambda _\mp(q_N^2)\ \varepsilon (q) q^{k-s_+-1})^{-1}\;.\end{align*}
}
\vskip 11pt


 \subsubsection{\boldmath Hecke $L$-series of the elliptic curve $E(\QQ)$}  

Considering the canonical diagonal injective mapping of $\QQ$ into the adele ring $\Aa$~, which can be seen from the discrete topology from $\Aa$ into $\QQ$~, \cite{Rob}, we can admit that the injection
\[ E_{\QQ\to\prod\limits_q} : \quad E(\QQ) \To E(\txt\prod\limits_q\FF_q)\;, \]
where $\prod\limits_q$ is taken over all primes $q$~, allows to transfer the study of $E(\QQ)$ onto a subset of $E(\prod_q\FF_q)$~.\vskip 11pt

But, if we realize that the kernel $\Ker(E^{-1}_{\QQ\to\prod\limits_q} )$ of the map  $E^{-1}_{\QQ\to\prod\limits_q} $ refers precisely to these primes ``~$p_\perp$~'' which do not  enter in the generation of the elliptic curve $E(\QQ)$~, then the Hecke $L$-series  $L\RL(s_{\mp})$  can be partitioned into two complementary subseries following: \[L\RL(s_{\mp})=L\RL(s_{\mp},E(\QQ))+L\RL(s_{\mp},E(\txt\prod\limits_{p_\perp}\FF_{p_\perp}))\] where
\Bi
\item $L\RL(s_{\mp},E(\QQ))$ is the $L$-subseries of the elliptic curve $E(\QQ)$ such that their Euler factors refer to the set of primes $\{p\}$ complementary to the set of primes $p_\perp$~:
\[ L\RL(s_{\mp},E(\QQ))=\sum_{n_{g}}\lambda _\mp (n^2_{g_N},m^2_N)n_{g}^{-s_\mp}\;, \quad
n_{g} \le n\;.\]

\item $L\RL(s_{\mp},E(\prod\limits_{p_\perp}\FF_{p_\perp}))$ is the ``virtual" subseries of $L\RL(s_{\mp})$ with regard to\linebreak $L\RL(s_{\mp},E(\QQ))$~;
\item $L\RL(s_{\mp},E(\prod\limits_{p_\perp}(\FF_{p_\perp})
))=\sum\limits_{n-n_{g}}\lambda_{\mp}((n_N-n_{g_N})^2,m_N^2)\ (n-n_{g})^{-s_{\mp}}$~, $(n_N-n_{g_N})\in\nit$~, $n_N\equiv n\cdot N$~, $n\ge n_g$~;
\item $\{q\}=\{p^\perp\}\cup\{p\}$~.
\Ei  \vskip 11pt  

\subsubsection{Proposition (Hyperbolic uniformization of arithmetic type)}

{\em Let
\begin{align*}
&L_R(s_-,E(\QQ))\otimes_D L_L(s_+,E(\QQ))\\
&\quad = \sum_{n_g}( \lambda _-(n^2_{g_N},m^2_{n_g})\ n_g^{-s_-} \otimes_D
\lambda _+(n^2_{g_N},m^2_{n_g})\ n_g^{-s_+} )\\
& \quad = \prod_{p\mid N} (1- \lambda _\mp(p^2_{N},m^2_{p})\ p^{-s_-} )^{-1} \ 
\prod_{p\nmid N} (1- \lambda _\mp(p^2_{N},m^2_{p})\ \varepsilon (p)\ p^{k-s_--1} )^{-1}\\
& \qquad  \prod_{p\mid N} (1- \lambda _\mp(p^2_{N},m^2_{p})\ p^{-s_+} )^{-1} \ 
\prod_{p\nmid N} (1- \lambda _\mp(p^2_{N},m^2_{p})\ \varepsilon (p)\ p^{k-s_+-1} )^{-1}
\end{align*}
be the restricted product of $L$-subseries attached to the elliptic curve $E(\QQ)$ and corresponding to the restricted cuspidal biform of weight two and level $N$~:
\begin{align*}
& f_R(z)_{\rm res}\otimes_D f_L(z)_{\rm res}\\
&\qquad = \sum_{n_g} \L( \lambda _-(n^2_{g_N},m^2_{n_g})\ e^{-2\pi in_gz}\otimes_D
\lambda _+(n^2_{g_N},m^2_{n_g})\ e^{2\pi in_gz} \R)\;, \quad z\in\CC\;,\end{align*}
covered by the restricted global elliptic bisemimodule:
\begin{align*}
& \phi _R(s_R)_{\rm res}\otimes_D \phi _L(s_L)_{\rm res}\\
& \qquad = \sum_{n_g} \sum_{m_{n_g}} \L( \lambda _-(n^2_{g_N},m^2_{n_g})\ e^{-2\pi in_gx}\otimes_D
\lambda _+(n^2_{g_N},m^2_{n_g})\ e^{2\pi in_gx} \R)\;, \quad x\in\RR\;.\end{align*}
Then, the modular representation of the elliptic curve $E(\QQ)$ over $\QQ$ will be worked out from the ``~$p$~'' sets of surjective mappings:
\[ \{ E_f(p_N,m_p)_R \otimes E_f(p_N,m_p)_L \}_{m_p}\To E(\FF_p)\;, \]
for all prime $p$ entering into the restricted eulerian product of $L_R(s_-,E(\QQ)) \otimes_D\linebreak
L_L(s_+,E(\QQ))$~, in such a way that the orbit spaces of  
$\{ E_f(p_N,m_p)_R \otimes E_f(p_N,m_p)_L \}_{m_p}$ are associated in the sense of proposition 3.1.4, with the elliptic curves $E(\FF_p)$~, for the above mentioned $p$ primes, having good reduction and conductor $N$~.  This modular representation of the elliptic curve $E(\QQ)$ corresponds to its hyperbolic uniformization of arithmetic type {\em \cite{Ma1}}, and, thus to the Shimura-Taniyama-Weil conjecture.}
\vskip 11pt

\subsubsection{Connection with Diophantine equations}

As it is well known, the Shimura-Taniyama-Weil conjecture is related to 
the problem of Diophantine equations especially by means of the Mordell-Weil group 
of the elliptic curve $E(\QQ)$~.  Indeed, Mordell proved that $E(\QQ)$ is a 
finitely generated abelian group, the generators of $E(\QQ)$ being the rational points on the 
curve or the rational solutions of the equation \cite{Gou}.  
So, in the context developed here, the generators of $E(\QQ)$ are precisely the 
$p$ sets of $m_p$ points associated with the $p$ surjective mappings 
$\{E_f(p_N,m_p)_R\otimes E_f(p_N,m_p)_L\}\To E(\FF_p)$ introduced in proposition 3.1.4.\vskip 11pt


\subsection{The Riemann conjecture}

\subsubsection{The trivial zeros and the zeta functions}  

The ``(pseudo-)unramified'' Hecke $L$-series (with $N=1$~), correspond to the classical $L$-series: they are in one-to-one correspondence with the classical zeta function $\zeta (s)$~, are denoted $L^{nr}\RL(s_{\mp})$ and given by:
\begin{align*}
L^{nr}_R(s_-) &= \sum_{n}\lambda^{nr}_\mp(n^2,m^2)\ n^{-s_-}\;, \\
L^{nr}_L(s_+) &= \sum_{n}\lambda^{nr}_\mp(n^2,m^2)\ n^{-s_+}\;,\end{align*}
where $\lambda ^{nr}_\mp(n^2,m^2)$ is introduced in corollary 2.1.6.\vskip 11pt

The Hamburger theorem \cite{Chan} tells us that $L_{R,L}^{nr}(s_{\mp})=a_1 \zeta_{R,L}(s_{\mp})$ where $a_1\in \CC$ and where $\zeta_{R,L}(s_{\mp})=\sum\limits_n n^{-s_{\mp}}$ is the classical zeta function associated with the right or the left case.

As it is done classically \cite{Ing}, let us introduce the right (resp. left) functions:
\begin{align*}
\xi_{R,L}(s_{\mp})&= (s_{\mp}-1)\Pi ^{-\half s_{\mp}}\Gamma ({\txt \half} s_{\mp}+1)\zeta_{R,L}(s_{\mp})\\
&= h_{R,L}(s_{\mp})\zeta_{R,L}(s_{\mp})\end{align*}
which satisfies the functional equation $\xi_{R,L}(1-s_{\mp})=\xi_{R,L}(s_{\mp})$~.

The poles of $h_{R,L}(s_{\mp})$~, i.e. of the gamma function $\Gamma ({\txt \half} s_{\mp}+1)$~, are simple ones at $s_{\mp}=-2,-4,-6,\cdots,-2n$~.  Since these are points at which $\xi_{R,L}(s_{\mp})$ is regular and not zero, they must be simple zeros of $\zeta_{R,L}(s_{\mp})$ \cite{Tit}.
\vskip 11pt

The trivial zeros $s_{\mp}=-2,-4,\cdots,-2n$ of $\zeta_{R,L}(s_{\mp})$ will be interpreted in the present context as being 
equal at a sign to $2f_{\o v_n}$ (resp. $2f_{v_n}$~) where $f_{v_n}$ is the global class residue degree of the $v_n$-th real place.  The factor 2 proceeds from the fact that the left (resp. right) places are defined in the upper (resp. lower) half space on which are defined curves isomorphic to $1D$-semitori: thus, these $1D$-semitori must be doubled leading to analytic continuation to the whole $s$-plane and their global class residue degrees must also be multiplied by two.
\vskip 11pt

More concretely, the trivial zeros of $\zeta \RL(s_{\mp})$ can be interpreted as follows.

As the (pseudo-)unramified cuspidal biform $f^{nr}_R(z)\otimes_D f^{nr}_L(z)$ is in one-to-one correspondence with the product, 
right by left, $\zeta _R(s_-)\otimes_D\zeta _L(s_+)=\sum\limits_n n^{-2x}$~, $x\in\Re(s_\mp)$~, of classical zeta functions, 
generating a zeta function $\sum\limits_nn^{-2x}$ over $\RR$~, it appears that the supercuspidal representation of the cuspidal 
biform $f^{nr}_R(z)\otimes_Df^{nr}_L(z)$ may degenerate into the pseudo-unramified simple real global elliptic bisemimodule:
\[ \phi ^{nr,0}_R(s_R)\otimes_D\phi ^{nr,0}_L(s_L)
= \sum_n(\lambda^{nr} (n^2)\ e^{-2\pi inx}\otimes _D \lambda^{nr} (n^2)\ e^{2\pi in x})\;,\]
the (pseudo-)unramification leading to consider the value $N=1$ for the conductor and the ``simplicity'' consisting in considering the irreducible curves $\pi _L(n)$ and $\pi _R(n)$ of global class $n$ having multiplicity one (see sections 2.20 and 2.22): so $m_R(n)=m_L(n)=1$~.

This global elliptic bisemimodule $\phi ^{nr,0}_R(s_R)\otimes _D\phi ^{nr,0}_L(s_L)$ can be interpreted geometrically as the direct sum of products of irreducible (pseudo-)unramified semitoric curves corresponding to each other by pairs of same level $n$ and constituting the automorphic representation space of the pseudo-unramified algebraic bilinear semigroup $\GL_2(F^{+,T,nr}_{\o v}\times F^{+,T,nr}_v)$ and the orbit space of $f^{nr}_R(z)\otimes_D f^{nr}_L(z)$~.

The interpretation of the trivial zeros of $\zeta \RL(s_\mp)$ can be given by means of the homomorphism $H_{\phi _{R\times L}
\to\zeta _{R\times L}}$ introduced in the following proposition.
\vskip 11pt 

\subsubsection{Proposition}

{\em The one-to-one correspondence between the double pseudo-unramified simple global elliptic bisemimodule $2\phi _R^{nr,0}(s_R)\otimes_D 2\phi ^{nr,0}_L(s_L)$ and the product, right by left, of zeta functions $\zeta _R(s_-)\otimes_D\zeta _L(s_+)$ is given by the homomophism:
\begin{align*}
H_{\phi _{R\times L}
\to\zeta _{R\times L}} : \quad
2\phi ^{nr,0}_R(s_R)\otimes_D2\phi _L^{nr,0}(s_L)
&= \sum_n(2\lambda^{nr}  (n^2)e^{-2\pi inx}\otimes_D2\lambda ^{nr} (n^2)e^{2\pi inx}\\
&\quad \To \zeta _R(s_-)\otimes_D \zeta _L(s_+)=\sum_n(n^{-s_-}\otimes_Dn^{-s_+})\end{align*}
such that the kernel $\Ker(H_{\phi _{R\times L}\to\zeta _{R\times L}})$ of  $H_{\phi _{R\times L}\to\zeta _{R\times L}}$ maps into the set of products, right by left, of the trivial zeros of $\zeta _R(s_-)$ and of $\zeta _L(s_+)$~.
}
\vskip 11pt 

{\parindent=0pt
\bpr  Indeed, the kernel $\Ker (H_{\phi _{R\times L}\to\zeta _{R\times L}})$ is given by the set of bipoints:
\[\{\sigma _{n_R}\times \sigma _{n_L}=4(\lambda^{nr} ) ^2(n^2)=2\lambda^{nr}  (n^2)(e^{-2\pi inx}\mid x=0)\times 2\lambda^{nr}  (n^2)
(e^{2\pi inx}\mid x=0))\}^\infty _{n=1}\]
localized on the real axis $\sigma $ in such a way that:
\Bi
\item the \lr point $\sigma _{n_L}=2\lambda^{nr}  (n^2)(e^{2\pi inx}\mid x=0)$
(resp. $\sigma _{n_R}=\linebreak 2\lambda^{nr}  (n^2)(e^{-2\pi inx}\mid x=0)$~) corresponds to the degeneracy of the irreducible (toric) curve
$ 2\lambda ^{nr} (n^2)e^{2\pi inx}$ (resp. $ 2\lambda^{nr} (n^2)e^{-2\pi inx}$~);
\item the bipoint $\sigma _{n_R}\times \sigma _{n_L}$ is defined as the product of the right point $\sigma _{n_R}$ by its left
equivalent $\sigma _{n_L}$~, and is equal to $\sigma _{n_R}\times \sigma _{n_L}=4(\lambda^{nr} )^2 (n^2)=4f^2_{v_n}=4\ n^2$~.
\Ei
So, as the kernel $\Ker (H_{\phi _{R\times L}\to\zeta _{R\times L}})$ of $ H_{\phi _{R\times L}\to\zeta _{R\times L}}$ is the set of bipoints $\{\sigma _{n_R}\times \sigma _{n_L}=4\ n^2\}^\infty _{n=1}$ corresponding to the trivial zeros of
$\zeta _R(s_-)\otimes \zeta _L(s_+)$~, we can see that the trivial zeros of $\zeta \RL(s_\mp)$ are the integers $s_\mp=-2,-4,\cdots,-2n$~.\epr}
\vskip 11pt 

The non-trivial zeros of $\zeta\RL(s_\mp)$ can then be generated from the corresponding trivial zeros by taking into account the following considerations.
\vskip 11pt 

\subsubsection{The non-trivial zeros of the zeta functions}

Let 
\[\alpha_{4n^2}=\BM 4n^2 & 0 \\ 0&1\EM\]
be the split Cartan subgroup element associated with the quadratic place $v^2_{2n}=\o v_{2n}\times v_{2n}$~. Let
\[ D_{4n^2;i^2}=\BM 1&i\\ 0&1\EM \BM 1&0\\ i&0\EM\]
be the coset representative of the Lie algebra of the decomposition group acting on $\alpha_{4n^2}$~: it corresponds to the  coset representative of an unipotent Lie algebra mapping in the topological Lie algebra $\gfrak\ell_2(F^{+,T,nr}_{\o v}\times F^{+,T,nr}_v)$~.  Let $\lambda^{nr} _+(4n^2,i^2)$ and $\lambda^{nr} _-(4n^2,i^2)$ be the eigenvalues of $(D_{4n^2,i^2}\cdot \alpha_{4n^2})$ interpreted as weights in section  2.16.
\vskip 11pt

The Lie algebra $(\gfrak\ell_2(F^{+,T,nr}_{\o v}\times F^{+,T,nr}_v)$ 
consists in vector fields on the Lie group $GL_2( F^{+,T,nr}_{\o v}\times F^{+,T,nr}_v )$~. Let
\[\Es_{4n^2}=\BM E_{4n^2} &0\\ 0&1\EM\]
be the infinitesimal generator of this Lie algebra $\gfrak\ell_2( F^{+,T,nr}_{\o v}\times F^{+,T,nr}_v )$ corresponding to the quadratic place $v^2_{2n}$~.
\vskip 11pt

Every root of this Lie algebra $ \gfrak\ell_2( F^{+,T,nr}_{\o v}\times F^{+,T,nr}_v )$ is determined by the (equivalent) eigenvalues $\lambda^{nr} _{\pm}(4n^2,i^2,E_{4n^2})$ of
\[ D_{4n^2,i^2}\cdot \Es_{4n^2}\cdot \alpha_{4n^2}=\left[ \BM 1&i \\ 0&1\EM \BM 1&0\\ i&1\EM\right] \BM E_{4n^2}& 0 \\ 0&1 \EM\BM 4n^2 & 0 \\ 0&1\EM \;.\]
They are given by:
\[ \lambda^{nr} _{\pm}(4n^2,i^2,E_{4n^2})=\F{1\pm i\sqrt{16n^2\cdot E_{4n^2}-1}}2\; \cdotp\]
The fact that the representation of the Lie algebra $\gfrak \ell_2( F^{+,T,nr}_{\o v}\times F^{+,T,nr}_v )$ is composed of  a tower of sections of a vertical tangent bundle explains why the coset representative $D_{4n^2,i^2}$ has  been chosen with $b=i \equiv \sqrt{-1}$~.
\\
$D_{4n^2,i^2}\cdot \Es_{4n^2}$ is in fact a coset representative of the Lie algebra of the decomposition group $D_{i^2}(\ZZ)_{|4n^2}$ noted $\Lie(D_{i^2}(\ZZ)_{|4n^2})$~.
\vskip 11pt

Then, we have the following proposition.
\vskip 11pt

 \subsubsection{Proposition}
{\em Let $D_{4n^2,i^2}\cdot \Es_{4n^2}$ be a coset representative of the Lie algebra of the decomposition group $D_{i^2}(\ZZ)_{|4n^2}$ and let $\alpha_{4n^2}$ be the corresponding split Cartan subgroup element. \\
Then, the products of the pairs of the trivial zeros of the Riemann zeta functions $\zeta_R(s_-)$ and $\zeta_L(s_+)$ are mapped into the products of the corresponding pairs of the non-trivial zeros as follows:
\begin{align*}
D_{4n^2,i^2} \cdot \Es_{4n^2}: \quad \det (\alpha_{4n^2}) &\To \det (D_{4n^2,i^2}\cdot \Es_{4n^2}\cdot \alpha_{4n^2})_{SS}\\
\{(-2n)\cdot (-2n)\} & \To \{\lambda^{nr} _+(4n^2,i^2,E_{4n^2})\cdot 
\lambda^{nr} _-(4n^2,i^2,E_{4n^2})\}\;, \tag*{$\forall\ n\in\NN\;,$}\end{align*}
where $(\ )_{SS}$ denotes the semisimple form of $D_{4n^2,i^2}\cdot \Es_{4n^2}\cdot \alpha_{4n^2}$~.}
\vskip 11pt

\bpr The squares of the trivial zeros $(-2n)^2$ can be interpreted as the squares of the global class residue degrees.\\
As $D_{4n^2,i^2}\cdot \Es_{4n^2}$ is of Galois type, it maps squares of trivial zeros $(-2n)^2$ into products of corresponding pairs of other zeros $\lambda^{nr} _+(4n^2,i^2,E_{4n^2})\cdot\lambda^{nr} _-(4n^2,i^2,E_{4n^2})$ which are non-trivial zeros since $\lambda^{nr} _{\pm}(4n^2,i^2,E_{4n^2})$ have real parts localized on the line $\sigma=\half$~.
 $E_{4n^2}\simeq \F{\hbar^2}{c^2}\ \omega^2_{2n}$ is the square of the energy of a subtorus of rank one in a maximal pseudo-unramified torus of global level $2n$ where $\omega_{2n}$ is the angular frequency of the subtorus.
\vskip 11pt

Remark that the correspondence between Riemann non-trivial zeros and eigenvalues of Hamiltonians has been postulated for a long time (see for example \cite{B-K}, \cite{K-S}).
\vskip 11pt

\subsubsection{Cuspidal new forms}

The non-trivial zeros $\{\lambda^{nr}  _+(4n^2,i^2,E_{4n^2})\}_{n=1}$ and 
$\{\lambda ^{nr} _-(4n^2,i^2,E_{4n^2})\}_{n=1}$ of $\zeta _R(s_-)$ and $\zeta _R(s_+)$ allow to consider the following double (pseudo-)unramified simple global elliptic bisemimodule
\[ \phi _R^{nr,i}(s_R) \otimes \phi _L^{nr,i}(s_L) 
= \sum_n( \lambda^{nr}  _-(4n^2,i^2,E_{4n^2})\ e^{-2\pi inx}\otimes_D \lambda^{nr}  _+(4n^2,i^2,E_{4n^2})\ e^{2\pi inx})\]
constituting, as orbit space, the (pseudo-)unramified supercuspidal representation of the (pseudo-)unramified cuspidal biform $f_R^{nr,{\rm new}}(z)
\otimes_D f_L^{nr,{\rm new}}(z)$ which is a cuspidal ``new'' biform. Indeed, the new forms
$f_R^{nr,{\rm new}}(z)$ and $f_L^{nr,{\rm new}}(z)$ are defined in the orthogonal complement space with respect to the old forms $f_R^{nr}(z)$ and $f_L^{nr}(z)$ considered above since $f_R^{nr,{\rm new}}(z)\otimes_D
f_L^{nr,{\rm new}}(z)$ constitutes a supercuspidal representation of the Lie algebra $\gfrak\ell_2
( F^{+,T,nr}_{\o v} \times F^{+,T,nr}_{v} )$~.
\vskip 11pt

\subsubsection{Corollary}

{\em The eigenvalues $\lambda^{nr}  _+(4n^2,i^2,E_{4n^2})$ and $\lambda^{nr}  _-(4n^2,i^2,E_{4n^2})$ for all $n\in \NN$ are non-trivial zeros of the Riemann zeta function $\zeta(s)=\sum\limits_n n^{-s}$~.}
\vskip 11pt

{\parindent=0pt
\bpr  The set $\{-2n\}$~, being trivial zeros of the Riemann zeta functions $\zeta_R(s_-)$ and $\zeta_L(s_+)$~, is also a set of trivial zeros of the classical Riemann function $\zeta(s)=\sum\limits_nn^{-s}$~.

So, the eigenvalues $\lambda^{nr}  _+(4n^2,i^2,E_{4n^2})$ and $\lambda^{nr}  _-(4n^2,i^2,E_{4n^2})$~, which are localized on the line $\sigma =\half$ and disposed symmetrically on this line with respect to $\tau =0$ if $s=\sigma +i\tau $~, constitute non-trivial zeros of the function $\zeta(s)$~.\epr}
\vskip 11pt

\subsection{The Birch and Swinnerton-Dyer conjecture}

The conjecture of Birch and Swinnerton-Dyer is closely related and dependent to the conjectures of Shimura-Taniyama-Weil and Riemann.  So, we refer to the two first sections of this chapter for the mathematical tools developed there and needed here.
\vskip 11pt

\stepcounter{subsubsection}

\paragraph{\thesubsubsection.} {\bf The basic Birch and Swinnerton-Dyer conjecture\/} for an elliptic curve over $\QQ$ asserts that the rank of $E(\QQ)$ is the order of vanishing of its $L$-function $L(s,E)$ at the point $s=1$ \cite{B-S-D}, \cite{Hus}, \cite{Rub}.\vskip 11pt

So, we have to study the non-trivial zeros of the $L$-subseries $L^{nr}\RL(s_\mp,E(\QQ))$~, i.e. all the zeros of 
$L^{nr}\RL(s_\mp,E(\QQ))$ lying on $\Re(s_\mp)=1$~.
\vskip 11pt

 \subsubsection{The trivial zeros of the restricted zeta functions}

In contrast with the Riemann conjecture, the pseudo-unramified $L$-subseries $L^{nr}\RL(s_\mp,E(\QQ))$ are assumed to be restricted in the sense that
$L^{nr}_R(s_-,E(\QQ))\otimes _DL^{nr}_L(s_+,E(\QQ))$\linebreak $=\sum_{n_{g}}\lambda ^{nr}_\mp(n^2_{g},[m^2])^2\ n^{-2x}_{g}$~, $n_g\le n$~, $x\in\Re(s_\mp)$~, defined over $\RR$~, is in one-to-one correspondence with the representation of the restricted cuspidal biform $f^{nr}_R(z)_{\rm res}\otimes_D f^{nr}_L(z)_{\rm res}$ degenerating into the restricted simple global elliptic bisemimodule:
\[\phi ^{nr,0}_R(s_R)_{\rm res}\otimes_D\phi ^{nr,0} _L(s_L)_{\rm res}
= \sum_{n_{g}}( \lambda^{nr}  (n^2_{g})\ e^{-2\pi in_{g}x}\otimes_D
\lambda ^{nr} (n^2_{g})\ e^{2\pi in_{g}x})\;.\]

This global elliptic bisemimodule can be interpreted geometrically as the direct sum of products  of irreducible pseudo-unramified semitoric curves corresponding to each other by pairs of same class $n_g$ and constituting the automorphic representation space of the (pseudo-)unramified algebraic bilinear semigroup $\GL_2(F^{+,T,nr}_{\o v_g}\times F^{+,T,nr}_{v_g})$ and the orbit space of $f^{nr}_R(z)_{\rm res}\otimes_D f^{nr}_L(z)_{\rm res}$~.
\vskip 11pt 

\subsubsection{Proposition}

{\em The one-to-one correspondence between the double restricted simple global elliptic bisemimodule
$2\phi ^{nr,0} _R(s_R)_{\rm res}\otimes_D 2\phi ^{nr,0} _L(s_L)_{\rm res}$ and the product, right by left, of restricted zeta functions
$\zeta ^{\rm res}_R(s_-)\otimes_D \zeta ^{\rm res}_L(s_+)$ is given by the homomorphism:
\begin{align*}
H^{\rm res}_{\phi_{R\times L}\to\zeta _{R\times L}}: \quad
& 2\phi ^{nr,0} _R(s_R)_{\rm res}\otimes_D 2\phi ^{nr,0} _L(s_L)_{\rm res}
= \sum_{n_{g}}(2  \lambda^{nr}  (n^2_{g})\ e^{-2\pi in_{g}x}\otimes_D
2\lambda^{nr}  (n^2_{g})\ e^{2\pi in_{g}x})\\
& \qquad \qquad \qquad\To \zeta ^{\rm res}_R(s_-)\otimes_D \zeta ^{\rm res}_L(s_+)
=\sum_{n_g}n_g^{-s_-}\otimes_Dn_g^{-s_+})=\sum_{n_g}n_g^{-2x}
\end{align*}
such that its kernel $\Ker(H^{\rm res}_{\phi_{R\times L}\to\zeta _{R\times L}})$ maps into the set of products, right by left, 
of trivial zeros of $\zeta ^{\rm res}_R(s_-)$ and of $ \zeta ^{\rm res}_L(s_+)$~.
}
\vskip 11pt 

{\parindent=0pt
\bpr This is an adaptation of proposition 3.2.2 to the restricted case.  Thus, the kernel
$\Ker(H^{\rm res}_{\phi_{R\times L}\to\zeta _{R\times L}})$ is the set of bipoints:
\[
\{\sigma _{n_{g_R}}\times \sigma _{n_{g_L}}
= 4(\lambda ^{nr} (n^2_{g}))^2= 2\lambda^{nr}  (n^2_{g})( e^{-2\pi in_{g}x}\mid x=0)\times
2\lambda^{nr}  (n^2_{g} ) (e^{2\pi in_{g}x}\mid x=0)\}\]
localized on the real axis $\sigma $ in such a way that the \lr point $\sigma _{n_{g_L}}$ (resp. $\sigma _{n_{g_R}}$~) corresponds to the degeneracy of the irreducible semitoric curve
$2\lambda ^{nr} (n^2_{g})\ e^{2\pi in_{g}x})$ (resp. $2
\lambda ^{nr} (n^2_{g})\ e^{-2\pi in_{g}x}$~).

Remark that $\{\sigma _{n_{g_R}}\times \sigma _{n_{g_L}}=4(\lambda ^{nr} (n^2_{g}))^2
=4n^2_g\}_{n_g}$ is the set of trivial zeros of $\zeta ^{\rm res}_R(s_-)\otimes_D \zeta ^{\rm res}_L(s_+)$~.
\epr}
\vskip 11pt 

\subsubsection{The non-trivial zeros of the restricted zeta functions}

As in section 3.2.3, let $\alpha _{4n ^2_g}= \L(\begin{smallmatrix} 4n ^2_g & 0\\ 0&1\end{smallmatrix}\R)$
 be the split Cartan subgroup element corresponding to the restricted quadratic place $v^2_{2n_g}$~.

Let $D_{4n^2_g;i^2}=\L(\begin{smallmatrix} 1& i\\ 0&1\end{smallmatrix}\R)\
\L(\begin{smallmatrix} 1& 0\\ i&1\end{smallmatrix}\R)$ be the coset representative of the Lie algebra of the decomposition group acting on
$\alpha _{4n^2_g}$~.

The infinitesimal generator of the Lie algebra $\gfrak \ell_2( F^{+,T,nr}_{\o v_g}\times F^{+,T,nr}_{v_g})$ at the quadratic place
$v^2_{2n_g}$ of the Lie group $\GL_2(F^{+,T,nr}_{\o v_g}\times F^{+,T,nr}_{v_g})$ is given by
\[ \varepsilon _{4n^2_g}=\BM E_{4n^2_g} & 0\\ 0&1\EM\;.\]
The roots of $\gfrak \ell_2( F^{+,T,nr}_{\o v_g}\times F^{+,T,nr}_{v_g})$ are determined by the eigenvalues
$\lambda^{nr}  _\pm(4n ^2_g,i^2,E_{4n^2_g})$ of:
\[D_{4n ^2_g}\cdot
\varepsilon _{4n^2_g}\cdot \alpha _{4n^2_g}\cdot M_{\half\to 1}
= \L[\BM 1&i\\ 0&1\EM\ \BM 1&0\\ i&1\EM\R]\ \BM E_{4n^2_g}& 0\\ 0&1\EM\ \BM 4n^2_g& 0\\ 0&1\EM\ \BM 1&0\\ 0&2\EM\;,\]
where the matrix $M_{\half\to 1}=\L(\begin{smallmatrix} 1&0\\ 0&2\end{smallmatrix}\R)$ maps the non-trivial zeros from $\Re(s^-_+)\equiv\sigma =\half$ to $\sigma =1$~, as it is envisaged in the Birch-Swinnerton-Dyer conjecture.  (This map may be of Galois type.)
They are given by:
\[\lambda^{nr}  _\pm(4n^2_g,i^2,E_{4n^2_g})
= 1\pm i\sqrt{8n^2_gE_{4n^2_g}-1}\;.\]
\vskip 11pt 

\subsubsection{Proposition}

{\em The roots of $\gfrak \ell_2( F^{+,T,nr}_{\o v_g}\times F^{+,T,nr}_{v_g})$ given by the eigenvalues
$\lambda^{nr}  _\pm(4n ^2_g,i^2,E_{4n^2_g})= 1\pm i\ \sqrt{8n^2_gE_{4n^2_g}-1}$ of
$D_{4n ^2_g;i^2}\cdot
\varepsilon _{4n^2_g}\cdot \alpha _{4n^2_g}\cdot M_{\half\to1}$ are the non-trivial zeros of the restricted Riemann zeta functions
$\zeta ^{\rm res}_R(s_-)$ and $\zeta ^{\rm res}_L(s_+)$~.
}
\vskip 11pt 

{\parindent=0pt
\bpr Referring to proposition 3.2.4, we see that the products of the pairs of trivial zeros of the restricted Riemann zeta functions
$\zeta ^{\rm res}_R(s_-)$ and $\zeta ^{\rm res}_L(s_+)$ are mapped into the products of the corresponding pairs of the non-trivial zeros according to:
\begin{align*}
D_{4n ^2_g}\cdot
\varepsilon _{4n^2_g}: \quad
\det ( \alpha _{4n^2_g})
&\To \det (D_{4n ^2_g}\cdot
\varepsilon _{4n^2_g}\cdot \alpha _{4n^2_g}\cdot M_{\half\to1})_{\rm ss}\;, \\
\{(-2n_g)(-2n_g)\}
&\To \{ \lambda ^{nr} _+(4n^2_g,i^2,E_{4n ^2_g})
\cdot \lambda^{nr}  _-(4n^2_g,i^2,E_{4n ^2_g})\}\;.\end{align*}

Since $D_{4n ^2_g}\cdot
\varepsilon _{4n^2_g}$ is of Galois type, it maps squares of trivial zeros $ (-2n_g)^2$ into products of corresponding pairs of other zeros $ \lambda^{nr}  _+(4n^2_g,i^2,E_{4n ^2_g})\cdot \lambda^{nr}  _-(4n^2_g,i^2,E_{4n ^2_g})$ which are non-trivial since their real parts are localized on the real line $\sigma =1$~.

Note that: $E_{4n^2_g N^2}=E_{4n^2_g}\cdot N^2  \simeq \F{\hbar^2}{c^2}\ \omega ^2_{2n_gN}$ is the square of the energy of a quantum (defined in section 2.7 as a $P_2(F^+_{[v_1]})$-subsemimodule of rank $N$~) in a maximal (pseudo-)ramified torus of class $2n_g$ where $\omega _{2n_gN}$ is the angular frequency of this quantum.

Remark that $E_{4(n_g-1)^2N^2}\ge E_{4n^2_gN^2}$~, which means that the energy of a quantum varies with respect to the (global) level in the sense that when the global level $n_g$ increases, the quantum angular frequency $\omega $ tends to decrease (for new arithmetic and algebraic developments of the quanta, see \cite{Pie1}).\epr
}
\vskip 11pt

In the context of the developments of this work, the Birch-Swinnerton-Dyer conjecture can be approached by the following statement:\vskip 11pt

\subsubsection{Proposition}  {\em Let $L\RL(s_{\mp},E(\QQ))\subset L\RL(s_{\mp})$ be the $L$-subseries attached to an 
elliptic curve over $\QQ$~. Then, the  rank of $E(\QQ)$ is the order of vanishing of 
$L^{nr}\RL(s_{\mp},E(\QQ))$ at $s=1$ if $L^{nr} \RL(s_{\mp},E(\QQ))$ is holomorphic at $s=1$~.}\vskip 11pt

\paragraph{Proof:}\mbox{} \Be 
\item Assume that there is a finite set of $P$ primes $p$ for which the mapping $E(\QQ)\to E(\prod\limits_p\FF_p)$ is a bijection. Then, to the $P$ Euler factors of $L^{nr} \RL(s_{\mp})$ correspond the $L$-subseries 
\[L^{nr} \RL(s_{\mp},E(\QQ))=\sum_{n_{g}}\lambda^{nr} _{\mp}(n^2_{g},m_{n_g}^2)\ 
n_{g}^{-s_{\mp}}\qquad \text{with $N_g$ integers ``~$n_g$~''}           \;, \]
in such a way that $N_g$ corresponds to the number of restricted places of the real semifield $F^+\RL$~: the are noted $\o v_g=\{\cdots,\o v_{n_g},\cdots,\o v_{N_g}\}$
(resp. $v_g=\{\cdots, v_{n_g},\cdots,\linebreak v_{N_g}\}$~).\vskip 11pt

 \item To $L^{nr} \RL(s_{\mp},E(\QQ))$ corresponds a restricted simple global elliptic $A_{R_g,L_g}$-subsemi\-module 
$\phi ^{nr,0} \RL(s\RL)_{\rm res}$ which is in one-to-one correspondence with the representation space $\Repsp(T^t_2(F^{+,T,nr}_{\o v_g}))$
(resp. $\Repsp(T_2(F^{+,T,nr}_{ v_g}))$~) of the algebraic semigroup
$T^t_2(F^{+,T,nr}_{\o v_g})$ (resp. $T_2(F^{+,T,nr}_{ v_g})$~) $\subset \GL_2(F^{+,T,nr}_{\o v_g}\times F^{+,T,nr}_{ v_g})$~. At every place $\o v_{n_g}$ (resp. $v_{n_g}$~) associated with the integer 
$n_{g}$~, there is a  one-dimensional irreducible curve $E_f(n_{g})\RL$~, isomorphic to a one-dimensional semitorus $T^1\RL[n_{g}]$ of class $n_g$ according to section 3.1. And, thus, at the $N_g$ places, we have the set $\{E_f(n_{g})\}_{n_g}$ of $N_g$ irreducible curves of class $n_g$~. \vskip 11pt

 \item Then, a natural surjection exists 
\[ S_{E_f\to E(\QQ)} : \quad \{E_f(n_{g})\}_{n_g}\To E(\QQ)\] which identifies the elliptic curve $E(\QQ)$ with the orbit space of $\{ E_f(n_{g})\}_{n_g}$ under the action of the semigroup $T^t_2(F^{+,T,nr}_{\o v_g})$ (resp. $T_2
(F^{+,T,nr}_{v_g})$~) such that $E(\QQ)$ has $N_g$ generating points.

Then, the  rank of $E(\QQ)$ is $r(E(\QQ))=N_g$~.\vskip 11pt

 \item Referring to the paper of N. Katz \cite{Kat}, the $L$-series $L\RL(s_{\mp})$ are holomorphic at the point $s=1$ as a consequence of a functional equation under $s\to 2-s$~. And, thus, if the $L$-series $L^{nr} \RL(s_{\mp},E(\QQ))$ are also holomorphic at $s=1$~, we can speak of the order of vanishing of $L^{nr} \RL(s_{\mp},E(\QQ))$ at $s=1$~.\vskip 11pt

 \item Taking into account the proposition 3.3.5 and section 3.3.2, it appears that the pairs of non-trivial zeros of $L^{nr} \RL(s_{\mp},E(\QQ))$ lie on $\Re(s)=1$ and are in one-to-one correspondence with the trivial zeros. 
  And, thus, to each integer $n_g$ of $L^{nr} \RL(s_{\mp},E(\QQ))$ corresponds a trivial zero and a pair of non-trivial zeros.\vskip 11pt

\item As there are $N_g$ non-trivial zeros of $L^{nr} \RL(s_{\mp},E(\QQ))$~, the order of vanishing of $L^{nr} \RL(s_{\mp},E(\QQ))$ at $s=1$ is $N_g$~.  Now, $N_g$ is also the rank of the elliptic curve $E(\QQ)$ according to 3. \Ee\vskip 11pt

\subsubsection[Connecting the Birch-Swinnerton-Dyer conjecture to the Shimura-Taniyama-Weil conjecture]{Connecting the Birch-Swinnerton-Dyer conjecture to the Shimura-\protect\linebreak Taniyama-Weil conjecture}

The conjecture of Shimura-Taniyama-Weil, dealing with the connection of the zeros of the $L$-series $L^{nr}\RL(s_\mp,E(\QQ))$ with the rank of the elliptic curve $E(\QQ)$~, seems to take into account only the restricted cuspidal biform $f^{nr}_R(z)_{\rm res}\otimes_D f^{nr}_L(z)_{\rm res}$ degenerating into the restricted simple global elliptic bisemimodule $\phi ^{nr,0}_R(s_R)_{\rm res}\otimes_D \phi ^{nr}_L(s_L)_{\rm res}$ (see section 3.3.2) in such a way that only one pair of semicircles $T^1_R[n_g]$ and $T^1_L[n_g]$ arises at each class $n_g$ according to proposition 3.3.6.

The general case would refer to the Shimura-Taniyama-Weil conjecture in such a way that the restricted (pseudo-)unramified cuspidal biform of weight two and level 1 (~$N=1$~)
\[
f^{nr}_R(z)_{\rm res}\otimes_D f^{nr}_L(z)_{\rm res}
=\sum_{n_g}\L( \lambda ^{nr}_-(n^2_g,m^2_{n_g})\ e^{-2\pi in_gz}
\otimes_D\lambda ^{nr}_+(n^2_g,m^2_{n_g})\ e^{2\pi in_gz}\R)\;, \quad z\in\CC\;, \]
be covered by the restricted pseudo-unramified global elliptic bisemimodule:
\[
\phi ^{nr}_R(s_R)_{\rm res}\otimes_D \phi ^{nr}_L(s_L)_{\rm res}
= \sum_{n_g} \sum_{m_{n_g}} \L( \lambda ^{nr}_-(n^2_g,m^2_{n_g})\ e^{-2\pi in_gx}
\otimes_D\lambda ^{nr}_+(n^2_g,m^2_{n_g})\ e^{2\pi in_gx}\R)\;, \quad x\in\RR \]
(see proposition 3.1.10), in such a way that the modular representation of the elliptic curve $E(\QQ)$ be given by $n_g$ sets of $m_{n_g}$ products, right by left, of semicircles 
$T^1_R[n_g,m_{n_g}]$ and $T^1_L[n_g,m_{n_g}]$~.  Then, the rank $r_{E(\QQ)}$ of $E(\QQ)$ would be given by
\[r_{E(\QQ)}=\bigoplus_{n_g} m_{n_g}\;\quad m_{n_g}\ge 1\;.\]

\end{document}